\documentclass[final,hidelinks,onefignum,onetabnum]{siamart250211}

\usepackage{filecontents}

\usepackage[nameinlink]{cleveref}
\usepackage{array}
\usepackage[caption=false]{subfig}
\usepackage{booktabs}
\usepackage{stmaryrd}
\usepackage{lipsum}
\usepackage{amsfonts}
\usepackage{graphicx}
\usepackage{epstopdf}
\usepackage{algorithmic}
\ifpdf
  \DeclareGraphicsExtensions{.eps,.pdf,.png,.jpg}
\else
  \DeclareGraphicsExtensions{.eps}
\fi

\newsiamremark{example}{Example}
\newsiamremark{remark}{Remark}

\headers{2D-FC for Domains with Corners}{O. Bruno, A. Yang}

\title{Two Dimensional Fourier Continuation\\ for Domains with Corners\thanks{Submitted to Editors 03/17/2026
\funding{The authors gratefully acknowledge support from the Air Force Office
of Scientific Research and the National Science Foundation under
contracts FA9550-21-1-0373, FA9550-25-1-0015, DMS-2109831, and the Citadel Global Fixed Income SURF Fellowship}}}

\author{Oscar P. Bruno\thanks{Computing and Mathematical Sciences, Caltech
  (\email{obruno@caltech.edu}).}
\and Allen Yang\thanks{School of Engineering and Applied Sciences, Harvard University (\email{ayang@fas.harvard.edu}).}}

\usepackage{amsopn}

\newcommand\e{\textrm{e}}
\newcommand\m{\textrm{m}}
\newcommand\bt{\textrm{b}}

\newcommand{\Intrl}[2]{\mathrm{int}_{#1}\!\left(#2\right)}
\newcommand{\Clrl}[2]{\mathrm{cl}_{#1}\!\left(#2\right)}
\newcommand{\Int}[1]{\mathrm{int}\!\left(#1\right)}
\newcommand{\Cl}[1]{\mathrm{cl}\!\left(#1\right)}

\makeatletter
\newcommand*{\addFileDependency}[1]{
  \typeout{(#1)}
  \@addtofilelist{#1}
  \IfFileExists{#1}{}{\typeout{No file #1.}}
}

\DeclareFontFamily{U}{tipa}{}
\DeclareFontShape{U}{tipa}{m}{n}{<->tipa10}{}
\newcommand{\arc@char}{{\usefont{U}{tipa}{m}{n}\symbol{62}}}%

\newcommand{\arc}[1]{\mathpalette\arc@arc{#1}}

\newcommand{\arc@arc}[2]{%
  \sbox0{$\m@th#1#2$}%
  \vbox{
    \hbox{\resizebox{\wd0}{\height}{\arc@char}}
    \nointerlineskip
    \box0
  }%
}
\makeatother

\usepackage{cite}
\ifpdf
\hypersetup{
  pdftitle={Two-Dimensional Fourier Continuation for Domains with Corners},
  pdfauthor={O. Bruno, A. Yang}
}
\fi

\begin{document}
\maketitle

\begin{abstract}
  This paper presents a fast ``two-dimensional Fourier Continuation''
  (2D-FC) method for the construction of biperiodic extensions of
  smooth, non-periodic functions defined over general two-dimensional
  (2D) domains, including domains with corners.  The algorithm
  operates with an $O(N \log N)$ computational cost, for an
  $N$-point discretization grid, and it achieves a user-prescribed
  $d$-th order of accuracy. The methodology can be generalized to
  non-smooth domains of arbitrary dimensionality, but such extensions
  are not considered in the present work. The usefulness and
  performance of the 2D-FC method are demonstrated through
  applications to the Poisson problem posed on bounded 2D domains
  with corners. One illustrative application concerns a Poisson
  problem on a non-smooth drop-shaped domain with a highly oscillatory
  forcing term; employing a discretization containing approximately
  four million degrees of freedom, the method produces the Poisson
  solution with accuracies of the order of machine precision in
  computing times on the order of one second on a single core of a
  present-day computer.
\end{abstract}


\section{Introduction}\label{sec:introduction}
This paper concerns the problem of the extension of non-periodic functions
$f$, defined on a given $n$-dimensional domain $\Omega$, to periodic
functions on a rectangular domain containing $\Omega$---with
particular attention to cases in which $\Omega$ has a non-smooth
boundary.  For definiteness only the 2D case ($n=2$) is considered
here, but the proposed methodology can be generalized to non-smooth
domains of arbitrary dimensionality. This paper addresses the problem
by suitably generalizing the smooth-domain extension-along-normals 2D
Fourier Continuation methodology introduced in
reference~\cite{bruno_two-dimensional_2022}. The generalization incorporates novel techniques to enable application to non-smooth domains, including a procedure for smoothly decomposing near-boundary function values using a near-boundary partition of unity. Crucially, it
also incorporates ``corner-patch'' parametrizations and discretization
strategies for convex and concave corners. These strategies enable the
robust extension of function values around near-corner boundary
regions---where the normal vector is not well defined and where ad-hoc
generalizations of the extension-along-normals concept used
in~\cite{bruno_two-dimensional_2022} may become problematic.

The new approach is highly effective: at a cost of $O(N\log N)$
operations (for an $N$-point discretization grid) it achieves
accuracies near machine precision in the approximation of functions on
domains containing multiple corners, including extreme convex and
concave cases, and at a user-prescribed $d$-th order of accuracy. For
example, near–machine–precision and high-order accuracy are reported
in this work for applications to corner problems with acute interior
angles of $1.8^\circ$, as well as corresponding obtuse interior angles
of $360^\circ - 1.8^\circ = 358.2^\circ$.

The potential usefulness and performance of the 2D-FC method are
illustrated through applications to the Poisson equation posed on
bounded two-dimensional domains with corners. In particular, a Poisson
problem is considered, on a non-smooth drop-shaped domain and with a
highly oscillatory forcing term. Employing discretizations containing
nearly four million degrees of freedom, the 2D-FC-based algorithm
yields near-machine-precision solution accuracies in computing times
on the order of one second on a single core of a present-day
computer. This compares favorably with the recent smooth-domain (but
applicable to sharp smooth boundaries) Poisson
solver~\cite{fortunato2025fully}.  Additional potentially impactful
applications of the proposed methodology may be envisioned, including,
possibly, as an enabling element in algorithms for the evaluation of singular
integral operators (cf.~\cite{bruno_cao}).

The problem of constructing periodic extensions has received
considerable attention in the recent literature, particularly due to its
relevance in the solution of partial differential
equations (PDEs); see
e.g.~\cite{bruno_two-dimensional_2022,amlani_fc-based_2016,ALBIN20116248}
and references therein. Contributions such as~\cite{amlani_fc-based_2016,ALBIN20116248} employ the
one-dimensional Fourier Continuation (FC) approach together with the
method of lines to tackle multidimensional PDEs in general domains for
a wide range of physical applications, including possibly high
frequencies, complex spatio-temporal variations, and shock waves. In
all cases the method provides solutions with low dispersion and high
accuracy, and it lends itself to highly efficient implementation in
high-performance computing systems.

Existing multidimensional extension approaches---including the method
of~\cite{bruno_two-dimensional_2022}, which the present work
generalizes, as well as RBF-based algorithms~\cite{FRYKLUND201857} and
PDE-based extension
approaches~\cite{stein2016immersed,askham2017adaptive}---typically
rely on boundary-smoothness assumptions that fail at corners; additional challenges posed by some of these methods are discussed
in~\cite{fortunato2025fully}. The 2D smooth-domain extension
method~\cite{matthysen_huybrechs_fa_frames}, in turn, has
remained largely of theoretical interest due to its significant
computational cost~\cite{bruno_two-dimensional_2022}. The near-boundary
partition-of-unity strategy and the specialized corner patches
introduced in the present contribution are particularly well suited
for use in conjunction with the extension-along-normals
strategy~\cite{bruno_two-dimensional_2022}, but they may also prove useful in combination with other two-dimensional extension methodologies. 

This paper is organized as follows:
Sec.~\ref{sec:fc_1d} reviews a certain one-dimensional
blending-to-zero (BTZ) procedure that underlies both the 1D-FC and the proposed 2D-FC methods; Sec.~\ref{sec:2dfc-prelim} outlines the overall 2D-FC framework for
domains with corners; and Sec.~\ref{sec:2d-fc} details the core
components of the corner-capable 2D-FC algorithm, including the
strategies employed for the construction of the parameterized
boundary-fitted patches, matching domains, partition of unity, and
extension domains, and the approaches used for interpolation onto the
Cartesian mesh. A variety of numerical results are presented in
Sec.~\ref{sec:numerical-results}, including, in Sec.~\ref{2dfcex},
numerical application of the 2D-FC approximation across a range of
challenging domains. The 2D-FC–based Poisson solver is described and
demonstrated numerically in Sec.~\ref{sec:poisson-problem}.
Conclusions, finally, are provided in Sec.~\ref{sec:conclusion}. The
supplementary material section presents details concerning both 
the geometry-processing algorithms employed and the method used for the near-boundary evaluation of the boundary integral operators utilized in the solution of the Poisson problem. C and MATLAB implementations of the algorithms presented in this work are publicly available at~\cite{yang_2dfc_c} and~\cite{yang_2dfc_matlab}.

\section{Background: 1D ``blending-to-zero'' FC algorithm}\label{sec:fc_1d}
Following~\cite[Sec. 2]{bruno_two-dimensional_2022}, this section presents a brief review of the 1D FC-Gram Fourier Continuation (1D-FC) method~\cite{amlani_fc-based_2016,ALBIN20116248,BRUNO20102009}, with an emphasis on one of its key components, namely, the {\em blending-to-zero procedure} (BTZ), which is also employed as an important element of the proposed 2D-FC method introduced in the present contribution. Throughout this paper the term ``smooth'' will be understood to mean ``infinitely differentiable''.

For a given smooth function $\phi$ defined in a 1D interval, say, for definiteness,
\begin{equation}
    \phi: [0, 1]\rightarrow\mathbb{C},
\end{equation} 
the BTZ algorithm produces a smooth extension of the function $\phi$, e.g. toward the right and over an extension interval $[1,b]$ with $b>1$, which equals $\phi$ on $[0,1]$ and which smoothly matches the function equal to zero for $x\geq b$. The resulting ``right blended to zero'' extended function $\phi^{rb}$  is produced by first obtaining, on the basis of a certain number $d$ of values of $\phi(x)$ with $x$ near 1 ($x\leq 1$),  a number $C$ of values of $\phi^{rb}(x)$ for $1< x\leq b$ which smoothly transition $\phi$ from its values near $x=1$ to the function which equals $0$ for all $x\geq b$. A leftward blending can be defined analogously: $C$ discrete values of some function $\phi^{\ell b}$, which extends $\phi$ to the left over the interval $[a,0]$ for some $a< 0$, are obtained on the basis of a certain number $d$ of values of $\phi(x)$ with $x$ near 0 ($x\geq 0$) such that they smoothly match $\phi$ to the function which equals zero for $x\leq a$. The left- and right-blending algorithms, which are described in what follows, are used in the 1D-FC method~\cite{amlani_fc-based_2016,ALBIN20116248} together with the FFT algorithm to produce 1D Fourier continuation function that converges to $\phi$ with order $d$ for $0\leq x\leq 1$~\cite{ALBIN20116248,amlani_fc-based_2016}. Like the 1D-FC algorithms, the algorithm presented in this paper utilizes the aforementioned $C$ discrete blending values along certain sets of 1D lines to produce 2D-FC extensions for functions defined in 2D domains with corners.

While the 2D-FC method proposed in this paper makes exclusive use of the $\ell$-BTZ algorithm, we nonetheless briefly review the historically relevant rightward BTZ algorithm ($r$-BTZ): the corresponding leftward version is obtained by a reduction to the rightward case (see Rem.~\ref{rmk:l-btz}). The $r$-BTZ algorithm produces the $C$ values of the extension function $\phi^{rb}$ on the basis of the vector $\Phi = (\Phi_0, \dots, \Phi_{N - 1})^t$ containing $N$ discrete values $\Phi_j =\phi(x_j)$ of $\phi$ on the 1D Cartesian mesh $\{x_0,x_1,\dots,x_{N-1}\}\subset [0,1]$, with $x_0=0$ and $x_{N-1} = 1$ and with step size $k = 1/(N-1)$. In fact, only the $d$ rightmost values $\Phi^r = (\Phi_{N-d}, \dots, \Phi_{N-1})^t$ of the vector $\Phi$ are used to obtain the desired $C$ right-extension values---which are then used to create the ($N+C$)-dimensional vector $\Phi^{rb} = (\Phi^{rb}_0, \dots, \Phi^{rb}_{N+C-1})^t$  satisfying $\Phi^{rb}_j = \Phi_j$ for $0 \leq j \leq N-1$ which  smoothly ``blends'' $\Phi$ to zero towards the right. The vector 
\begin{equation}\label{eq:ext_vec}
    \Phi^{rb}_C = (\Phi^{rb}_N, \dots, \Phi^{rb}_{N+C-1})^t
\end{equation}
of extension values is produced by means of a certain least-squares procedure that is outlined in what follows. 

The 1D $r$-BTZ algorithm precomputes $C$ discrete blending-to-zero values in the rightward ``extension interval'' 
\begin{equation}\label{eq:Iextr}
    I^{\e, r} = (x_{N-1}, x_{N+C-1}]
\end{equation}
for each polynomial $g_j\in G_d$ where $G_d = \{g_0(x),g_1(x),\dots,g_{d-1}(x) \}$ denotes the orthogonal Gram polynomial basis---that results from an application of the Gram-Schmidt orthogonalization method to the standard basis $\{ 1, x,x^2,\dots, x^{d-1}\}$ of the space of polynomials of degree less than or equal to $(d-1)$, in the order of increasing degree, and with respect to the scalar product
\begin{equation}
    (g, h) = \sum_{j=N-d}^{N-1} g(x_j)h(x_j).
\end{equation}

The vectors of values of the Gram polynomials on the rightward ``matching mesh'' 
\begin{equation}\label{eq:Dmatchr}
    D^{\m, r} = \{x_{N-d}, ..., x_{N-1}\}
\end{equation}
are the columns of the $QR$ factorization $\boldsymbol{P} = \boldsymbol{QR},$
 of the Vandermonde matrix $\boldsymbol{P}$ corresponding to the point set $D^{\m, r}$~\cite{amlani_fc-based_2016, golub2012matrix}: the $ij$-th element of $\boldsymbol{Q}$ contains
the value of the $j$-th Gram polynomial at the $i$-th grid point of the discretization. For added accuracy, in the FC applications the necessary QR factorizations are usually obtained by means of the stabilized Gram-Schmidt procedure, a practice that is also followed in the present paper. Note that, calling $\phi_p$ the polynomial of degree $d-1$ that interpolates a given function $\phi$ on $D^{\m, r}$ the $j$-th coefficient in the expansion of $\phi_p$ in the basis $G_d$ is given by $(g_j, \phi_p) = (\boldsymbol{Q}^T \Phi^r)_j$.
The $C$ $r$-BTZ values in the extension interval $I^{\e, r}$  for a given function $\phi$ are produced on the basis of corresponding {\em precomputed} blending-to-zero values for each Gram polynomial $g_j$. Thus, letting $A\in\mathbb{C}^{C\times d}$ denote the ``extension matrix'' whose $j$-th column contains the precomputed blending-to-zero values for $g_j$, the vector~\eqref{eq:ext_vec} of extension values is given by
\begin{equation}\label{eq:b_to_z}
  \Phi^{rb}_C =  \boldsymbol{A}\boldsymbol{Q}^{T} \boldsymbol{\Phi}^r.
\end{equation}
To complete the algorithm it therefore suffices to specify the procedure that yields the $g_j$ blending-to-zero values that make up the columns of $A$.

The $r$-BTZ algorithm obtains the necessary BTZ values for each Gram polynomial $g_j$ by constructing a band-limited trigonometric polynomial of the form
\begin{equation}\label{eq:auxi_trig_poly}
  g^c_j(x) = \sum_{m = -J} ^ {J} a_m^j e ^ {\frac{2\pi i m (x-x_{N-d})}{(d + 2C + Z - 1)k}},
\end{equation}
with $a^j_m$ chosen so that $g^c_j$ (a)~Approximates $g_j$ on the ``rightward-matching interval'' 
\begin{equation}\label{eq:Imatchr}
    I^{\m, r} = [x_{N-d}, x_{N-1}] \supset D^{\m, r};
\end{equation}
and, (b)~Approximates the zero function at a number  $Z$ of equispaced ``zero-matching points'' in the ``zero-matching interval'' $[x_{N+C-1}, x_{N+C-1}+(Z-1)k]$, in the least squares sense. To achieve higher accuracy approximations, the $r$-BTZ method utilizes sets of matching points, at both the right matching interval and the zero-matching interval, whose step size is finer than $k$, with step size $k/n_\mathrm{os}$, for a suitably selected integer oversampling factor $n_\mathrm{os} > 1$; the value $n_\mathrm{os} = 20$ is used throughout this paper. The values of the Gram polynomials on the finer grid are given by $\boldsymbol{Q}_{\textrm{os}} = \boldsymbol{P}_{\textrm{os}} \boldsymbol{R} ^ { - 1},$
where $\boldsymbol{P}_\textrm{os}$ is a $(n_\textrm{os}(d-1)+1)\times d$ Vandermonde matrix corresponding to the oversampled grid.
The coefficients $a^j_m$ in~\eqref{eq:auxi_trig_poly} are then selected to equal the entries of the minimizing vector $\boldsymbol{a}^j$ of the problem
\begin{align}\label{eq:min_prob}
\boldsymbol{a}^j = {(a^j_{ - J}, \dots, a^j_J)} ^ T = \textrm{argmin}
  \left\|
  \boldsymbol{B} _ {\textrm{os}}\boldsymbol{a} -
  \begin{pmatrix}
    \boldsymbol{q} ^ j_{\textrm{os}}\\
    \boldsymbol{0}
  \end{pmatrix}
  \right\|_{2}, 
\end{align}
where $\boldsymbol{q} ^ j_{\textrm{os}}$ denotes the $j$-th column of $\boldsymbol{Q}_{\textrm{os}}$, where $\boldsymbol{0}$ denotes the $(n_\textrm{os}(d-1)+1)$-dimensional column zero vector, and where the $m$-th column of the matrix $\boldsymbol{B}_\mathrm{os}$  contains the values of the $m$-th exponential factor in the sum~\eqref{eq:auxi_trig_poly} evaluated at each point in the union of right-matching and zero-matching point sets, in increasing order. The minimizing coefficients are computed by means of the SVD-based least-squares approach~\cite{golub2012matrix}. Once the
 vectors $\boldsymbol{a}^j$ of coefficients ($1 \leq j\leq d-1$) have been obtained, the values of the resulting Fourier expansion~\eqref{eq:auxi_trig_poly} on the rightward ``extension mesh''
 \begin{equation}\label{eq:Dext}
     D^{\e, r} = \{x_N, \dots, x_{N+C-1}\} \subset I^{\e, r}
 \end{equation}
are computed; these are the needed BTZ values that make up the $j$-th column of $\boldsymbol{A}$. Using this matrix, the BTZ values $\Phi^{rb}_C$ for any given vector $\Phi^r$ are given by~\eqref{eq:b_to_z}.

Since,  for a given order $d$, the BTZ values generated by these matrices do not depend on the point spacing $k$, the matrices $\boldsymbol{A}$ and $\boldsymbol{Q}$ can be precomputed and stored on disk to be input as needed in any application that requires the BTZ process.
\begin{remark}
    \label{rmk:l-btz}
    The \textit{leftward} blending-to-zero algorithm ($\ell$-BTZ) is easily reduced to the rightward problem, as described in what follows. Similar to the $r$-BTZ method, the $\ell$-BTZ algorithm uses the $d$ leftmost values $\Phi^{\ell} = (\Phi_0, \dots, \Phi_{d-1})$ to produce a ($N+C$)-dimensional vector $\Phi^{\ell b} = (\Phi_{-C}, \dots, \Phi_{N-1})^t$ which satisfies $\Phi^{\ell b}_j = \Phi_j$ for $0 \leq j\leq N-1$ and blends $\phi$ to zero towards the left. To compute the $C$ extension values $\Phi^{\ell b}_C = (\Phi^{\ell b}_{-C}, \dots, \Phi^{\ell b}_{-1})^t$, we introduce the ``order reversion'' matrix $R^e \in \mathbb{C}^{e\times e}$ with $e \in \mathbb{N}$ such that 
    \[R^e(g_0, g_1, \dots, g_{e-2}, g_{e-1})^t = (g_{e-1}, g_{e-2}, \dots, g_1, g_0).\]
    Then, 
    \begin{equation}\label{eq:reversion}
    \Phi_C^{\ell b} = R^C\boldsymbol{A}\boldsymbol{Q}^TR^d\Phi^\ell.
    \end{equation}
    Since the 2D-FC method in this paper uses only the leftward $\ell$-BTZ procedure, we simplify notation by omitting the $\ell$ superscript in the following definitions: we denote the leftward matching interval and matching mesh
    \begin{equation}\label{eq:Im}
        I^{\m} = I^{\m, \ell} = [x_0, x_{d-1}] \quad \textrm{and}\quad D^{\m} = D^{\m, \ell} = \{x_0, x_1, \dots, x_{d-1}\} \subset I^{\m};
    \end{equation}
    and, similarly, we denote the leftward extension interval and extension mesh by
    \begin{equation}\label{eq:Ie}
        I^{\e} = I^{\e, \ell} = [x_{-C}, x_0] \quad \textrm{and}\quad D^{\e} = D^{\e, \ell} = \{x_{-C}, x_{-C+1}, \dots, x_{-1}\} \subset I^{\e}.
    \end{equation}
\end{remark}
\begin{remark}\label{rmk:refinebtz}
  In the course of this work, and in order to prevent accuracy loss in the various interpolation steps when high precision (say, with errors smaller than $10^{-6}$) is required, it was necessary to obtain BTZ values in the interval $(x_{N-1}, x_{N+C-1}]$ (resp. $[x_{-C}, x_0)$) on an equispaced grid finer than the grid $\{x_{N}, x_{N+1} \dots, x_{N+C-1}\}$ (resp. $\{x_{-C}, x_{-C+1}, \dots, x_{-1}\}$); cf.~\cite[Rem. 2.2]{bruno_two-dimensional_2022} where a similar comment can be found. In the $r$-BTZ case, for example, once the minimizer $\boldsymbol{a}^j$ in \eqref{eq:min_prob} has been obtained, the fine-grid values can be obtained by computing the values of the function $g^c_j(x)$ in~\eqref{eq:auxi_trig_poly}  on a desired finer equispaced mesh containing a number $C_\rho$ of fine-grid points (say $C_\rho= \rho C$, with integer refinement factor $\rho> 0$). Using such values, a fine-grid continuation matrix $\boldsymbol{A}_\rho\in \mathbb{C}^{C_\rho\times d}$ is constructed and the necessary BTZ function values on the $C_\rho$ fine-grid points
\begin{equation}
    D^{\e, r}_{C_\rho} = \{\widetilde{x}_{N}, \widetilde{x}_{N + 1}, \dots, \widetilde{x}_{N+C_\rho-1}\},\quad \textrm{where}\quad \widetilde{x}_{N+j} = x_N + \frac{k}{\rho}j,
\end{equation}
are given by $\Phi_{C_\rho}^{{rb}} = \boldsymbol{A}_\rho \boldsymbol{Q}^T\boldsymbol{\phi}_r$. In the $\ell$-BTZ case, the corresponding fine-grid values $\Phi_{C_\rho}^{\ell b}$ on the fine-grid points
  \begin{equation}\label{Derho}
      D^{\e}_{C_\rho} = D^{\e, \ell}_{C_\rho} = \{\widetilde{x}_{-\rho C}, \widetilde{x}_{-\rho C + 1}, \dots, \widetilde{x}_{-1}\},\quad\textrm{where}\quad\widetilde{x}_{-\rho C+j} = x_{-C}+\frac{k}{\rho} j,
  \end{equation}
  are obtained by applying the $r$-BTZ procedure in conjunction with the relation~\eqref{eq:reversion}; that is, $\Phi_{C_\rho}^{\ell b} = R^{C_\rho}\boldsymbol{A}_\rho\boldsymbol{Q}^TR^d\Phi_\ell$. For later use, we additionally include the boundary matching-mesh point $\widetilde{x}_{0} = x_{0}$  in the fine extension mesh~\eqref{Derho}. We furthermore use the notation $\Phi^{\ell b}_{C_{\rho}}$ to denote both the vector of BTZ extension values and the corresponding grid function $\Phi_{C_\rho}^{\ell b}: D^{\mathrm{e}}_{C_{\rho}} \to \mathbb{C}$ that assigns to each point of the extension mesh its associated blended-to-zero value.

\end{remark}
\begin{figure}[htbp]
\vspace{-0.2cm}
\centering
\includegraphics[width=0.9\linewidth,]{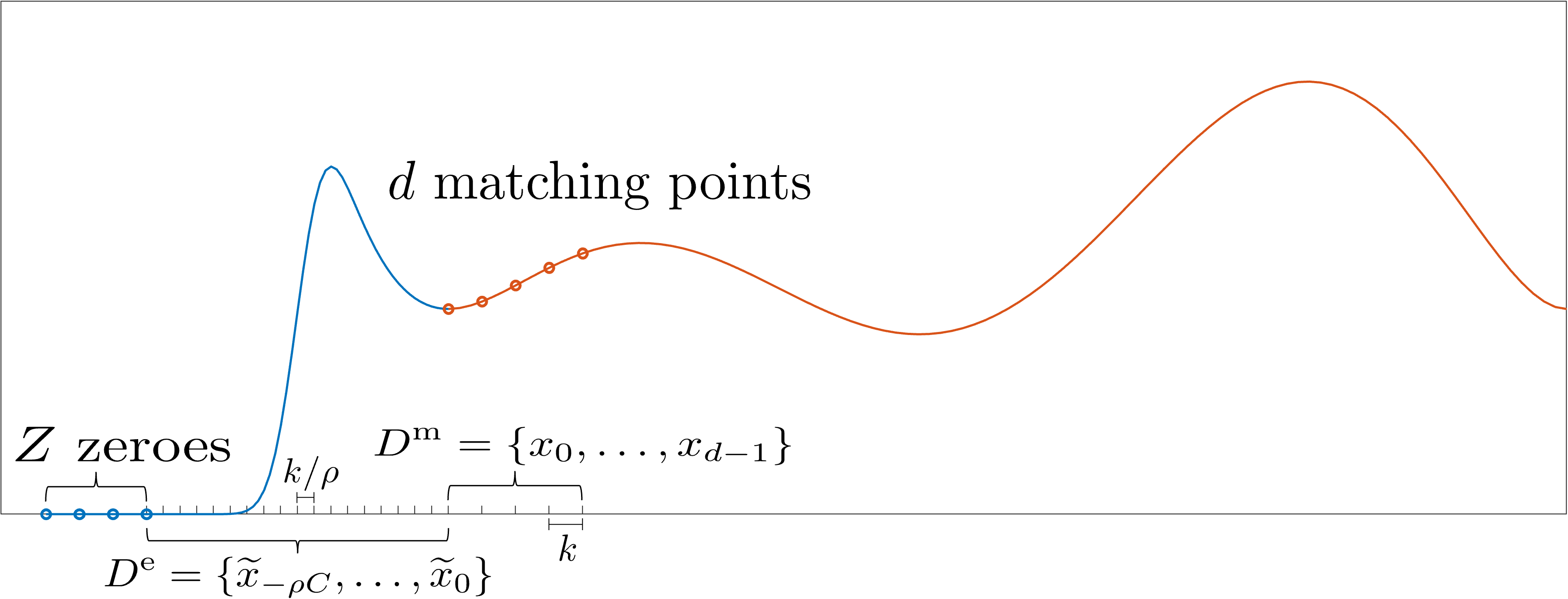}
\vspace{-0.3cm}
\caption{Illustration of the $\ell$-BTZ procedure. The orange line represents the original function, while the blue line shows the 1D $\ell$-BTZ extension function. By construction, the matching interval is $I^\m = [x_0, x_{d-1}]$ and the extension interval is $I^\e = [\widetilde{x}_{-\rho C}, \widetilde{x}_0]$.}
    \label{fig:lBTZ}
    \vspace{-0.5cm}
\end{figure}

\section{2D-FC method on domains with corners: preliminary description}\label{sec:2dfc-prelim}
This section presents an outline of the proposed FC method for 2D domains $\Omega \subset \mathbb{R}^2$ with a piecewise smooth boundary $\Gamma = \partial \Omega$ that may contain convex and/or concave corners (that is, corners whose angles interior to $\Omega$ are less than $180^\circ$ and larger than $180^\circ$, respectively); a detailed description then follows in Sec.~\ref{sec:2d-fc}. 

Like~\cite{bruno_two-dimensional_2022}, the proposed algorithm constructs a two-dimensional Fourier continuation for a smooth function  
\begin{equation}\label{eq:f}
    f:\overline{\Omega} \rightarrow \mathbb{C},
\end{equation}  
defined in the closure $\overline{\Omega}$ of the domain $\Omega$, by extending $f$ outwards from a region near the domain boundary and smoothly blending the extension to zero at a certain distance from the  boundary. Following~\cite{bruno_two-dimensional_2022}, for regions adjacent to smooth portions of the boundary $\Gamma$ and away from corner points the new algorithm generates the required extensions by employing the 1D-BTZ algorithm reviewed in Sec.~\ref{sec:fc_1d} along lines normal to $\Gamma$.  Clearly, a different strategy is required at and near corner points, 
since the normal vectors inevitably lead to intersections of normal lines around concave corners, and to regions unswept by normals around convex corners. 

\begin{figure}[htbp]
\vspace{-0.1cm}
\centering
\includegraphics[width=0.7\linewidth,]{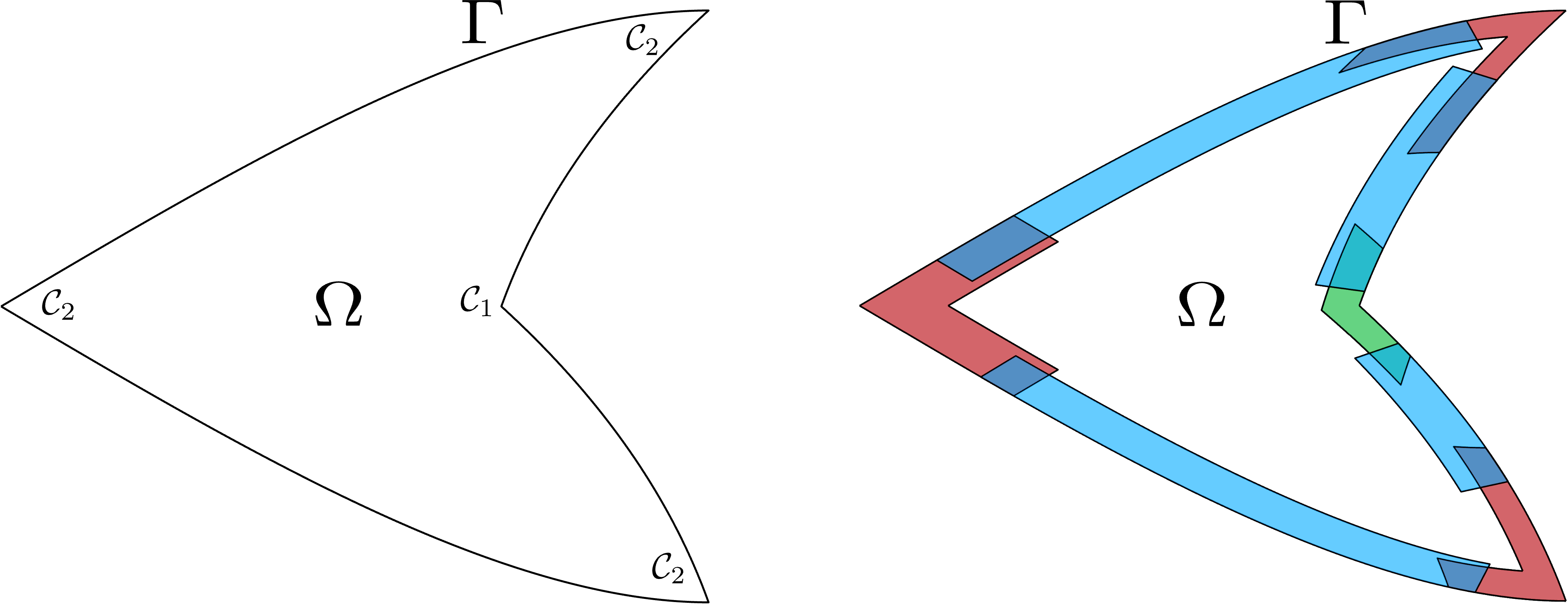}
\vspace{-0.25cm}
\caption{Left: 2D domain $\Omega$ with a piecewise smooth boundary $\Gamma$ that has one $\mathcal{C}_1$-type corner and three $\mathcal{C}_2$-type corners. Right: patches used in an overlapping-patch decomposition of a neighborhood of $\Gamma$ contained in $\Omega$, including, (a)~Smooth-boundary $\mathcal{S}$-type patches (blue); (b)~A concave corner $\mathcal{C}_1$-type patch (green); and, (c)~Three convex corner $\mathcal{C}_2$-type patches (red).}
    \label{fig:patchdecompprelim}
    \vspace{-0.2cm}
\end{figure}

\begin{figure}[htbp]
\vspace{-0.5cm}
\centering
\includegraphics[width=0.4\linewidth,]{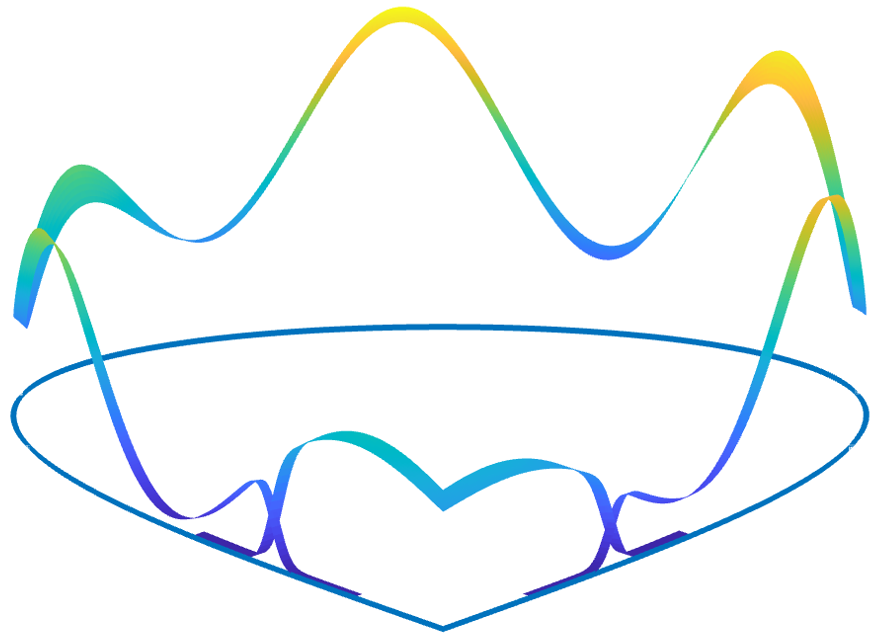}
\vspace{-0.2cm}
\caption{Example of patch-localized windowed functions for a decomposition consisting of one $\mathcal{S}$-type patch and one $\mathcal{C}_2$-type patch. Note that the properties of the partition of unity ensure that the sum of all the windowed functions coincides with the original function~\eqref{eq:f}.}
    \label{fig:POUprelim}
    \vspace{-0.3cm}
\end{figure}

\begin{figure}[htbp]
\vspace{-0.3cm}
\centering
\includegraphics[width=0.6\linewidth,]{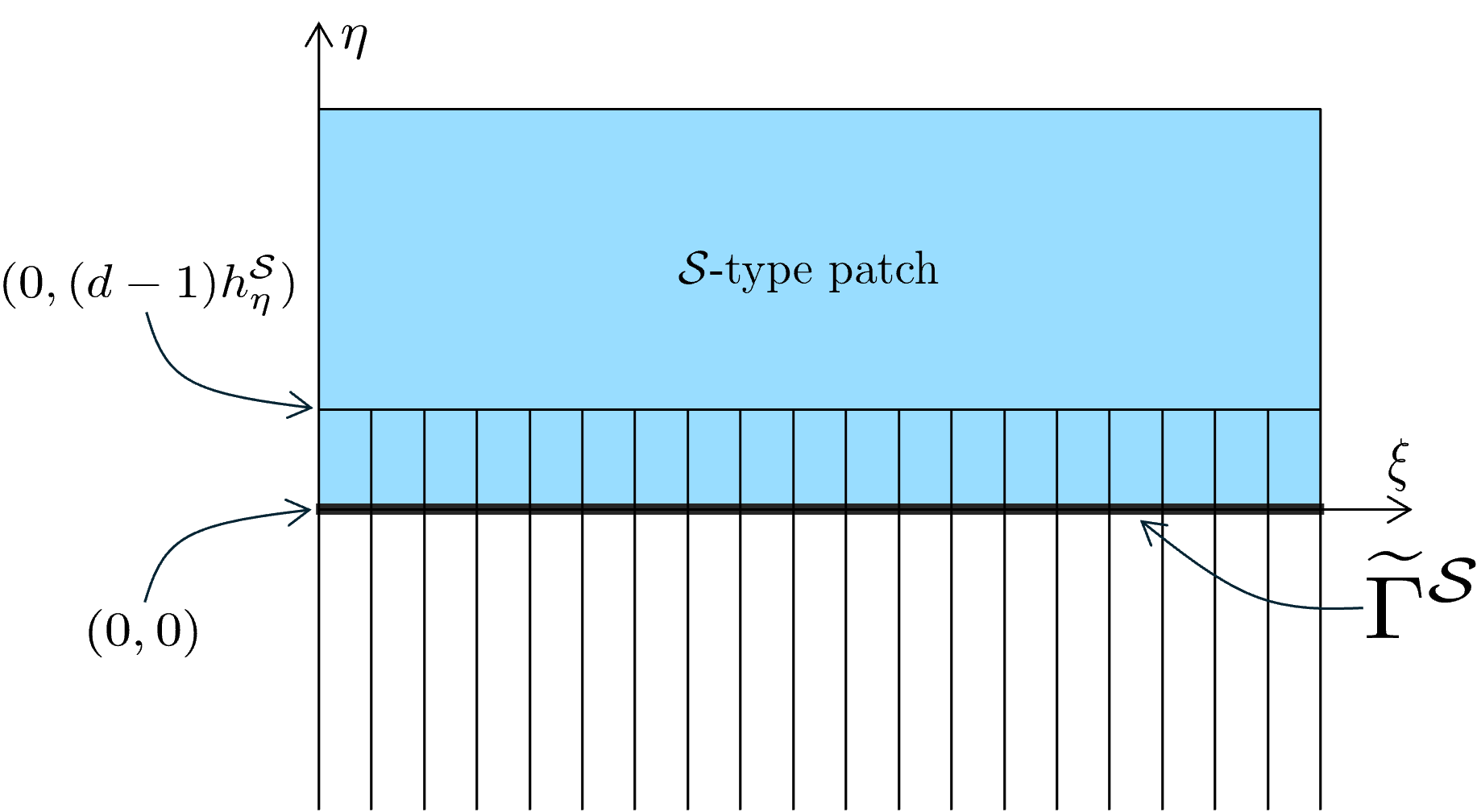}
\vspace{-0.3cm}
\caption{Following~\cite{bruno_two-dimensional_2022}, the proposed algorithm obtains the BTZ extension values on $\mathcal{S}$-type patches by applying the 1D-BTZ procedure along lines normal to the boundary; cf. Fig.~\ref{fig:SBTZdetail}.}
    \label{fig:SBTZprelim}
    \vspace{-0.3cm}
\end{figure}

\begin{figure}[htbp]
\vspace{-0.1cm}
\centering
\includegraphics[width=\linewidth,]{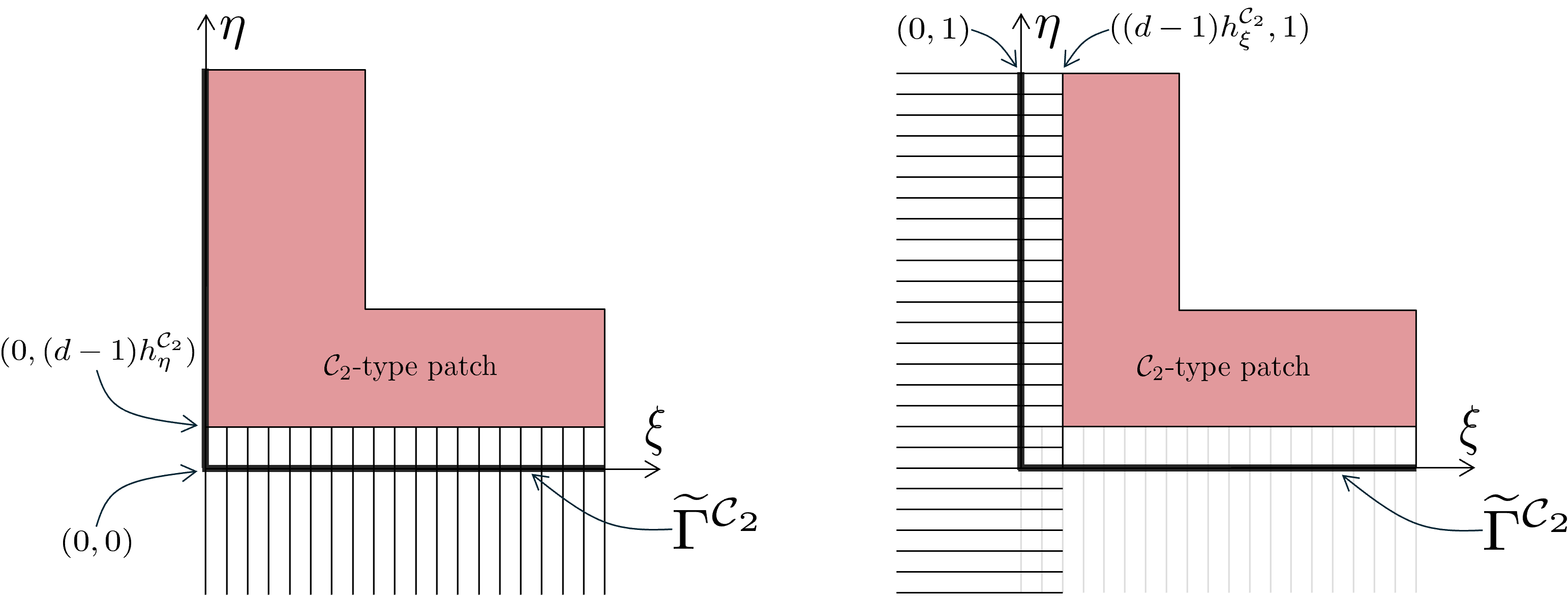}
\vspace{-0.4cm}
\caption{For $\mathcal{C}_2$-type patches the algorithm proceeds in two stages. First, 1D-BTZ is applied along the lines shown in the left panel of the figure. Then, a second 1D-BTZ is performed along the lines shown in the right panel, where both the original patchwise function values and the extension values from the first stage are used in the matching meshes for this second step; cf. Fig.~\ref{fig:C2BTZdetail}.}
    \label{fig:C2BTZprelim}
    \vspace{-0.1cm}
\end{figure}

\begin{figure}[htbp]
\centering
\includegraphics[width=\linewidth,]{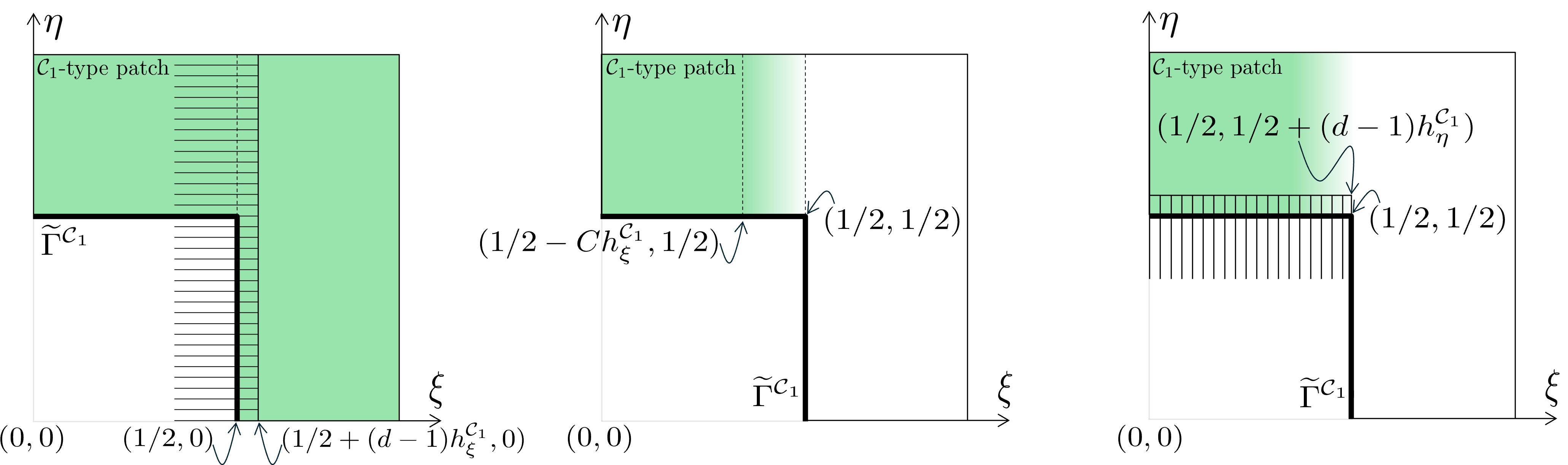}
\vspace{-0.3cm}
\caption{For $\mathcal{C}_1$-type patches, the algorithm proceeds in three steps. Step~1: 1D-BTZ is applied along the discrete set of horizontal lines shown in the leftmost panel, using only the patchwise function values on the subregion consisting of points to the right of the dashed line. Step~2: The intermediate extension obtained in Step~1 is subtracted from the original function, resulting in a modified function that vanishes within this subregion (middle panel, where the white color represents zero function values). Step~3: 1D-BTZ is applied along the set of lines shown in the rightmost panel. Since the subregion is now zero, the resulting extension from this step blends smoothly to zero as it approaches the subregion, and this extension is then added to the one from the first step to yield the final smooth extension; cf. Fig.~\ref{fig:C1BTZdetail}.
}
    \label{fig:C1BTZprelim}
    \vspace{-0.5cm}
\end{figure}

To address these difficulties, the new algorithm first decomposes a neighborhood of the boundary $\Gamma$ into a collection of {\em boundary-fitted patches}, including $\mathcal{S}$-type patches along smooth portions of the boundary; $\mathcal{C}_1$-type patches around concave corners; and $\mathcal{C}_2$-type patches around convex corners. An example of such a decomposition is displayed in Fig.~\ref{fig:patchdecompprelim}; the corresponding parametrization domains are displayed, respectively, in light blue, green, and red colors in Fig.~\ref{fig:SBTZprelim}, and on the left panels of Figs.~\ref{fig:C2BTZprelim} and~\ref{fig:C1BTZprelim}. On each patch, suitably windowed restrictions of $f$ near the boundary are locally extended by employing the corresponding patch parametrizations, as outlined later in this section, and they are subsequently combined to form a single smooth global extension.
To produce patchwise extension functions that can then be combined into an overall extension function, the algorithm employs a partition of unity (POU) subordinated to a certain open covering $\widehat{\Omega}_p^{\mathcal{R}, \m}$ ($\mathcal{R} = \mathcal{S}$, $\mathcal{C}_1$, $\mathcal{C}_2$) of a neighborhood $W^\m$  of the boundary $\Gamma$ within $\overline{\Omega}$, as detailed in Section~\ref{subsec:POU}. Each one of the smooth windowing functions is at first applied, by multiplication, to the restriction of the original function $f$ to $W^\m$---thus producing local components of the near-boundary values of $f$, one for each patch, and each one vanishing outside of the corresponding set $\widehat{\Omega}_p^{\mathcal{R}, \m}$---as illustrated in Fig.~\ref{fig:POUprelim}. Clearly, the sum of these local components equals $f$ throughout $W^\m$---thus guaranteeing that, when the BTZ extensions are added, they naturally assemble into a single globally smooth BTZ extension function.

Three different patchwise 2D-BTZ procedures are employed by the algorithm to produce 2D BTZ extensions---one for each patch type---as illustrated in Figs.~\ref{fig:SBTZprelim}--\ref{fig:C1BTZprelim}. The procedure for $\mathcal{S}$-type patches is straightforward: the BTZ extension values are obtained by applying the 1D-BTZ algorithm employing $d$ points on the $h_\eta^\mathcal{S}$-mesh along lines normal to the boundary in parameter space, as illustrated in Fig.~\ref{fig:SBTZprelim}, and then lifting the resulting values to physical space and interpolating onto the physical Cartesian grid. (The notations $h_\xi^{\mathcal{R}}$ and $h_\eta^\mathcal{R}$ used in Figs.~\ref{fig:SBTZprelim}--\ref{fig:C1BTZprelim}, with $\mathcal{R} = \mathcal{S}$, $\mathcal{C}_1$, and $\mathcal{C}_2$, are short versions of the notations $h_{\xi, p}^\mathcal{R}$ and $h_{\eta, p}^{\mathcal{R}}$ introduced in Sec.~\ref{subsec:patchdisc} for the meshsizes used in the $p$-th $\mathcal{R}$-type patch. The symbol $\widetilde{\Gamma}^{\mathcal{R}}$, also used in the figures, denotes the inverse image, under the patch parametrization, of the portion of the physical-domain boundary $\Gamma$ contained in that patch.)  A similar procedure is employed for $\mathcal{C}_2$-type patches, but in two stages: first, 1D-BTZ is applied along the lines shown in the left panel of Fig.~\ref{fig:C2BTZprelim}, and, then, a second 1D-BTZ step is performed along the lines shown in the right panel, where both the original patchwise function values and the extension values from the first stage are used in the matching meshes for the second step. The procedure is somewhat more involved in the case of the $\mathcal{C}_1$-type patches: it proceeds in three steps. In the first step, 1D-BTZ is applied along the discrete set of lines shown in the leftmost panel of Fig.~\ref{fig:C1BTZprelim}, using only the patchwise function values on the subregion consisting of points to the right of the dashed line. In the second step, this intermediate extension is subtracted from the original function values, producing a modified function that is zero within this subregion (middle panel, where white represents zero). In the third step, 1D-BTZ is applied along the set of lines shown in the rightmost panel. Since the subregion is now zero, the resulting extension from this step blends smoothly to zero as it approaches the subregion, and this extension is then added to the one from the first step to yield the final smooth extension. (Note that the straightforward procedure of using the 1D-BTZ functions horizontally without an initial zeroing step would generically give rise to 2D-BTZ functions that do not continuously extend the values of the function $f$.)

Finally, the 2D-FC scheme is completed by interpolating the values of the BTZ extension functions from the discrete patch meshes onto the Cartesian mesh. The interpolated values at points in overlap regions are then summed---which, on account of the POU-based decomposition of the near-boundary values of $f$ produces the desired overall 2D-BTZ of $f$; the 2D-FC extension function $f^c$ is then obtained through an application of the FFT algorithm, as detailed in Sec.~\ref{subsec:interpolationontomesh}.

\section{2D-FC method for domains with corners: detailed description}\label{sec:2d-fc}
This section now presents the essential details of the 2D-FC algorithm introduced in Sec.~\ref{sec:2dfc-prelim}.
\begin{remark}\label{rmk:rel-top}
Throughout the remainder of this paper, for any set $B \subseteq \mathbb{R}^2$ we denote by
$\Cl{B}$ (or equivalently $\overline{B}$) the topological closure of $B$ in $\mathbb{R}^2$, and by
$\Int{B}$ its topological interior. 
Additionally, for any set $B \subseteq \overline{\Omega}$ (cf.~\eqref{eq:f}), we further denote by
$\Clrl{\overline{\Omega}}{B}$ and $\Intrl{\overline{\Omega}}{B}$ the closure and interior of $B$
relative to $\overline{\Omega}$, respectively; that is,
$\Clrl{\overline{\Omega}}{B} = \Cl{B} \cap \overline{\Omega}$ and
$\Intrl{\overline{\Omega}}{B} = \overline{\Omega} \setminus \Clrl{\overline{\Omega}}{\overline{\Omega}\setminus B}$.
\end{remark}
\subsection{Boundary-fitted patches and discretizations}\label{subsec:2d-fc-domain-decomp}
The proposed algorithm relies on the decomposition of a certain neighborhood $U$ of $\Gamma$ within $\overline{\Omega}$ (wherein $U$ is a topological relative open subset of $\overline{\Omega}$ which, in particular, satisfies $\Gamma\subset U$) as a union of a number $P = P_\mathcal{S}+P_{\mathcal{C}_1} + P_{\mathcal{C}_2}$ of overlapping boundary-fitted patches, each one of which is a relative open subset of $\overline{\Omega}$, including a number $P_\mathcal{S}$ of $\mathcal{S}$-type patches $\Omega_p^\mathcal{S}$, a number $P_{\mathcal{C}_1}$ of $\mathcal{C}_1$-type patches  $\Omega_p^{\mathcal{C}_1}$, and a number $P_{\mathcal{C}_2}$ of $\mathcal{C}_2$-type patches $\Omega_p^{\mathcal{C}_2}$, such that, in all, we have the decomposition
\begin{equation}\label{eq:decomp}
      U = \left(\bigcup_{p=1}^{P_\mathcal{S}}{\Omega_p^{\mathcal{S}}}\right )\cup\left(\bigcup_{p=1}^{P_{\mathcal{C}_1}}{\Omega_p^{\mathcal{C}_1}}\right)\cup\left(\bigcup_{p=1}^{P_{\mathcal{C}_2}}{\Omega_p^{\mathcal{C}_2}}\right) \subset \overline{\Omega}; \quad \Gamma \subset U.
\end{equation}
Each patch $\Omega_p^\mathcal{R}$ ($\mathcal{R} = \mathcal{S}, \mathcal{C}_1, \mathcal{C}_2$; $p = 1, \dots, P_\mathcal{R}$) serves as the domain of a  ``patchwise localized function'' $f_p^\mathcal{R}$  (a suitably windowed restriction of the given function $f$) that is employed to construct patchwise 2D-BTZ extensions. Each patch $\Omega_p^\mathcal{R}$ is constructed as a subset of the image $V_p^\mathcal{R}$ of a $C^k$-smooth parametrization
\begin{equation}\label{eq:param}
    \mathcal{M}_p^{\mathcal{R}}: \widetilde{V}_p^{\mathcal{R}} \rightarrow V_p^\mathcal{R} \quad \textrm{with}\quad\Omega_p^{\mathcal{R}} \subset {\Intrl{\overline{\Omega}}{V_p^\mathcal{R}\cap\overline{\Omega}}}
\end{equation}
($k\in\mathbb{N}$), where both the parameter space $\widetilde{V}_p^\mathcal{R}$ (defined below) and its image $V_p^\mathcal{R}$ are closed subsets of $\mathbb{R}^2$. Since, per~\eqref{eq:decomp} $U \supset \Gamma$, it follows, in particular, that 
\begin{equation}\label{eq:VpRcover}
    \left(\bigcup_{p=1}^{P_\mathcal{S}}{\Intrl{\overline{\Omega}}{V_p^{\mathcal{S}}\cap\overline{\Omega}}}\right )\cup\left(\bigcup_{p=1}^{P_{\mathcal{C}_1}}{\Intrl{\overline{\Omega}}{V_p^{\mathcal{C}_1}\cap\overline{\Omega}}}\right)\cup\left(\bigcup_{p=1}^{P_{\mathcal{C}_2}}{\Intrl{\overline{\Omega}}{V_p^{\mathcal{C}_2}\cap\overline{\Omega}}}\right) \supset \Gamma;
\end{equation}
note that e.g. $V_p^\mathcal{S}$ equals $V_p^\mathcal{R}$ in~\eqref{eq:param} with $\mathcal{R} = \mathcal{S}$.
Further, the parametrization and its domain and image are constructed in such a way that
\begin{equation}\label{quote0}
\begin{minipage}{0.8\textwidth}
$V_p^\mathcal{R}$ has a non-empty intersection with its two nearest neighbors only.
\end{minipage}
\end{equation}
In the case $\mathcal{R} = \mathcal{S}$ the parameter space $\widetilde{V}_p^\mathcal{S}$ is defined by
\begin{equation}\label{eq:SRtilde}
    \widetilde{V}_p^{\mathcal{S}} = [0, 1] \times [b_{p}^\mathcal{S}, 1] \quad \textrm{for some}\quad b_{p}^\mathcal{S} < 0.
\end{equation}

\begin{figure}[htbp]
\vspace{-0.4cm}
\centering
\includegraphics[width=0.7\linewidth,]{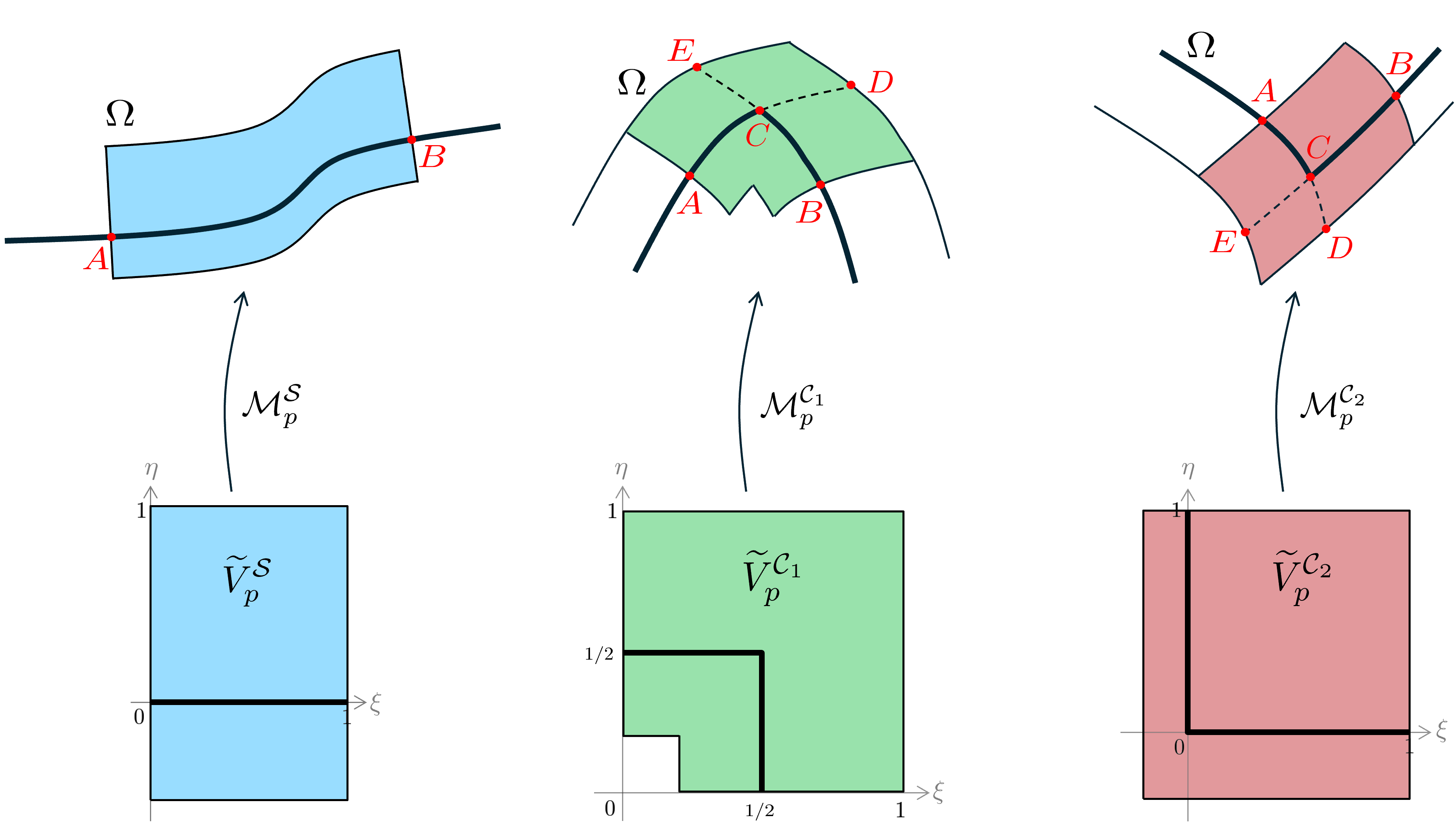}
\vspace{-0.1cm}
\caption{Parametrizations for an $\mathcal{S}$-type patch (left pair of images), a $\mathcal{C}_1$-type patch (middle pair), and a $\mathcal{C}_2$-type patch (right pair). The bottom row shows the parameter domains $\widetilde{V}_p^\mathcal{R}$ for $\mathcal{R} = \mathcal{S}$, $\mathcal{C}_1$, $ \mathcal{C}_2$, and the top row displays their corresponding images $V_p^\mathcal{R}$, in blue, green, and red, respectively.}
    \label{fig:MpR}
    \vspace{-0.2cm}
\end{figure}

Similarly, for $\mathcal{R} = \mathcal{C}_1$ and $\mathcal{C}_2$, we define
\begin{equation}\label{eq:C1Rtilde}
    \widetilde{V}^{\mathcal{C}_1}_p := \left(\left[a_p^{\mathcal{C}_1}, 1\right]\times \left[0, 1\right]\right)  \cup \left(\left[0, 1\right]\times \left[b_p^{\mathcal{C}_1}, 1\right]\right) \quad \textrm{for some}\quad 0 < a_p^{\mathcal{C}_1}, b_p^{\mathcal{C}_1} < \frac{1}{2},
\end{equation}
and
\begin{equation}\label{eq:C2Rtilde}
    \widetilde{V}^{\mathcal{C}_2}_p := [a_p^{\mathcal{C}_2}, 1] \times [b_p^{\mathcal{C}_2}, 1] \quad \textrm{for some} \quad a_p^{\mathcal{C}_2}, b_p^{\mathcal{C}_2} < 0.
\end{equation}
The parametrizations~\eqref{eq:param} for each type $\mathcal{R} = \mathcal{S}, \mathcal{C}_1, \mathcal{C}_2$, together with their associated domains $\widetilde{V}_p^\mathcal{R}$ and images $V_p^\mathcal{R}$, are illustrated in Fig.~\ref{fig:MpR}, and methodologies for their explicit construction are provided in Sec.~\ref{subsec:MpR}.

Employing the notations and conventions introduced above in this section, the open covering~\eqref{eq:decomp} is constructed by defining the $\mathcal{R}$-type patch $\Omega_p^\mathcal{R} \subset V_p^\mathcal{R}$ ($\mathcal{R} = \mathcal{S}, \mathcal{C}_1,$ $\mathcal{C}_2$) as the image under the mapping $\mathcal{M}_p^{\mathcal{R}}(\xi, \eta)$ of a certain subset $\widetilde{\Omega}^\mathcal{R}_p$ of $\widetilde{V}_p^\mathcal{R}$,
\begin{equation}\label{eq:omegapR}
    {\Omega}^\mathcal{R}_p = \mathcal{M}_p^\mathcal{R}({\widetilde{\Omega}}^\mathcal{R}_p) \subset \Intrl{\overline{\Omega}}{V_p^\mathcal{R}\cap\overline{\Omega}},
\end{equation}
where the parameter-space patches $\widetilde{\Omega}_p^\mathcal{R}$ are defined for each one of the patch types $\mathcal{R} = \mathcal{S}$, $\mathcal{C}_1$ and $\mathcal{C}_2$ in terms of certain suitably selected parameters. In detail, we set
 \begin{equation}\label{omegaS}
    {\widetilde{\Omega}}^\mathcal{S}_p = (0, 1) \times \left[0, \eta_p^\mathcal{S}\right), \quad 0 < \eta_p^\mathcal{S} < 1,
\end{equation}
\begin{equation}\label{OmegaC1}
    \widetilde{\Omega}^{\mathcal{C}_1}_p = \left(\left(0, 
    \frac{1}{2}\right] \times \left[\frac{1}{2}, \eta_p^{\mathcal{C}_1}\right)\right) \cup \left(\left[\frac{1}{2}, \xi_p^{\mathcal{C}_1}\right) \times \left(0, \frac{1}{2}\right]\right), \quad \frac{1}{2} < \xi_p^{\mathcal{C}_1}, \eta_p^{\mathcal{C}_1} < 1,
\end{equation}
and
\begin{equation}\label{OmegaC2}
    \widetilde{\Omega}^{\mathcal{C}_2}_p = \left(\left[0, 
    1\right) \times \left[0, \eta_p^{\mathcal{C}_2}\right)\right) \cup \left(\left[0, \xi_p^{\mathcal{C}_2}\right) \times \left[0, 1\right)\right), \quad 0 < \xi_p^{\mathcal{C}_2}, \eta_p^{\mathcal{C}_2} < 1,
\end{equation}
with values the parameters $\eta_p^\mathcal{S}$, $\xi_p^{\mathcal{C}_1}$, $\eta_p^{\mathcal{C}_1}$,  $\xi_p^{\mathcal{C}_2}$, and $\eta_p^{\mathcal{C}_2}$ sufficiently small,  to guarantee that, (i)~The closure $\Clrl{\overline{\Omega}}{\Omega_p^\mathcal{R}}$ is contained entirely in the relative open set $\Intrl{\overline{\Omega}}{V_p^\mathcal{R}\cap\overline{\Omega}} \cup \Intrl{\overline{\Omega}}{V_{p'}^{\mathcal{R}'}\cap\overline{\Omega}} \cup \Intrl{\overline{\Omega}}{{V}_{p''}^{\mathcal{R}''}\cap\overline{\Omega}}$; (ii)~The closure $\Clrl{\overline{\Omega}}{\Omega_p^\mathcal{R}}$ does not intersect the parts of the boundaries of $V_{p'}^{\mathcal{R}'}$ and ${V}_{p''}^{\mathcal{R}''}$ within $\Omega$ (Figure~\ref{fig:MpR}) that are ``parallel'' to $\Gamma$ (or, in the nomenclature of Sec.~\ref{subsec:patchrequirements}, Rem.~\ref{rmk:intedges}, that are ``tangential'' to $\Gamma$); and that, (iii)~For patches ${\Omega_p^\mathcal{R}}$ neighboring a $\mathcal{C}_2$-type patch ${\Omega_{p'}^{\mathcal{C}_2}}$, the closure $\Clrl{\overline{\Omega}}{\Omega_p^\mathcal{R}}$ does not intersect the ``diagonal'' of the parametrized region ${V_{p'}^{\mathcal{C}_2}}$,  along which the ``projection-to-boundary'' map is discontinuous. Full details concerning the selection of the parameters $\eta_p^\mathcal{S}$, $\xi_p^{\mathcal{C}_1}$, $\eta_p^{\mathcal{C}_1}$,  $\xi_p^{\mathcal{C}_2}$, and $\eta_p^{\mathcal{C}_2}$ are provided in  Sec.~\ref{subsec:patchrequirements} and, in particular, points~(a) and~(b) (eq.~\eqref{quote2}) in that section.

Like the 1D-BTZ procedure, where extensions are generated on the basis of function values on the matching mesh $D^\m \subset I^\m$ (see~\eqref{eq:Im}), the 2D-BTZ procedure defines, for each pair $(p, \mathcal{R})$, a patchwise matching domain $\Omega_p^{\mathcal{R}, \m} \subset \Omega_p^{\mathcal{R}}$ and, via a uniform discretization in parameter space, the corresponding patchwise matching mesh $H_{p}^{\mathcal{R}, \m} \subset \Clrl{\overline{\Omega}}{\Omega_p^{\mathcal{R}, \m}}$. The global matching domain is then obtained as the union of all patchwise matching domains,
\begin{equation}\label{eq:Um}
    U^\m = \left(\bigcup_{p=1}^{P_\mathcal{S}}{\Omega_p^{\mathcal{S}, \m}}\right )\cup\left(\bigcup_{p=1}^{P_{\mathcal{C}_1}}{\Omega_p^{\mathcal{C}_1, \m}}\right)\cup\left(\bigcup_{p=1}^{P_{\mathcal{C}_2}}{\Omega_p^{\mathcal{C}_2, \m}}\right) \subset U \subset \overline{\Omega},
\end{equation}
and the corresponding global matching mesh is defined analogously:
\begin{equation}\label{eq:Hm}
      H^\m = \left(\bigcup_{p=1}^{P_\mathcal{S}}{H_p^{\mathcal{S}, \m}}\right )\cup\left(\bigcup_{p=1}^{P_{\mathcal{C}_1}}{H_p^{\mathcal{C}_1, \m}}\right)\cup\left(\bigcup_{p=1}^{P_{\mathcal{C}_2}}{H_p^{\mathcal{C}_2, \m}}\right) \subset \Clrl{\overline{\Omega}}{U^\m}\subset \overline{\Omega}.
\end{equation}
Additional details concerning the construction of the matching domains and meshes $\Omega_p^{\mathcal{R}, \m}$ and $H_p^{\mathcal{R}, \m}$ are presented in Sec.~\ref{subsec:patchdisc}. The 2D-FC algorithm produces the desired 2D-BTZ extension on the basis of the values of~\eqref{eq:f} on the mesh $H^\m \subset \overline{\Omega}$.

\subsection{Partition of unity (POU)}\label{subsec:POU}
The algorithm employs a partition of unity
\begin{equation}\label{eq:POU}
    \mathcal{U} = \{w^\mathcal{R}_p : \mathcal{R} = \mathcal{S},\mathcal{C}_1\mbox{ or } \mathcal{C}_2 \mbox{ and } 1\leq p\leq P_\mathcal{R} \}.
\end{equation}
over $W^\m$ (see Sec.~\ref{sec:2dfc-prelim} and, in particular, Fig.~\ref{fig:POUprelim}) to construct patchwise extension functions that assemble into a smooth global extension. Each smooth window function $w_p^\mathcal{R}$ is defined on a relative open neighborhood $W^\m$ of $\Clrl{\overline{\Omega}}{U^\m}$ (eq.~\eqref{eq:Um}) and is supported within a relative-open ``patchwise extended matching domain'' $\widehat{\Omega}_p^{\mathcal{R},\m}$ satisfying $\Omega_p^{\mathcal{R},\m} \subset \widehat{\Omega}_p^{\mathcal{R},\m} \subset V_p^\mathcal{R}$. Thus, the POU~\eqref{eq:POU} is subordinated to the relative open covering
\begin{equation}\label{Umhat}
    W^\m 
    =\left(\bigcup_{p=1}^{P_\mathcal{S}}\widehat{\Omega}_p^{\mathcal{S},\m}\right)
    \cup\left(\bigcup_{p=1}^{P_{\mathcal{C}_1}}\widehat{\Omega}_p^{\mathcal{C}_1,\m}\right)
    \cup\left(\bigcup_{p=1}^{P_{\mathcal{C}_2}}\widehat{\Omega}_p^{\mathcal{C}_2,\m}\right)
    \subset \overline{\Omega},
\end{equation}
of $W^\m$.  
Explicit constructions of the neighborhood $W^\m$, the domains $\widehat{\Omega}_p^{\mathcal{R},\m}$, and the partition of unity~\eqref{eq:POU} are provided in Sec.~\ref{subsec:POUoverall}.
Note that, by definition, we have
\begin{equation}\label{w_norm}
    \sum_{\stackrel{\mathcal{R} = \mathcal{S}, \mathcal{C}_1, \mathcal{C}_2}{p = 1, \dots, P_\mathcal{R}}} w_p^\mathcal{R}(\mathbf{x}) = 1 \quad \mbox{for all} \quad \mathbf{x} \in \Clrl{\overline{\Omega}}{U^\m}.
\end{equation}

The partition of unity $\mathcal{U}$ is an important element in the 2D-FC algorithm. The 2D-FC method proceeds by first multiplying~\eqref{w_norm} by $f(\mathbf{x})$, to obtain
\begin{equation}\label{eq:fpR}
   f(\mathbf{x})=\sum_{\stackrel{\mathcal{R} = \mathcal{S}, \mathcal{C}_1, \mathcal{C}_2}{1\leq p \leq P_\mathcal{R}}}f_p^\mathcal{R}(\mathbf{x})\quad \textrm{where}\quad f_p^\mathcal{R}(\mathbf{x}) = w_p^\mathcal{R}(\mathbf{x})f(\mathbf{x}), \quad \mathbf{x}\in H^\m,
\end{equation}
and subsequently producing the BTZ extension of $f(\mathbf{x})$ as the sum of the corresponding extensions of the functions $f_p^\mathcal{R}$ for all $(p, \mathcal{R})$.

\subsection{Patchwise extension domains, meshes, and blending-to-zero}\label{subsec:patchwiseBTZ}
Like the smooth-domain 2D-FC algorithm~\cite{bruno_two-dimensional_2022},
the method proposed in this paper extends the function $f$ in~\eqref{eq:f}
into an exterior ``global extension domain''
\begin{equation}\label{exteriordomain}
    U^{\e} = \left(\bigcup_{p=1}^{P_\mathcal{S}}{\Omega_p^{\mathcal{S}, \e}}\right )\cup\left(\bigcup_{p=1}^{P_{\mathcal{C}_1}}{\Omega_p^{\mathcal{C}_1, \e}}\right)\cup\left(\bigcup_{p=1}^{P_{\mathcal{C}_2}}{\Omega_p^{\mathcal{C}_2, \e}}\right) \subset \mathbb{R}^2\setminus \Omega,
\end{equation}
where each ``patchwise extension domain'' $\Omega_{p}^{\mathcal{R}, \e} \subset V_p^\mathcal{R}$ is constructed in such a way that its boundary segment $\Omega_{p}^{\mathcal{R}, \e} \cap \Gamma$ coincides with the corresponding segment $\Omega_{p}^{\mathcal{R}} \cap \Gamma$ of the original patch $\Omega_{p}^\mathcal{R}$. The set $U^\e$ is the 2D analog of the 1D extension interval~\eqref{eq:Ie}. The associated ``global extension mesh'' $H_{\rho}^{\e}$---the natural 2D counterpart of the 1D extension mesh~\eqref{eq:Ie}---is consequently defined by
\begin{equation}\label{exteriormesh}
    H^\e_\rho = \left(\bigcup_{p=1}^{P_\mathcal{S}}{H_{p, \rho}^{\mathcal{S}, \e}}\right )\cup\left(\bigcup_{p=1}^{P_{\mathcal{C}_1}}{H_{p, \rho}^{\mathcal{C}_1, \e}}\right)\cup\left(\bigcup_{p=1}^{P_{\mathcal{C}_2}}{H_{p, \rho}^{\mathcal{C}_2, \e}}\right) \subset \Cl{U^\e},
\end{equation}
where each ``patchwise extension mesh'' $H_{p, \rho}^{\mathcal{R}, \e}$ provides a discretization of $\Cl{\Omega_p^{\mathcal{R}, \e}}$. In view of Figs.~\ref{fig:SBTZprelim},~\ref{fig:C1BTZprelim}, and~\ref{fig:C2BTZprelim}, the extension values for the function $f$ are obtained by repeated application of the $\rho$-refined leftward $\ell$-BTZ procedure (Rems.~\ref{rmk:l-btz} and~\ref{rmk:refinebtz}) along coordinate curves in parameter space $\widetilde{V}_p^\mathcal{R}$ using function values in the corresponding patchwise matching mesh $H_p^{\mathcal{R}, \m}$. As in the 1D case, each patchwise extension domain $\Omega_p^{\mathcal{R}, \e}$ is jointly determined by the 1D-BTZ parameter $C$ and the meshsize in the corresponding patchwise matching mesh. The associated ($\rho$-refined) patchwise extension mesh $H_{p,\rho}^{\mathcal{R},\e}$ is then obtained as the image under the parametrization $\mathcal{M}_p^\mathcal{R}$ of suitably $\rho$-refined meshes in parameter space, as detailed in Sec.~\ref{subsec:emeshes}. 

Utilizing the preceding constructions together with the POU-windowed functions $f_p^\mathcal{R}: H^{\m}\rightarrow \mathbb{C}$ (eq.~\eqref{eq:fpR}), and more specifically, the windowed discrete function $f_p^\mathcal{R}$ over the set of points in the patchwise matching mesh $H_p^{\mathcal{R}, \m}$, the patchwise BTZ grid functions $f_p^{\mathcal{R}, \bt}: H^{\e}_\rho \rightarrow \mathbb{C}$ can be produced, via repeated degree-$d$ 1D-BTZ as described in Sec.~\ref{sec:2dfc-prelim} (and illustrated in Figs.~\ref{fig:SBTZprelim}--\ref{fig:C2BTZprelim}). The grid function $f_p^{\mathcal{R}, \bt}$ is defined on the entire mesh $H^\e_\rho$,  and takes on the value zero at all points outside its corresponding extension domain $H_{p, \rho}^{\mathcal{R}, \e}$. The procedures for constructing $f_p^{\mathcal{R}, \bt}$ in the cases $\mathcal{R} = \mathcal{S}, \mathcal{C}_1,$ and $\mathcal{C}_2$ are described in more detail in Secs.~\ref{subsec:SpatchBTZ}--\ref{subsec:C2patchBTZ}.

\subsection{Interpolation onto the Cartesian mesh}\label{subsec:interpolationontomesh}
For a given mesh size $h>0$, the 2D-FC algorithm produces a Fourier continuation extension $f^c$ over a rectangle
$R = [a_1,b_1]\times[a_2,b_2]$, where $h = (b_1-a_1)/(n_x-1) = (b_2-a_2)/(n_y-1)$ for certain integers $n_x$ and $n_y$, whose interior
$\Int{R} = (a_1,b_1)\times(a_2,b_2)$ satisfies
\begin{equation}\label{eq:Rsup}
\Int{R} \supset \Cl{U^\e} \cup \Omega.
\end{equation}
Rectangles $R$ whose boundaries lie arbitrarily close to the boundary $\partial U^\e$ may be used, since the BTZ extension
vanishes smoothly on $\partial U^\e$. Clearly, the rectangle $R$ is  discretized by the fitted Cartesian mesh
\begin{equation}\label{eq:mesh}
H = \{(a_1+ih,\; a_2+jh)\mid 0\le i\le n_x-1,\; 0\le j\le n_y-1\}\subset R.
\end{equation}

Once the BTZ values on $H_{p, \rho}^{\mathcal{R}, \e}$ have been computed for each pair $(p, \mathcal{R})$ following the procedure outlined in Sec.~\ref{subsec:patchwiseBTZ},  the 2D-FC scheme can be completed by means of the following procedure:
(a) For each $(p,\mathcal{R})$, interpolate the BTZ extension
$f_p^{\mathcal{R},\bt}: H_\rho^\e \to \mathbb{C}$ onto the exterior Cartesian grid points $H \cap \Omega_p^{\mathcal{R},\e}$, and redefine $f_p^{\mathcal{R},\bt}$ as a function $f_p^{\mathcal{R},\bt}: H \cap U^\e \rightarrow \mathbb{C}$ equal to the interpolated values on $H \cap \Omega_p^{\mathcal{R},\e}$ and equal to zero on $H \cap U^\e \setminus H \cap \Omega_p^{\mathcal{R},\e}$; (b)~Sum, at each point on the overall grid $H\cap U^\e$, the interpolated contributions from the functions $f^{\mathcal{R}, \bt}_p$ for all $(p, \mathcal{R})$---which results in the Cartesian values of the overall extension function $f^\bt$ in $U^\e$; (c)~Set $f^\bt(\mathbf{x}) = f(\mathbf{x})$ for all $\mathbf{x} \in H \cap \Omega$; and (d)~Evaluate,  by means of the FFT algorithm, the desired 2D-FC extension
\begin{equation}\label{eq:fc}
    f^c: R \rightarrow \mathbb{C},\qquad f^c(x,y) = \sum_{\ell = -n_x / 2 + 1}^{n_x / 2}
              \sum_{m = -n_y / 2 + 1}^{n_y / 2} \hat{f}^c_{\ell, m}e^{2\pi i\left(\frac{\ell x}{L_x} + \frac{m y}{L_y}\right)}
\end{equation}
with Fourier coefficients $\hat{f}^c_{\ell, m}$ equal to the discrete Fourier transform of the grid function $f^\bt$ in the periodicity domain $R$ (eq.~\eqref{eq:Rsup}).
The grid function $f^\bt$ coincides with the grid values of $f$ within $\Omega$,  it smoothly extends them outside $\Omega$, and it smoothly vanishes outside $\Omega \cup U^\e$---and, therefore, its truncated Fourier expansion $f^c$ closely approximates $f$ within $\Omega$.

The interpolation step in point~(a) must be carefully designed, as both efficiency and accuracy are crucial in this context---given 
the relatively large size of the set $H \cap U^\e$ of Cartesian grid points. The interpolation is performed using the parameter-space 
Cartesian-mesh grid function $\widetilde{f}_p^{\mathcal{R},\bt}$ introduced in Sec.~\ref{SM3}, which constitutes the parameter-space version of the BTZ grid function $f_p^{\mathcal{R},\bt}$ described in Sec.~\ref{subsec:patchwiseBTZ}:
\[
\left(
\widetilde{f}_p^{\mathcal{R},\bt}
\circ
\left(\mathcal{M}_p^{\mathcal{R}}\right)^{-1}
\right)(\mathbf{x})
=
f_p^{\mathcal{R},\bt}(\mathbf{x}), \quad\textrm{for all}\quad \mathbf{x}\in H_{p,\rho}^{\mathcal{R},\e}.
\]
With this notation, the proposed interpolation algorithm 
consists of the following three steps:

 \begin{enumerate}
    \item Identify the global Cartesian points contained in the exterior patch $H \cap \Omega_p^{\mathcal{R}, e}$;
    \item Compute the inverse mapping $ (\xi, \eta) = \left(\mathcal{M}_p^{\mathcal{R}}\right)^{-1}(\mathbf{x}) $ for all $ \mathbf{x} \in H \cap \Omega_p^{\mathcal{R}, \e} $; and,
    \item Interpolate $\widetilde{f}_p^{\mathcal{R},\bt}$ onto the inverse images (Step~2) of the physical-space Cartesian grid points $H \cap \Omega_p^{\mathcal{R},\e}$. The required interpolation from the Cartesian parameter-space extension meshes is performed using the iterated 1D scheme described in~\cite[Sec.~3.6.1]{numer-recipes}, with polynomial degree $M-1$.
\end{enumerate}
Further details concerning Steps~1 and~2 are provided in 
Secs.~\ref{subsec:id-int} and~\ref{subsec:computing-inv-im}, 
respectively, thereby completing the description of the 
algorithmic implementation of point~(a). The implementation 
of Points~(b), (c), and~(d) is straightforward and requires 
no further discussion.

\section{Numerical Results}\label{sec:numerical-results}
This section presents several illustrations demonstrating the character of the proposed 2D-FC method. For quick reference, we note that, as indicated in the preceding sections, the 2D-FC algorithm relies on specific choices for each of the following parameters:
\begin{itemize}
\item[--] $d$: Number of points in each 1D matching mesh, which also equals the order of accuracy of the overall 2D-FC method (Sec.~\ref{sec:fc_1d}).
\item[--] $C$: Number of 1D continuation points (Sec.~\ref{sec:fc_1d}).
\item[--] $Z$: number of 1D zero matching points (Sec.~\ref{sec:fc_1d}).
\item[--] $n_\textrm{os}$: oversampling factor in the 
  matching procedure (Sec.~\ref{sec:fc_1d}).
\item[--] $\rho$: refinement factor for the 1D extension meshes $D^\e$ (Rem.~\ref{rmk:refinebtz}).
\item[--] $R$: rectangle whose interior contains $\overline{U^\e} \cup \Omega$ (Sec.~\ref{sec:2d-fc}).
\item[--] $h$: meshsize of the uniform spatial grid $H$ (Sec.~\ref{sec:2d-fc}).
\item[--] $M-1$: interpolating polynomial degree (Sec.~\ref{subsec:interpolationontomesh}).
\end{itemize}
In addition, the 2D-FC algorithm requires specification of the 
boundary-fitted patches, their parametrizations, and the associated 
discretizations (see Sec.~\ref{subsec:2d-fc-domain-decomp}). In all numerical examples, parameters and patches were selected and constructed, respectively, in accordance with Rems.~\ref{rmk:param_sel} and~\ref{rmk:geometry_sel} below. All runs were performed on a single core of a 2.45 GHz AMD Genoa 9534 processor. The 2D-FC computing times reported throughout this section were from a C implementation of the algorithm.
\begin{remark}[Parameter selections]\label{rmk:param_sel}
  The parameter values $C = 27$, $n_\textrm{os} = 20$, $Z = 12$, and $\rho = 6$ were used in the evaluation of the matrices $\boldsymbol{A}_\rho$ and $\boldsymbol{Q}$ (see Sec.~\ref{sec:fc_1d} and Rem.~\ref{rmk:refinebtz}). The interpolating-polynomial degree $(M-1) = (d + 2)$, in turn, was used for the various matching-point numbers $d$ considered. For each domain considered, a graphical representation of the patch decomposition used is provided.
  \end{remark}
  \begin{remark}[Patch decomposition and discretization]\label{rmk:geometry_sel}
    In the examples below, we consider domains $\Omega$ whose boundaries $\Gamma$ are composed of $N$ smooth curve segments joined at $N$ corner points. The patch decomposition consists of $N$ $\mathcal{S}$-type patches corresponding to the smooth boundary segments, and a total of $N$ corner patches (of type either $\mathcal{C}_1$ or $\mathcal{C}_2$). The patch meshsizes (called $h_{\xi, p}^{\mathcal{R}}$ and $h_{\eta, p}^{\mathcal{R}}$ , in Sec.~\ref{subsec:patchdisc}, $\mathcal{R} = \mathcal{S}$, $\mathcal{C}_1$,  $\mathcal{C}_2$) are chosen to resolve the underlying function within a prescribed error tolerance. Once these meshes are set, the meshsize $h$ of the uniform grid $H$ is selected to ensure adequate resolution of the BTZ extensions. In practice, $h$ was taken to be between one and six times finer than the smallest patch meshsize in the $\mathbf{x}$-variable. The larger refinement factors correspond to patches with the smallest interior (for $\mathcal{C}_2$-type) or exterior (for $\mathcal{C}_1$-type) corner angles, where the BTZ widths become particularly small. 
  \end{remark}

\begin{remark}[Extra vanishing values]
    As in the 1D and smooth 2D settings (cf.~\cite[Remark~3.2]{bruno_two-dimensional_2022}), the grid $H$ may be enlarged prior to the FFT step by appending additional discretization points at which the function is set to zero. This padding allows the total number of grid points in each Cartesian direction to be chosen as a power of two (or, more generally, a product of small prime powers), thereby enabling particularly efficient evaluations via the fast Fourier transform. We employ this padding strategy in all numerical examples considered in this work.
\end{remark}

\begin{figure}[htbp]
    \vspace{-0.25cm}
    \centering
    \includegraphics[width=0.95\linewidth]{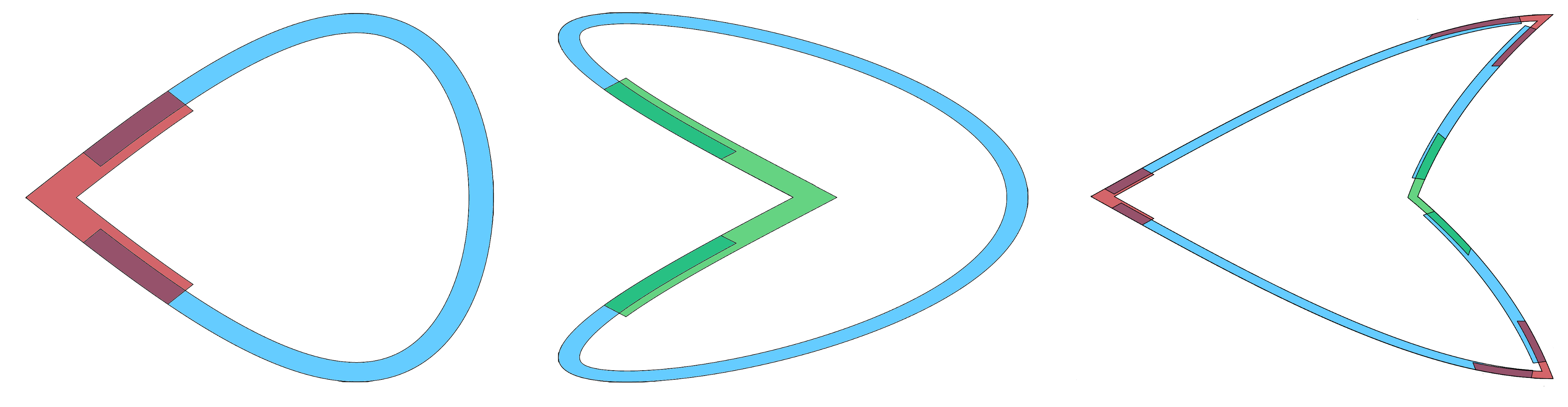}
    \vspace{-0.3cm}
    \caption{Patch decompositions for domains in Exs.~\ref{ex:2dfcteardrop} and~\ref{ex:guitarbase2dfc}: teardrop domain ($\alpha = 1/2$, left), boomerang domain ($\alpha = 3/2$, middle), and guitar-base domain (right).}
    \label{fig:ex1-2-decomp}
\end{figure}

\begin{table}
    {\scriptsize
    \captionsetup{position=top}
    \caption{Convergence of the 2D-FC approximations in Exs.~\ref{ex:2dfcteardrop} and~\ref{ex:guitarbase2dfc}. The FFT cost depends only on the total number of discretization points, which is fixed across different 2D-FC degrees.}\label{tbl:ex1-2-tbl}
    \begin{center}
    \subfloat[Teardrop domain, $\alpha = 1/2$,  (Ex.~\ref{ex:2dfcteardrop}).\label{tbl:ex1teardrop}]{
        \setlength{\tabcolsep}{2pt}
        \makebox[0.95\linewidth][c]{
        \begin{tabular}{@{}c c c|c c c c|c c c c|c@{}}
            \toprule
            & & & \multicolumn{4}{c|}{$d=4$} & \multicolumn{4}{c}{$d=8$}\\
            $h$ & $N_x$ & $N_y$ &
            $\varepsilon_\infty^\mathrm{Rel}$ & $\varepsilon_2^\mathrm{Rel}$ & Ord. & $T_\mathcal{B}$ &
            $\varepsilon_\infty^\mathrm{Rel}$ & $\varepsilon_2^\mathrm{Rel}$ & Ord. & $T_\mathcal{B}$ & $T_\mathcal{F}$\\
            \midrule
            $1\times 10^{-2}$   & 384   & 288   & $4.9\times 10^{-4}$ & $5.6\times 10^{-5}$ & ---  & 0.26 & $8.0\times 10^{-5}$ & $2.7\times 10^{-5}$ & ---  & 0.33 & $2.0\times 10^{-3}$ \\
            $5\times 10^{-3}$   & 576   & 486   & $3.2\times 10^{-5}$ & $2.5\times 10^{-6}$ & 4.5  & 0.44 & $1.1\times 10^{-7}$ & $3.8\times 10^{-9}$ & 12.8 & 0.48 & $4.0 \times 10^{-3}$\\
            $2.5\times 10^{-3}$ & 972   & 972   & $2.0\times 10^{-6}$ & $1.1\times 10^{-7}$ & 4.5  & 0.59 & $1.8\times 10^{-10}$ & $9.6\times 10^{-12}$ & 8.6  & 0.69 & $1.6 \times 10^{-2}$\\
            $1.25\times 10^{-3}$& 1728  & 1728  & $1.3\times 10^{-7}$ & $4.7\times 10^{-9}$ & 4.5  & 0.97 & $1.9\times 10^{-11}$ & $8.4\times 10^{-13}$ & 3.5  & 1.1 & $6.3 \times 10^{-2}$\\
            $6.25\times 10^{-4}$& 3456  & 3456  & $8.3\times 10^{-9}$ & $2.1\times 10^{-10}$ & 4.5  & 1.9  & --- & --- & --- & --- & $2.7 \times 10^{-1}$\\
            $3.125\times 10^{-4}$& 6561 & 6561  & $5.2\times 10^{-10}$& $9.2\times 10^{-12}$ & 4.5  & 4.9  & --- & --- & --- & --- & 1.3\\
            $1.5625\times 10^{-4}$& 13122 & 13122  & $3.4\times 10^{-11}$& $5.2\times 10^{-13}$ & 4.1  & 11.5 & --- & --- & --- & --- & 4.7\\
            \bottomrule
        \end{tabular}
        }
    }
    \newline
    \subfloat[Boomerang domain, $\alpha = 3/2$ (Ex.~\ref{ex:2dfcteardrop}).\label{tbl:ex1boomerang}]{
        \setlength{\tabcolsep}{2pt} 
        \scriptsize
        \makebox[0.95\linewidth][c]{
        \begin{tabular}{@{}c c c|c c c c|c c c c|c@{}}
            \toprule
            & & & \multicolumn{4}{c|}{$d = 5$} & \multicolumn{4}{c}{$d = 7$} & \\
            $h$ & $N_x$ & $N_y$ &$\varepsilon_\infty^\mathrm{Rel}$& $\varepsilon_2^\mathrm{Rel}$ & Ord. & $T_\mathcal{B}$ &$\varepsilon_\infty^\mathrm{Rel}$& $\varepsilon_2^\mathrm{Rel}$& Ord. & $T_\mathcal{B}$ & $T_\mathcal{F}$\\
            \midrule
        $1\times 10^{-2}$ & $192$ & $282$ & $2.6\times 10^{-5}$ & $3.6\times 10^{-6}$ & --- & $0.19$& $6.0\times 10^{-6}$ & $4.6\times 10^{-7}$ & --- & $0.21$ & $1.2\times 10^{-3}$\\
        $5\times 10^{-3}$ & $324$ & $486$ & $8.1\times  10^{-7}$ & $6.6\times 10^{-8}$ & $5.8$ & $0.36$& $1.3\times 10^{-8}$ & $1.1\times 10^{-9}$ & $8.7$ & $0.42$& $3.0 \times 10^{-3}$\\
        $2.5\times 10^{-3}$ & $648$ & $864$ &$2.5\times 10^{-8}$& $1.4\times 10^{-9}$ & $5.6$ & $0.57$ & $9.9 \times 10^{-11}$ & $5.7 \times 10^{-12}$ & $7.6$ & $0.60$ & $8.3 \times 10^{-3}$ \\
        $1.25\times 10^{-3}$ & $1152$ & $1728$ &$7.6\times 10^{-10}$& $3.1\times 10^{-11}$ & $5.5$ & $0.88$ & $5.6\times 10^{-11}$ & $9.6\times 10^{-13}$ & $2.6$ & $0.93$& $3.8\times 10^{-2}$\\
        $6.25\times 10^{-4}$ & $2916$ & $3888$  &$2.4\times 10^{-11}$& $6.8\times 10^{-13}$ & $5.5$ & $2.1$ & \:---\: & \:---\: & \:---\: & \:---\: & $2.7 \times 10^{-1}$\\
            \bottomrule
        \end{tabular}
        }
    }
    \newline
    \subfloat[Guitar-base domain (Ex.~\ref{ex:guitarbase2dfc}).\label{tbl:ex2}]{
        \makebox[0.95\linewidth][c]{
        \begin{tabular}{c c c|c c c c|c}
            \toprule
          $h$ & $N_x$ & $N_y$ & $\varepsilon_\infty^\mathrm{Rel}$ & $\varepsilon_2^\mathrm{Rel}$& Ord. & $T_\mathcal{B}$ & $T_{\mathcal{F}}$ \\
            \midrule
            $5\times 10^{-4}$ & $3072$ & $4374$ & $3.2\times 10^{-9}$ & $1.3 \times 10^{-11}$ & ---& $3.8$ & $3.5 \times 10^{-1}$\\
            $2.5\times 10^{-4}$ & $5832$ & $8748$ & $6.6\times 10^{-11}$ & $1.6 \times 10^{-13}$ & $6.3$& $9.1$ & $1.4$\\
            \bottomrule
          \end{tabular}
        }
    }
    \end{center}
    }
    \vspace{-0.5cm}
\end{table}
\vspace{-0.2cm}
\subsection{2D-FC approximation examples}\label{2dfcex}
All errors reported in this section---namely, the absolute maximum error 
$\varepsilon_\infty^\mathrm{Abs}$, the relative $\ell_2$ error 
$\varepsilon_2^\mathrm{Rel}$, and the relative maximum error 
$\varepsilon_\infty^\mathrm{Rel}$---were computed on the restriction to $\Omega$ of the Cartesian grid of mesh size $h/2$ fitted to the rectangle $R$. The approximations on the finer $h/2$-mesh-size grid were obtained by FFT zero padding~\cite{numer-recipes}.
\begin{example}[Performance and efficiency of the 2D-FC method]\label{ex:2dfcteardrop} 
Our first set of examples demonstrates the performance of the 2D-FC approximation of various degrees $d$. For these tests the function is given by the expression
\begin{equation}
    f(x,y)
    = \exp(0.5(x^2+y^2))(\sin(10\pi x)+\cos(10\pi y))
\end{equation}
on two different domains. The first test utilizes the ``teardrop'' domain, characterized by an interior corner angle
$\alpha\pi$, whose boundary is parameterized by
\begin{equation}
    \mathbf{r}(t) =
    \begin{pmatrix}
        2 \sin \left(\pi t\right)\\
        -\beta \sin(2\pi t)
    \end{pmatrix},
    \quad \beta = \tan \left(\frac{\alpha \pi}{2}\right), \quad 0 < \alpha < 1, \quad 0 \leq t \leq 1.
\end{equation} 
The second is the ``boomerang'' domain, also possessing an interior corner angle $\alpha\pi$, given by
\begin{equation}
    \mathbf{r}(t) =
    \begin{pmatrix}
        -\frac{2}{3} \sin \left(3\pi t\right)\\
        \beta \sin(2\pi t)
    \end{pmatrix},
    \quad \beta = \tan \left(\frac{\alpha \pi}{2}\right), \quad 1 < \alpha < 2, \quad 0 \leq t \leq 1.
\end{equation}
Specifically, we consider the teardrop and boomerang structures with $\alpha = 1/2$ and $\alpha = 3/2$, respectively, which, along with the patch decompositions used, are displayed on the left and middle panels of Fig.~\ref{fig:ex1-2-decomp}. Tbls.~\ref{tbl:ex1teardrop} and~\ref{tbl:ex1boomerang} report the time $T_\mathcal{B}$ required to construct the blending-to-zero function on the grid $H$, and the computing time $T_\mathcal{F}$ needed to evaluate the FFT of the BTZ function in order to obtain the continuation function~\eqref{eq:fc} (both measured in seconds, and averaged over five separate runs), along with the relative errors $\varepsilon_\infty^\mathrm{Rel}$ and $\varepsilon_2^\mathrm{Rel}$. The tables demonstrate the high accuracy and $\mathcal{O}(h^d)$ convergence rate (or higher) of the 2D-FC algorithm. The BTZ construction times scale essentially linearly with $1/h$, while the FFT runtimes exhibit $h^{-2}\log h$ complexity---reflecting the rate of growth of the size of the 2D discretization mesh.
\end{example}

\begin{figure}[htbp]
    \vspace{-0.5cm}
    \centering
    \includegraphics[width=0.8\linewidth]{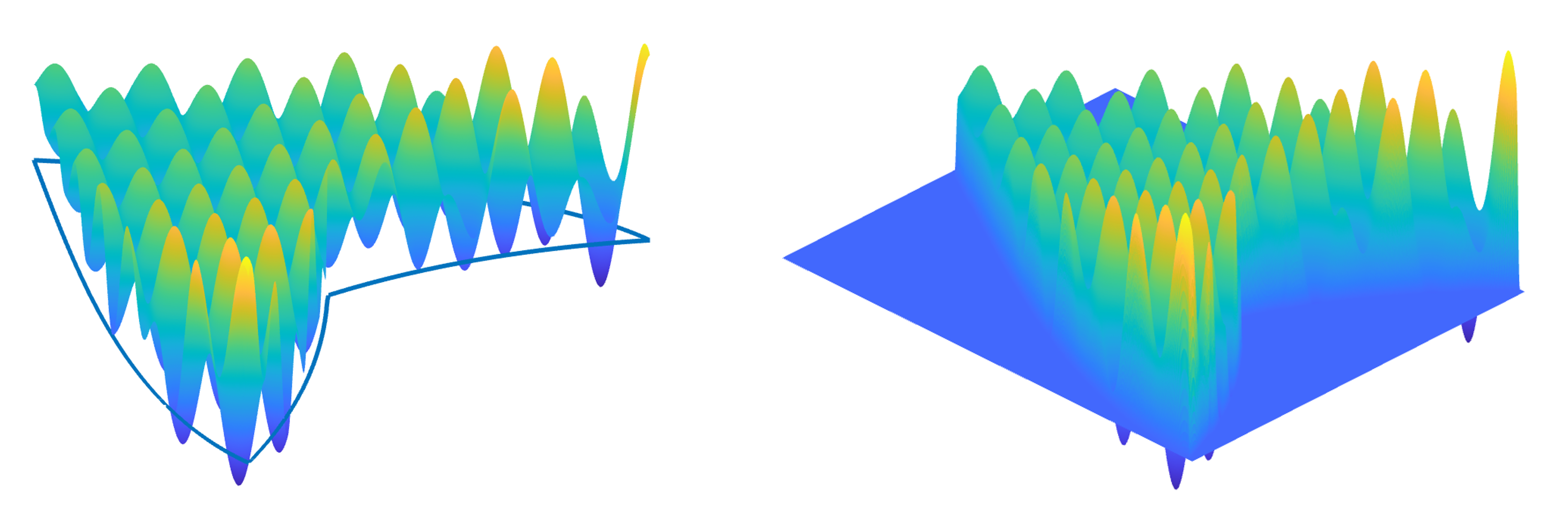}
    \vspace{-0.5cm}
    \caption{Illustration of the 2D-BTZ procedure for the guitar-base domain. Left: original function. Right: blended-to-zero extension. Note the narrow region in which the blending to zero takes place.}
    \label{fig:ex2-illustration}
    \vspace{-0.2cm}
\end{figure}

\begin{example}[2D-FC approximation on a complex domain\label{ex:guitarbase2dfc}]
We next consider the ``guitar-base'' domain displayed in Fig.~\ref{fig:patchdecompprelim}, whose boundary is the union of four smooth curves $\mathbf{r}_j:[0,1]\to\mathbb{R}^2$ ($1\leq j\leq 4)$. Here\looseness = -1
\begin{equation}
\mathbf{r}_1(t) = 
\begin{pmatrix} 
    2\sin(\pi(0.25t))\\ 
    -\sin(2\pi (0.25t))
\end{pmatrix}\quad \textrm{and}\quad \mathbf{r}_4(t) = 
    \begin{pmatrix} 
        2\sin(\pi(0.25t+0.75))\\ 
        -\sin(2\pi (0.25t+0.75)) 
    \end{pmatrix},
\end{equation}
and  $\mathbf{r}_2$ and $\mathbf{r}_3$ are quadratic Bézier curves, with endpoints $\mathbf{a}_2 = \mathbf{r}_1(1)$, $\mathbf{b}_2 = (1,0)^t$ and $\mathbf{a}_3 = (1,0)^t$, $\mathbf{b}_3 = \mathbf{r}_4(0)$ respectively, given by $\mathbf{r}_2(t) = (1-t)^2\mathbf{a}_2 + 2(1-t)t\,\mathbf{c}_2 + t^2\mathbf{b}_2$ and
$\mathbf{r}_3(t) = (1-t)^2\mathbf{a}_3 + 2(1-t)t\,\mathbf{c}_3 + t^2\mathbf{b}_3,$
with $ \mathbf{c}_2 = \tfrac{1}{2}(\mathbf{a}_2+\mathbf{b}_2) + 0.01\,\boldsymbol{\nu}_2,\;
\mathbf{c}_3 = \tfrac{1}{2}(\mathbf{a}_3+\mathbf{b}_3) - 0.01\,\boldsymbol{\nu}_3,\;
\boldsymbol{\nu}_j =
\begin{pmatrix}
0 & -1\\
1 & 0
\end{pmatrix}
\frac{\mathbf{b}_j - \mathbf{a}_j}{\|\mathbf{b}_j - \mathbf{a}_j\|}.$
Tbl.~\ref{tbl:ex2} demonstrates the performance of the 2D-FC approximation of degree $d = 6$ for the function
\begin{equation}
f(x, y) = 4 + (1+x^2+y^2)\bigl(\sin(10.5\pi x - 0.5) + \cos(10\pi y - 0.5)\bigr),
\end{equation}
on the guitar-base domain under consideration; the right panel of Fig.~\ref{fig:ex1-2-decomp} includes the patch decomposition used and Fig.~\ref{fig:ex2-illustration} additionally includes illustrations of the function $f$ and its BTZ extension. Note that very small errors are obtained on account of the fine patch decomposition used and necessary accompanying mesh fineness; naturally, larger errors would result should, say, more oscillatory functions were to be approximated. In any case, these results demonstrate that the 2D-FC method maintains high accuracy even for general domains containing multiple corners.
\end{example}

\begin{figure}[htbp]
    \vspace{-0.2cm}
    \centering
    \includegraphics[width=\linewidth]{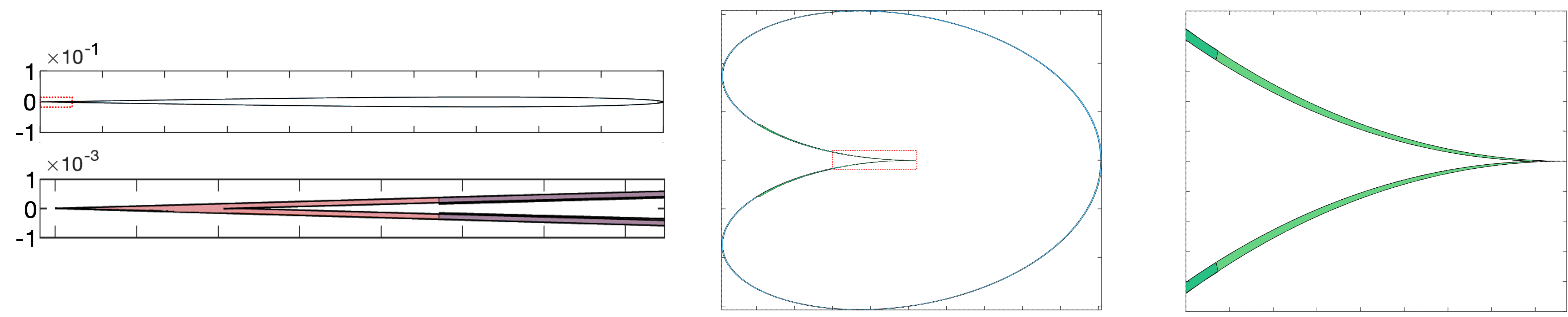}
    \vspace{-0.2cm}
    \caption{Patch decomposition of sharp-corner domains. Left: teardrop domain and its patch decomposition with equal axis scaling (top) and a magnified view of the corner region for clarity (bottom). The patch decomposition in the top-left panel is barely visible since the patches are extremely thin. Middle and right: sharp heart domain, including patches (middle) and corner region, rescaled  for clarity (right). Note the extreme sharpness of the $\mathcal{C}_1$-type patch (green), which necessitates the small step size used in Tbl.~\ref{tbl:ex3heart}. (Use of larger patches might result in geometric inconsistencies.)}
    \label{fig:ex3-decomp}
    \vspace{-0.4cm}
\end{figure}

\begin{example}[Domains with sharp $\mathcal{C}_1$- and $\mathcal{C}_2$-type corners]\label{ex:sharpC1C2}
We conclude this section with a set of two tests involving domains with sharp corner geometries, namely, extreme instances of the ``teardrop'' geometry introduced in Ex.~\ref{ex:2dfcteardrop} and the corresponding ``heart'' geometry  given by
\begin{equation}
    \mathbf{r}(t) = 
    \begin{pmatrix}
        \cos\left(\left(1+\alpha\right)\pi t\right) & -\sin((1+\alpha)\pi t)\\
        \sin((1+\alpha)\pi t) & \cos((1+\alpha)\pi t)
    \end{pmatrix}
    \begin{pmatrix}
        \beta \\ 1
    \end{pmatrix}
    -
    \begin{pmatrix}
        \beta \\ \cos(\pi t)
    \end{pmatrix}, 
    \quad 0 \leq  t \leq 1,
\end{equation}
with $\beta = \tan\left(\frac{\alpha\pi}{2}\right)$. In this example, we consider the rather extreme cases $\alpha = 1/100$ for the teardrop domain and $\alpha = 199/100$ for the heart domain. These curves form $\mathcal{C}_2$-type and $\mathcal{C}_1$-type corners, respectively, with an interior (resp. exterior) angle of $1/100 \pi$ radians $=1.8^\circ$. These ``sharp'' teardrop and heart domains, together with their corresponding patch decompositions, are displayed in Fig.~\ref{fig:ex3-decomp}. In both cases we employ the 2D-FC method of degree~7. For variety two different functions, namely 
\begin{equation}
    f(x,y) = -\big((x+1)^2 + (y+1)^2\big)\sin(\pi(x-0.1))\cos(\pi y),
\end{equation}
and
\begin{equation}
    f(x,y) = (y-1)^2\cos(4\pi x).
\end{equation}
were used as test functions for the $\mathcal{C}_2$-type teardrop domain and the $\mathcal{C}_1$-type heart domain, respectively. The results of these tests are summarized in Tbls.~\ref{tbl:ex3tear} and~\ref{tbl:ex3heart}, respectively. Clearly, the 2D-FC approach is still effective for domains with extreme corner sizes.
\begin{remark}
    The pronounced geometric sharpness of the domains
    considered in this section requires the use of large discretizations, as indicated in Table~\ref{tbl:ex3}, when the
    2D-FC algorithm is applied directly to the complete domain; computing times on the order of 1 or 2 minutes were necessary to produce the results in that table. In practice, this limitation can be avoided
    through a domain decomposition strategy, in which FC approximations are separately
    constructed in localized regions near sharp corners---where fine discretizations
    are required---and in bulk regions away from corners, where substantially coarser
    meshes suffice. The resulting collection of FC expansions can then be employed within a given application. For example, in the context of the Poisson problem
    (cf.~Section~\ref{sec:poisson-problem}), these local FC approximations would be used
    to compute particular solutions on each subdomain, which subsequently would be coupled
    via surface potentials supported on the subdomain interfaces---a strategy that could be incorporated as an addition to the boundary-condition enforcement stage of the Poisson solver algorithm described in Section~\ref{sec:poisson-problem}. No such domain decomposition strategy is employed in the present example, however, which
    is intended solely as an illustration of the ability of the 2D-FC method to
    accurately capture sharply cornered geometries---without regard to the computing cost that is otherwise entailed.
\end{remark}
\end{example}

\begin{table}
    \vspace{-0.25cm}
    \scriptsize
    \captionsetup{position=top}
    \caption{Convergence of the 2D-FC approximations in the sharp teardrop and heart domains in Ex.~\ref{ex:sharpC1C2}.}\label{tbl:ex3}
    \vspace{-0.25cm}
    \begin{center}
    \subfloat[Teardrop domain ($\alpha = 1/100$).\label{tbl:ex3tear}]{
        \makebox[0.95\linewidth][c]{
        \begin{tabular}{c c c |c c c c}
            \toprule
          $h$ & $N_x$ & $N_y$ & $\varepsilon_\infty^\mathrm{Rel}$ & $\varepsilon_2^\mathrm{Rel}$& Ord. \\
            \midrule
            $2.5\times 10^{-5}$ & $82944$ & $1458$ & $1.4 \times 10^{-9}$ & $2.9\times 10^{-12}$ & ---\\
            $1.25\times 10^{-5}$ & $165888$ & $2916$ & $1.1 \times 10^{-11}$ & $1.8\times 10^{-14}$ & $7.3$\\
            \bottomrule
          \end{tabular}
        }
    }
    \newline
    \subfloat[Heart domain ($\alpha = 199/100$).\label{tbl:ex3heart}]{
        \makebox[0.95\linewidth][c]{
        \begin{tabular}{c c c |c c}
            \toprule
          $h$ & $N_x$ & $N_y$ & $\varepsilon_\infty^\mathrm{Rel}$ & $\varepsilon_2^\mathrm{Rel}$ \\
            \midrule
            $2\times 10^{-4}$ & $10368$ & $16384$ & $1.2\times 10^{-10}$ & $2.3 \times 10^{-13}$\\
            \bottomrule
          \end{tabular}
        }
    }
    \end{center}
    \vspace{-0.7cm}
\end{table}

\subsection{Illustrative application: A 2D-FC-based Poisson solver}\label{sec:poisson-problem}
This section demonstrates a 2D-FC-based method for the numerical
solution of the two-dimensional Dirichlet Poisson problem
\begin{equation}\label{eq:def_poisson}
  \begin{cases}
    \Delta u(x, y) = f(x,y), &(x, y)\in\Omega,\\
    u(x, y) = g(x, y), &(x, y)\in \Gamma, 
  \end{cases}
\end{equation}
which generalizes the smooth-domain formulation~\cite{bruno_two-dimensional_2022} to domains containing corner points.
As in the smooth-domain case, the solution $u$ of~\eqref{eq:def_poisson} is expressed as the sum
\begin{equation}\label{eq:total_soln}
  u = u_p + v,
\end{equation}
where $u_p$ denotes a ``particular solution'' of the Poisson equation, $\Delta u_p = f$, computed numerically by means of the 2D-FC method, and where $v$ denotes a ``homogeneous solution'', that is, a solution of Laplace’s equation. The latter is obtained by means of a boundary integral formulation designed to enforce the prescribed Dirichlet boundary data: $v$ corrects the boundary values of $u_p$ so that the total solution $u = u_p + v$ satisfies $u|_{\Gamma} = g$. Accordingly, $v$ satisfies
\begin{equation}\label{eq:def_laplace}
  \begin{cases}
    \Delta v(x, y) = 0, &(x, y)\in\Omega,\\
    v(x, y) = g_{\hom}(x, y), &(x, y)\in \Gamma,
  \end{cases}
\end{equation}
where
\begin{equation}\label{eq:laplace_bc}
  g_{\hom}(x, y) = g(x, y) - u_p(x, y)\big|_{\Gamma}.
\end{equation}

A particular solution $u_p$ for~\eqref{eq:def_poisson} can be  obtained by employing the 2D-FC expansion $f^c(x,y)$ (eq.~\eqref{eq:fc}) of the right-hand side $f(x,y)$; one such solution is given by
\begin{align}\label{eq:poisson_ps}
  u_p(x, y) &= \hat{f}^c_{0, 0} (x^2 + y^2)/4 +
              \sum_{\ell = -N_x / 2 + 1}^{N_x / 2}
              \sum_{m = -N_y / 2 + 1}^{N_y / 2}
              b_{\ell,m} e^{2\pi i\left(\frac{\ell x}{L_x} + \frac{m y}{L_y}\right)},
\end{align}
where
\begin{align}\label{eq:poisson_ps_coeffs}
  b_{\ell,m} =
  \begin{cases}
    0, & (\ell, m) = (0, 0),\\[3pt]
    \displaystyle\frac{-\hat{f}^c_{\ell, m}}
    {(2\pi \ell / L_x)^2 + (2\pi m / L_y)^2}, & (\ell, m)\neq (0, 0).
  \end{cases}
\end{align}
Since the Fourier coefficients $b_{\ell,m}$ are proportional to the 2D-FC coefficients $\hat{f}^c_{\ell,m}$, with proportionality factors that decay quadratically in $(\ell, m)$, it follows that, as demonstrated below in Exs.~\ref{ex:osci-poisson} and~\ref{ex:guitarbase-poison}, the overall numerical solution $u$ converges with order $\mathcal{O}(h^{d+2})$ when a 2D-FC scheme of order $\mathcal{O}(h^{d})$ is used to compute $u_p$---provided that the subsequent evaluation of the homogeneous solution $v$ is performed with sufficient accuracy.

The boundary values $u_p(\boldsymbol{x})$, $\boldsymbol{x}\in\Gamma$, required in~\eqref{eq:laplace_bc}, are obtained by means of a fast algorithm---since direct evaluation of~\eqref{eq:poisson_ps} would lead to an $\mathcal{O}(N^2)$ computational cost. Instead, the proposed algorithm first evaluates the right-hand side of~\eqref{eq:poisson_ps} for all $\boldsymbol{x}\in H$ via FFT, and subsequently computes $u_p$ at boundary points
$\boldsymbol{r}\in\Gamma$ by means of iterated one-dimensional interpolation, following the procedure described in~\cite[Sec.~3.6.1]{numer-recipes}.
To maintain the global $\mathcal{O}(h^{d+2})$ accuracy, polynomial interpolants of degree $(M-1)\ge (d+2)$ (where the same parameter $M$ is used as in the iterated 1D interpolation of Sec.~\ref{subsec:interpolationontomesh}) are utilized in both interpolation directions.

The numerical solution of~\eqref{eq:def_laplace} is obtained efficiently
via the  boundary integral formulation~\cite[eq. (6.35)]{book-kress-lie-2014}. Following~\cite[Sec.~3.6]{colton_inverse_2019} we address the singular behavior of the integral-equation solution near corner points by introducing the change of variables
\begin{equation}\label{eq:graded-map}
  t = w(s), \quad s \in [0,1]
\end{equation}
where the mapping~$w:[0,1]\to[0,1]$ is smooth, monotonically increasing, and
satisfies $w'(0) = w'(1) = 0$, thereby clustering discretization points
near the corner locations.  In practice, for a prescribed grading parameter $p$ ($p \ge 2$), we employ the change of variables~\cite[Sec.~3.6]{colton_inverse_2019}
\begin{equation}\label{eq:graded-explicit}
    w(s) = \frac{v(s)^p}{v(s)^p + \big(v(1 - s)\big)^p}
\end{equation}
along each smooth segment of the boundary $\Gamma$, where each segment is assumed to be parameterized by a vector function defined over the interval $[0,1]$.
The auxiliary function $v$ is given by
\begin{equation}\label{eq:graded-v}
    v(s) = \left(\frac{1}{p} - \frac{1}{2}\right)(1 - 2s)^3 + \frac{1}{p}(2s - 1) + \frac{1}{2}.
\end{equation}

Employing the parametrizations $\mathbf{r}_i(t)$ for the smooth segment $\Gamma_i$, and incorporating the transformation~\eqref{eq:graded-map}, the boundary integral equation for the Laplace Dirichlet problem takes the form
\begin{multline}\label{eq:graded-integral}
    \frac{\varphi(\mathbf{r}_i(w(s)))}{2} - \sum_{j}\int_0^1 K\bigl(\mathbf{r}_i(w(s)), \mathbf{r}_j(w(s'))\bigr)w'(s') \varphi (\mathbf{r}_j(w(s'))) J_j(w(s'))ds'\\ = -g_{\hom}(\mathbf{r}_i(w(s)))
\end{multline}
with $J_j(t) = \|\nabla \mathbf{r}_j(t)\|$. Since we are interested in the re-scaled density $\widehat{\varphi}(\mathbf{r}_i(w(s))) = w'(s) \varphi(\mathbf{r}_i(w(s)))$, we multiply both sides of~\eqref{eq:graded-integral} by $w'(s)$ to get the equivalent formulation
\begin{multline}\label{eq:rescaled-ie}
    \frac{\widehat{\varphi}(\mathbf{r}_i(w(s)))}{2} - w'(s)\sum_{j}\int_0^1 K\bigl(\mathbf{r}_i(w(s)), \mathbf{r}_j(w(s'))\bigr)\widehat{\varphi} (\mathbf{r}_j(w(s'))) J_j(w(s'))ds' \\= -w'(s) g_{\hom}(\mathbf{r}_i(w(s))).
\end{multline}
Discretizing the resulting integral equation by employing the Nystr\"om method described in~\cite[Sec.~12.2]{book-kress-lie-2014}, with a uniform grid in the variable $s$ of mesh size $h$, yields a numerical approximation that converges with an error of order $\mathcal{O}(h^{p-1})$.

Once the particular and homogeneous components $u_p$ and $v$ have been obtained over the Cartesian mesh $H \cap \Omega$, the total solution $u$ of the Poisson problem on that mesh is given by~\eqref{eq:total_soln}; evaluation at other selected sets of points can be produced via evaluation of the Fourier-series and integral-equation solutions, possibly employing relevant fast algorithms if required. For evaluation points $x\not\in\Gamma$ but near $\Gamma$ the integral kernel exhibits near-singular behavior. This difficulty is addressed by a refinement procedure similar to the one described in~\cite[SM1]{bruno_two-dimensional_2022}, which combines local mesh refinement together with the polynomial interpolation of degree $M -1$ that is otherwise used for interpolation in the 2D-FC algorithm (see Sec.~\ref{sec:poisson-eval} for details). 
The numerical convergence rate of the method is primarily governed by the order $d$ of the 2D-FC algorithm employed and the grading parameter $p$. 

\begin{figure}[htb]
    \vspace{-0.5cm}
    \centering
    \includegraphics[width=0.8\linewidth]{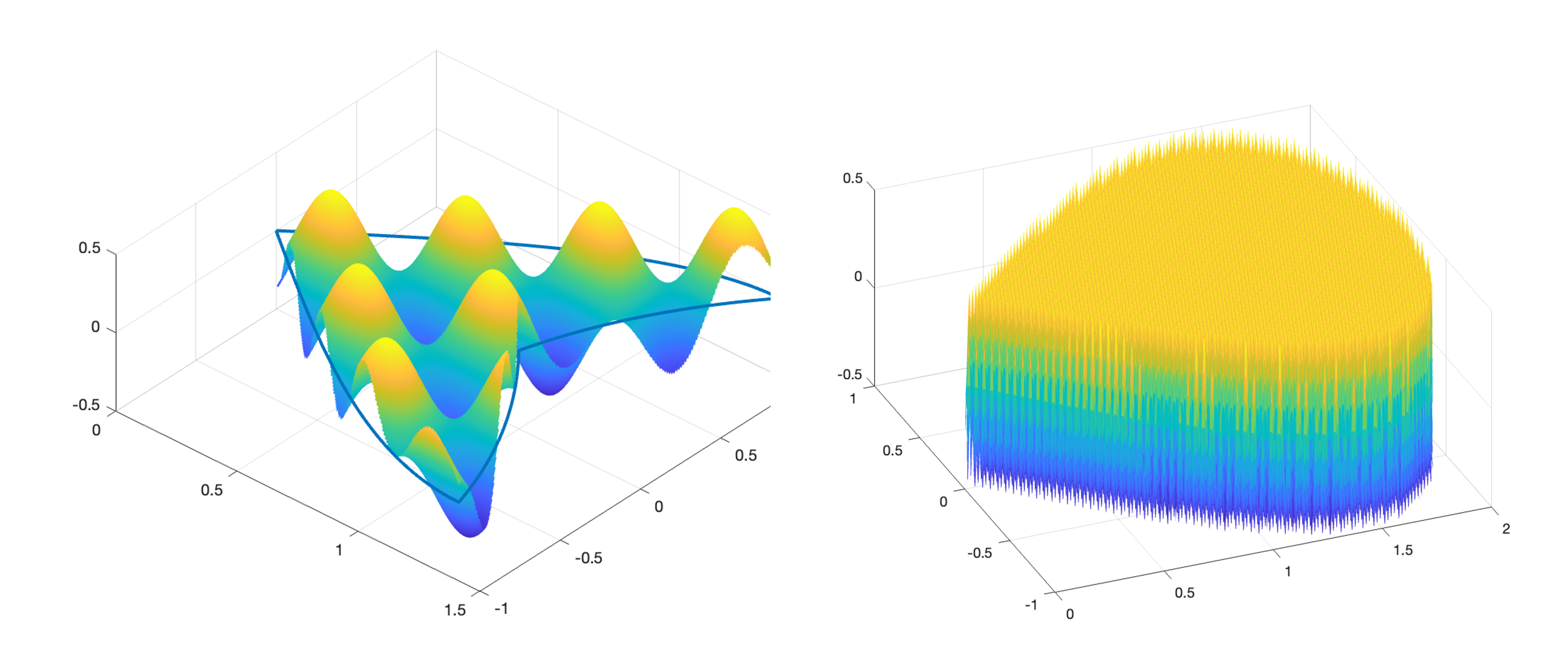}
    \vspace{-0.6cm}
    \caption{
    Numerical solutions of the Poisson equation corresponding to the finest and most accurate discretizations reported in Table~\ref{tbl:poisson-ex}. Left panel: Poisson solution on the guitar-base domain in Ex.~\ref{ex:guitarbase-poison}, obtained from a 2D-FC approximation of degree $d = 7$ and evaluated on a mesh of step size $h = 2\times 10^{-2}$. Right panel: solution to the highly oscillatory Poisson problem in Ex.~\ref{ex:osci-poisson}, obtained from a 2D-FC approximation of degree $d = 10$ and evaluated on a mesh of step size $h = 1\times 10^{-2}$.}
    \label{fig:poisson-exfig}
\end{figure}

Illustrative results for two challenging Poisson problems are presented below. In all cases, the method exhibits fast run times, high accuracy, and rapid convergence. For both examples, and for each chosen value of $d$, we set $p = d+3$ so that the homogeneous solver attains the same $\mathcal{O}(h^{d+2})$ accuracy as the particular solution, and we use the 2D-FC parameter choices specified in Rem.~\ref{rmk:param_sel}. The absolute maximum errors $\varepsilon_\infty^\mathrm{Abs}$ and relative $\ell_2$ errors $\varepsilon_2^\mathrm{Rel}$ are evaluated on the Cartesian mesh $H \cap \Omega$ (see equation~\cref{eq:mesh}). For evaluation of the solution, the Laplace solver contributes a solution error $\varepsilon^\mathrm{int}$, whose estimation is described in Sec.~\ref{sec:poisson-eval}, and which is made to be negligible relative to the overall error based on empirical testing by employing suitable parameter values. Computing times  $T$ required by the 2D-FC stage of the Poisson solver (2D-BTZ and FFT) are reported in Table~\ref{tbl:poisson-ex}. For simplicity, the integral equations used to solve the Laplace problems were treated using an algorithm that does not incorporate acceleration techniques~\cite{Bauinger2021, 2006FMMRokhlin, 1996AIMFFTBelszynski, 2001FFTKunyansky}. Despite the absence of such acceleration, the computational cost of these solves was negligible, with solution times on the order of a fraction of a second for the problems considered in this paper. In contrast, without integral-equation acceleration, the {\em evaluation of the integral potentials} that yield the Laplace solutions represents a non-negligible cost---e.g. on the order of a few seconds to a few tens of seconds on the plotting meshes required to produce the plots in Figure~\ref{fig:poisson-exfig} for the most accurate solutions obtained---but are expected to be reduced to a level comparable to that of the 2D-FC computations, and, indeed, become insignificant for large meshes, once a suitable integral-equation acceleration algorithm is incorporated. 

\begin{table}
    {\scriptsize
    \captionsetup{position=top}
    \caption{Performance of the proposed solver for the Poisson problem on the guitar-base domain considered in Ex.~\ref{ex:guitarbase-poison} and the highly oscillatory problem on the teardrop domain considered in Ex.~\ref{ex:osci-poisson}.}\label{tbl:poisson-ex}
    \vspace{-0.25cm}
    \begin{center}
    \subfloat[Guitar-base results ($d=5$). 
    \label{tbl:poissonguitarbase}]{
        \setlength{\tabcolsep}{2pt}
        \makebox[0.95\linewidth][c]{
        \begin{tabular}{@{}c c c c | c c c c@{}}
            \toprule
            $h$ & $N_x$ & $N_y$ & $\varepsilon^\mathrm{int}$& $\varepsilon_\infty^{\textrm{Abs}}$ & $\varepsilon_2^{\textrm{Rel}}$ & Ord. & $T$ \\
            \midrule
            $4\times 10^{-3}$ & $576$ & $648$ & $10^{-8}$ &$6.5\times10^{-6}$& $3.4\times10^{-7}$ & --- & 0.51 \\
            $2\times 10^{-3}$ & $972$ & $1152$& $10^{-11}$ &$1.6\times10^{-8}$ & $5.1\times10^{-10}$ & $9.4$ & 0.71 \\
            $1\times 10^{-3}$ & $1728$ & $2187$ & $10^{-14}$ & $1.7\times10^{-11}$ & $4.5\times10^{-13}$ & $10.1$ & 1.2 \\
            \bottomrule
          \end{tabular}
        }
    }
    \newline
    \subfloat[Highly-oscillatory Poisson problem on the teardrop domain. A convergence rate close to the expected $\mathcal{O}(h^{12})$ estimate can be observed.\label{tbl:poissontear}]{
    \makebox[0.95\linewidth][c]{
    \begin{tabular}{c c c c | c c c c}
        \toprule
      $h$ & $N_x$ & $N_y$   & $\varepsilon^\mathrm{int}$ & $\varepsilon_\infty^\mathrm{Abs}$ & $\varepsilon_2^\mathrm{Rel}$& Ord. & $T$\\
        \midrule
        $5 \times 10^{-3}$ & $486$ & $486$ & $10^{-7}$ & $4.6\times 10^{-4}$ & $9.2 \times 10^{-5}$ &--- & 0.51\\
        $2.5\times 10^{-3}$ & $972$ & $864$ & $10^{-9}$ & $2.1\times 10^{-7}$ & $3.4 \times 10^{-8}$ & $11.4$ & 0.75\\
        $1.25\times 10^{-3}$ & $1728$ & $1728$ & $10^{-13}$ & $5.7\times 10^{-11}$ & $8.8\times 10^{-12}$ & $11.9$ & $1.2$\\
        \bottomrule
      \end{tabular}
    }
    }
    \end{center}
    }
\end{table}


\begin{example}[Poisson problem for guitar-base domain]\label{ex:guitarbase-poison}
    Our first example for the Poisson problem~\eqref{eq:def_poisson} is posed on the guitar-base domain introduced in Ex.~\ref{ex:guitarbase2dfc}, using the patch decomposition displayed on the right panel of Fig.~\ref{fig:ex1-2-decomp}. For this example we employ the right hand side
    \begin{equation}
            f = (3\pi)^2 \sin(3\pi x-1)\sin(3\pi y-1)
    \end{equation}
   and we enforce the boundary conditions induced by the exact solution
    \begin{equation}
        -\frac{1}{2}\sin(3\pi x-1)\sin(3\pi y-1).
    \end{equation}
    The errors for this problem, with $d=5$ 2D-FC, are reported in Tbl.~\ref{tbl:poissonguitarbase}; once again high-order convergence and high accuracies are observed.
\end{example}

\begin{example}[Highly oscillatory Poisson problem on teardrop domain]\label{ex:osci-poisson}
In our second example, we solve~\eqref{eq:def_poisson} on the teardrop domain in  Ex.~\ref{ex:2dfcteardrop} (with the patch decomposition illustrated in the left panel of Fig.~\ref{fig:ex1-2-decomp}), and for the highly oscillatory right-hand side 
\begin{equation}
    f = (40\pi)^2 \sin(40\pi x-1)\sin(40\pi y-1),
\end{equation}
together with boundary conditions taken from the exact solution 
\begin{equation}
    -\frac{1}{2}\sin(40\pi x-1)\sin(40\pi y-1).
\end{equation}
With $d = 10$, the results in Tbl.~\ref{tbl:poissontear} exhibit convergence near the 12-th order expected from theory. In spite of the strong oscillatory behavior of both the solution and the forcing term and the added complexity induced by the corner point, the observed errors are comparable to those reported for the analogous {\em smooth-domain} boundary-value problem considered in~\cite[Table 4]{bruno_two-dimensional_2022}.
\end{example}

\section{Conclusions}\label{sec:conclusion}
This paper has developed and validated a fast two-dimensional Fourier
Continuation (2D-FC) methodology for the construction of biperiodic
extensions of smooth, non-periodic functions defined on general
two-dimensional domains, including domains with corners.
The results presented demonstrate that the proposed algorithm achieves
the intended $O(N \log N)$ computational complexity for an $N$-point
discretization, while delivering a user-prescribed $d$-th order of
accuracy. In particular, the patch-based construction, partition-of-unity
blending, and extension strategies introduced in this paper enable the method to maintain high-order accuracy in the presence of curved boundaries and corner points. Applications to Poisson problems on bounded domains with corners were presented, which illustrate the practical effectiveness of the 2D-FC approach in application contexts: the 2D-FC
method is capable of producing Poisson solutions with accuracies on the order of machine precision, with computing times on the order of one second on a single core of a modern processor.

\bibliographystyle{siamplain}
\bibliography{references}

\begin{thebibliography}{10}

\bibitem{ALBIN20116248}
{\sc N.~Albin and O.~P. Bruno}, {\em A spectral {FC} solver for the compressible {N}avier–{S}tokes equations in general domains {I}: {E}xplicit time-stepping}, Journal of Computational Physics, 230 (2011), pp.~6248 -- 6270.

\bibitem{amlani_fc-based_2016}
{\sc F.~Amlani and O.~P. Bruno}, {\em An {FC}-based spectral solver for elastodynamic problems in general three-dimensional domains}, Journal of Computational Physics, 307 (2016), pp.~333--354.

\bibitem{askham2017adaptive}
{\sc T.~Askham and A.~J. Cerfon}, {\em An adaptive fast multipole accelerated poisson solver for complex geometries}, Journal of Computational Physics, 344 (2017), pp.~1--22.

\bibitem{Bauinger2021}
{\sc C.~Bauinger and O.~Bruno}, {\em ``{I}nterpolated {F}actored {G}reen {F}unction'' method for accelerated solution of scattering problems}, Journal of Computational Physics, 430 (2021).

\bibitem{1996AIMFFTBelszynski}
{\sc E.~{Bleszynski}, M.~{Bleszynski}, and T.~{Jaroszewicz}}, {\em Aim: Adaptive integral method for solving large-scale electromagnetic scattering and radiation problems}, Radio Science, 31 (1996), pp.~1225--1251.

\bibitem{bruno_cao}
{\sc O.~P. Bruno and J.~Cao}, {\em Fast singular-kernel convolution on general non-smooth domains via truncated fourier filtering}, In preparation,  (2026).

\bibitem{2001FFTKunyansky}
{\sc O.~P. Bruno and L.~A. Kunyansky}, {\em A fast, high-order algorithm for the solution of surface scattering problems: Basic implementation, tests, and applications}, Journal of Computational Physics, 169 (2001), pp.~80--110.

\bibitem{BRUNO20102009}
{\sc O.~P. Bruno and M.~Lyon}, {\em High-order unconditionally stable {FC-AD} solvers for general smooth domains {I}. {B}asic elements}, Journal of Computational Physics, 229 (2010), pp.~2009 -- 2033.

\bibitem{bruno_two-dimensional_2022}
{\sc O.~P. Bruno and J.~Paul}, {\em Two-{Dimensional} {Fourier} {Continuation} and {Applications}}, SIAM Journal on Scientific Computing, 44 (2022), pp.~A964--A992.

\bibitem{yang_2dfc_c}
{\sc O.~P. Bruno and A.~Yang}, {\em {2DFC C}}.
\newblock \url{https://github.com/ayang923/2dfc-c}, 2026.

\bibitem{yang_2dfc_matlab}
{\sc O.~P. Bruno and A.~Yang}, {\em {2DFC MATLAB}}.
\newblock \url{https://github.com/ayang923/2dfc-matlab}, 2026.

\bibitem{2006FMMRokhlin}
{\sc H.~Cheng, W.~Y. Crutchfield, Z.~Gimbutas, L.~F. Greengard, J.~F. Ethridge, J.~Huang, V.~Rokhlin, N.~Yarvin, and J.~Zhao}, {\em A wideband fast multipole method for the helmholtz equation in three dimensions}, Journal of Computational Physics, 216 (2006), pp.~300--325.

\bibitem{colton_inverse_2019}
{\sc D.~L. Colton, R.~Kress, and R.~Kress}, {\em Inverse acoustic and electromagnetic scattering theory}, vol.~93, Springer, 2019.

\bibitem{fortunato2025fully}
{\sc D.~Fortunato, D.~B. Stein, and A.~H. Barnett}, {\em A fully adaptive, high-order, fast poisson solver for complex two-dimensional geometries}, arXiv preprint arXiv:2501.17967,  (2025).

\bibitem{FRYKLUND201857}
{\sc F.~Fryklund, E.~Lehto, and A.-K. Tornberg}, {\em Partition of unity extension of functions on complex domains}, Journal of Computational Physics, 375 (2018), pp.~57 -- 79.

\bibitem{golub2012matrix}
{\sc G.~H. Golub and C.~F.~V. Loan}, {\em Matrix Computations}, Matrix Computations, Johns Hopkins University Press, 2012.

\bibitem{book-kress-lie-2014}
{\sc R.~Kress}, {\em Linear Integral Equations}, Applied Mathematical Sciences, Springer Science and Business Media, third~ed., 2014.

\bibitem{matthysen_huybrechs_fa_frames}
{\sc R.~Matthysen and D.~Huybrechs}, {\em Function approximation on arbitrary domains using fourier extension frames}, SIAM Journal on Numerical Analysis, 56 (2018), pp.~1360--1385.

\bibitem{numer-recipes}
{\sc W.~H. Press, S.~A. Teukolsky, W.~T. Vetterling, and B.~P. Flannery}, {\em Numerical Recipes 3rd Edition: The Art of Scientific Computing}, Cambridge University Press, 3~ed., 2007.

\bibitem{stein2016immersed}
{\sc D.~B. Stein, R.~D. Guy, and B.~Thomases}, {\em Immersed boundary smooth extension: a high-order method for solving pde on arbitrary smooth domains using fourier spectral methods}, Journal of Computational Physics, 304 (2016), pp.~252--274.

\end{thebibliography}


\begin{thebibliography}{1}

\bibitem{bruno_two-dimensional_2022}
{\sc O.~P. Bruno and J.~Paul}, {\em Two-{Dimensional} {Fourier} {Continuation} and {Applications}}, SIAM Journal on Scientific Computing, 44 (2022), pp.~A964--A992.

\bibitem{folland2020introduction}
{\sc G.~B. Folland}, {\em Introduction to partial differential equations}, Princeton University Press, Princeton, NJ, second~ed., 1995.

\bibitem{galetzka2017simple}
{\sc M.~Galetzka and P.~Glauner}, {\em A simple and correct even-odd algorithm for the point-in-polygon problem for complex polygons}, in 12th International joint conference on computer vision, imaging and computer graphics theory and applications (VISIGRAPP 2017), 2017.

\bibitem{book-kress-lie-2014}
{\sc R.~Kress}, {\em Linear Integral Equations}, Applied Mathematical Sciences, Springer Science and Business Media, third~ed., 2014.

\end{thebibliography}
\end{document}


\maketitle

The overall strategy presented in the main body of this paper admits
multiple practical implementations, particularly in its geometric
realization. This section details the specific implementation choices
adopted in this work. Section~\ref{SM1} describes the patch structure
employed, including patch parametrization and
discretization. Section~\ref{sec:SMPOU} outlines the construction of
the Partition of Unity. Section~\ref{SM3} presents the approach used
to effect blending to zero of a given function, including the
extension meshes employed. Section~\ref{SM4} describes the efficient
procedure used for interpolation from the extension meshes to the global
Cartesian mesh. Section~\ref{sec:poisson-eval}, finally, concerns the
evaluation of the boundary integral operators arising in the
application of the 2D-FC method to the Poisson problem, with
particular emphasis on near-boundary and near-corner evaluation.

\section{Parametrized boundary-fitted patches and discretizations}\label{SM1}
As indicated above, this section concerns the construction of patches and their parametrization and associated discretization.
\begin{remark}\label{rmk:partial-notaiton}
    With reference to Rem.~\ref{rmk:rel-top}, for each set $B\subset \mathbb{R}^2$, $\Bnd{B} = \Cl{B} \setminus \Int{B}$ denotes the boundary of $B$ in the topology of $\mathbb{R}^2$. Similarly, for any set $B \subset \overline{\Omega}$,  $\Bndrl{\overline{\Omega}}{B} = \Clrl{\overline{\Omega}}{B} \setminus\Intrl{\overline{\Omega}}{B}$ denotes the boundary of $B$ in the  relative topology of $\overline{\Omega}$.
\end{remark}
\begin{figure}[phtb]
\centering
\includegraphics[width=0.45\linewidth,]{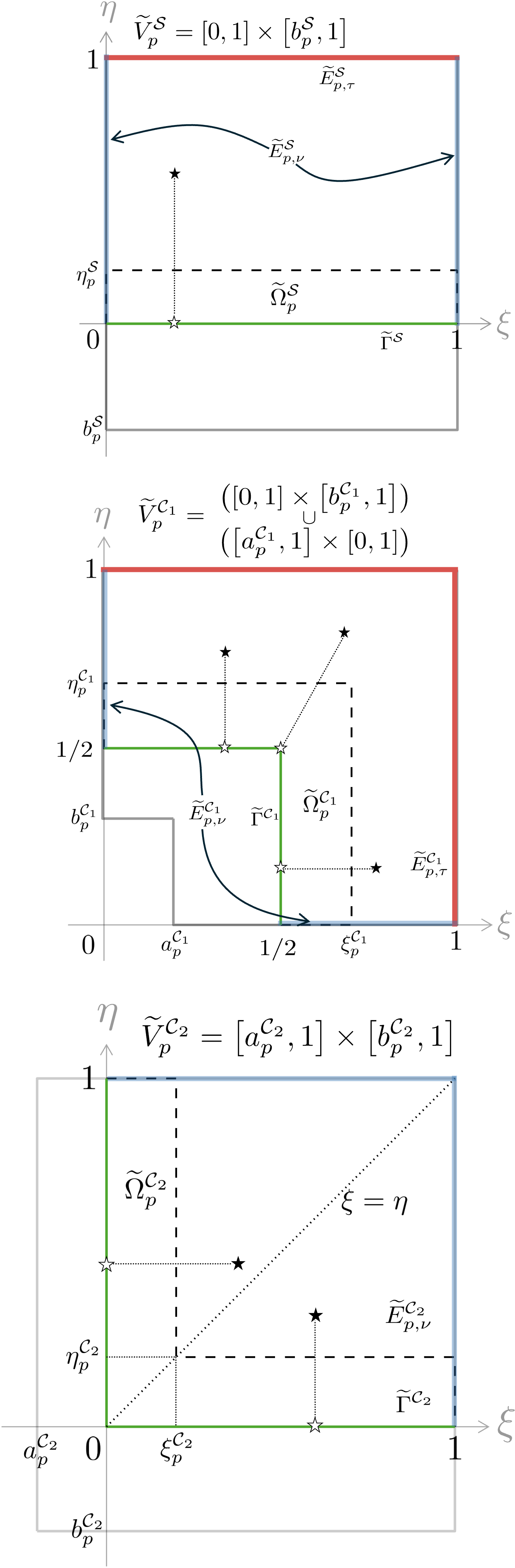}
\caption{Parametrization domain $\widetilde{V}_p^\mathcal{R}$ for the patch types $\mathcal{R} = \mathcal{S}$ (top panel), $\mathcal{R} = \mathcal{C}_1$ (middle panel), and $\mathcal{R} = \mathcal{C}_2$ (bottom panel), including, in each case the parameter-space patch $\widetilde{\Omega}_p^\mathcal{R}$, the parameter-space boundary portion $\Gamma_p^\mathcal{R}$ (green), the transverse edges $\widetilde{E}_{p,\nu}^\mathcal{R}$ (blue), and, for $\mathcal{R} = \mathcal{S}$ and $\mathcal{R} = \mathcal{C}_1$, the tangential edges $\widetilde{E}_{p,\tau}^\mathcal{R}$ (red).
}
    \label{fig:VpS}
\end{figure}
\subsection{Patch parametrizations}\label{subsec:MpR}
This section provides the explicit constructions of the parametrizations $\mathcal{M}_p^\mathcal{R}(\xi, \eta)$ (eq.~\eqref{eq:param}) for the three types of patches $\mathcal{R} = \mathcal{S}$, $\mathcal{C}_1$, and $\mathcal{C}_2$. In order to obtain the necessary parametrization $ \mathcal{M}_p^{\mathcal{R}}(\xi, \eta)$ in the case $\mathcal{R} = \mathcal{S}$ along a smooth arc $\arc{AB}\ \subset \Gamma$ joining the points $A,B\in\Gamma$ we utilize the underlying parametrization of $\Gamma$ and, employing an appropriate rescaling, we first obtain from it a smooth invertible parametrization $\ell_A^B: [0, 1] \rightarrow \mathbb{R}^2$ ($\ell_A^B = \ell_A^B(\xi)$) such that
\begin{equation}
    \ell_A^B([0, 1]) =\ \arc{AB}.
\end{equation}
The corresponding $\mathcal{S}$-type patch and its parametrization are then obtained by extruding along the direction of the normal to $\arc{AB}$ pointing toward the interior of $\Omega$. In detail, calling $\nu(\xi)$  the unit inward-pointing normal vector to $\arc{AB}\ \subset\Gamma$ at the point $\ell_A^B(\xi)$ ($\nu:[0, 1]\to\mathbb{R}^2$), the $\mathcal{S}$-patch parametrization is given by
\begin{equation}\label{eq:MpS}
    \mathcal{M}_p^\mathcal{S}(\xi, \eta) = \ell_A^B(\xi) + \eta H \nu(\xi)
\end{equation}
for some adequately selected parameter $H > 0$; see Fig.~\ref{fig:MpR}. 

In the case $\mathcal{R} = \mathcal{C}_1$, in turn, it is necessary to obtain a parametrization of a neighborhood of a concave corner point $C \in \Gamma$. To do this we note that, in some neighborhood of $C$, $\Gamma$ equals  the union of two smooth arcs $\arc{AC}$ and $\arc{BC}$ that meet at $C$. In view of the assumed smoothness of the arcs $\arc{AC}$ and $\arc{BC}$ there exist  ``smooth extensions'' $\arc{AD} \supset \arc{AC}$ and $\arc{BE} \supset \arc{BC}$ of the curves past the corner point $C$, as illustrated in Fig.~\ref{fig:MpR}, which are utilized as elements in the patch-construction process. (In practice such extensions can be obtained from closed form expressions for various sections of the surface $\Gamma$ which are, e.g.,  commonly found in most CAD standards (Computer Aided Design), or, otherwise, via an application of the 1D FC method.) We assume, as we may, that the extended curves are given by smooth parametrizations  $\ell_A^D: [0, 1] \rightarrow \mathbb{R}^2$ and $\ell_B^E: [0, 1] \rightarrow \mathbb{R}^2$ such that
\begin{align}
    &\ell_A^D([0, \frac{1}{2}]) = \arc{AC}, \quad \ell_A^D([\frac{1}{2}, 1]) = \arc{CD}\\
    &\ell_B^E([0, \frac{1}{2}]) = \arc{BC}, \quad \ell_B^E([\frac{1}{2}, 1]) = \arc{CE}.
\end{align}
Using these curves we then obtain the desired $\mathcal{C}_1$-patch parametrization $\mathcal{M}_p^{\mathcal{C}_1}: \widetilde{V}^{\mathcal{C}_1}_p \rightarrow V_p^{\mathcal{C}_1}$ (cf. eq.~\eqref{eq:C1Rtilde})  given by 
\begin{equation}\label{eq:MpC1}
    \mathcal{M}_p^{\mathcal{C}_1}(\xi, \eta) = \ell_A^D(\xi) + \ell_B^E(\eta) - C;
\end{equation}
see Fig.~\ref{fig:MpR}, where the corner point is given by $C =  \mathcal{M}_p^{\mathcal{C}_1}\left(\frac{1}{2}, \frac{1}{2}\right)$.

In the case $\mathcal{R} = \mathcal{C}_2$, finally, letting once again $C \in \Gamma$ denote the corresponding corner point and assuming the smooth arcs $\arc{AC}$ and $\arc{BC}$ have smooth extensions $\arc{AD} \supset \arc{AC}$ and $\arc{BE} \supset \arc{BC}$ that are parametrized by smooth and invertible maps $\ell_D^A: [a_{p}^{\mathcal{C}_2}, 1] \rightarrow \mathbb{R}^2$ and $\ell_E^B: [b_{p}^{\mathcal{C}_2}, 1] \rightarrow \mathbb{R}^2$ satisfying
\begin{align}
    &\ell_D^A([a_p^{\mathcal{C}_2}, 0]) = \arc{DC}, \quad \ell_D^A([0, 1]) = \arc{CA}\\
    &\ell_E^B([b_p^{\mathcal{C}_2}, 0]) = \arc{EC},\quad \ell_E^B([0, 1]) = \arc{CB}
\end{align}
we obtain the $\mathcal{C}_2$-patch parametrization $\mathcal{M}_p^{\mathcal{C}_2}: \widetilde{V}^{\mathcal{C}_2}_p \rightarrow V_p^{\mathcal{C}_2}$ (cf. eq.~\eqref{eq:C2Rtilde})  given by
\begin{equation}\label{eq:MpC2}
    \mathcal{M}_p^{\mathcal{C}_2}(\xi, \eta) = \ell_D^A(\xi) + \ell_E^B(\eta) - C;
\end{equation}
see Fig.~\ref{fig:MpR}; in this case the corner point is given by $C =  \mathcal{M}_p^{\mathcal{C}_2}\left(0, 0\right)$.
\begin{remark}\label{rmk:MpRsmooth}
    Clearly, the parametrizations~\eqref{eq:MpS},~\eqref{eq:MpC1} and~\eqref{eq:MpC2} satisfy the $C^k$ parametrization smoothness requirement ($k\in\mathbb{N}$) introduced at the beginning of Sec.~\ref{subsec:2d-fc-domain-decomp} provided the boundary $\Gamma$ is itself $C^{k+1}$. Further, since the Jacobians of the transformations~\eqref{eq:MpS},~\eqref{eq:MpC1} and~\eqref{eq:MpC2} do not vanish along the patch boundaries, an application of the inverse function theorem shows (cf. e.g.~\cite[Prop. 0.2]{folland2020introduction}) that, provided $k\geq 1$, the corresponding inverse functions are also continuously differentiable up to order $k$ in suitable transverse neighborhoods of the patch boundaries. 
\end{remark}

\subsection{Additional geometric conditions ensuring smooth POU construction on and near $\Gamma$}\label{subsec:patchrequirements}
In order to establish certain necessary additional requirements ensuring smooth POU construction we first introduce three auxiliary concepts. The first one is the parameter-space ``patchwise boundary portion''
\begin{equation}\label{eq:gammaRp_param}
    \widetilde{\Gamma}^\mathcal{R} =
    \begin{cases}
        (0, 1) \times \{0\}, & \mathcal{R} = \mathcal{S},\\
        \left(\left(0, \frac{1}{2}\right] \times \left\{\frac{1}{2}\right\}\right)
        \cup 
        \left(\left\{\frac{1}{2}\right\} \times \left(0, \frac{1}{2}\right]\right), 
        & \mathcal{R} = \mathcal{C}_1,\\
        \left([0, 1) \times \{0\}\right) 
        \cup 
        \left( \{0\} \times [0, 1)\right), 
        & \mathcal{R} = \mathcal{C}_2,
    \end{cases}
\end{equation}
which in view of eqs.~\eqref{eq:MpS}, \eqref{eq:MpC1}, \eqref{eq:MpC2}, see also Fig.~\ref{fig:MpR}, satisfies $\mathcal{M}_p^\mathcal{R}(\widetilde{\Gamma}^\mathcal{R}) = \Gamma \cap V_p^\mathcal{R}$. The second concept alluded to above is the ``boundary projection along coordinate lines'' mapping $\mathcal{P}_p^\mathcal{R}: V_p^\mathcal{R}\cap \overline{\Omega}\rightarrow \Gamma_p^\mathcal{R}$ which, for a smooth patch ($\mathcal{R}=\mathcal{S}$), is given by $\mathcal{P}_p^\mathcal{S}\left(\mathcal{M}_p^\mathcal{S}(\xi, \eta)\right) = \mathcal{M}_p^\mathcal{S}(\xi, 0)$, and which for other patch types is defined in Rem.~\ref{rmk:ptoGamma} below. The third and final concept concerns ``interior edges'' within $V_p^\mathcal{R} \cap \overline{\Omega}$, and associated subclasses of edges, namely ``tangential interior edges'', which run parallel to the boundary $\Gamma$, and ``transverse interior edges'', which intersect it at a non-vanishing angle. The precise definitions of these sets are provided in Rem.~\ref{rmk:intedges}.

The parameter-space boundary portions~\eqref{eq:gammaRp_param} are used to define the physical patches $\Omega_p^{\mathcal{R}}$ in such a way that, letting 
\begin{equation}\label{eq:gammaRp}
    \Gamma_p^\mathcal{R} = \Gamma \cap \Omega_p^\mathcal{R},
\end{equation}
the relations
\begin{equation}\label{eq:gammapReq}
    \Gamma_p^\mathcal{R} = \mathcal{M}_p^{\mathcal{R}}(\widetilde{\Gamma}^\mathcal{R}) = \Gamma \cap V_p^\mathcal{R}
\end{equation}
are satisfied. In view of~\eqref{eq:VpRcover}, the images of the sets $\widetilde{\Gamma}^\mathcal{R}$ form a relatively open covering of the boundary $\Gamma$, 
\begin{equation}\label{eq:gammaRpcov}
    \Gamma = 
    \left( \bigcup_{p=1}^{P_\mathcal{S}} \Gamma_p^\mathcal{S} \right) 
    \cup 
    \left( \bigcup_{p=1}^{P_{\mathcal{C}_1}} \Gamma_p^{\mathcal{C}_1} \right) 
    \cup 
    \left( \bigcup_{p=1}^{P_{\mathcal{C}_2}} \Gamma_p^{\mathcal{C}_2} \right)
    \subset U,
\end{equation}
which can be used to construct the desired smooth POU. Further, with reference to~\eqref{quote0}, to guarantee the smoothness and proper construction of the window functions in Secs.~\ref{subsec:POU} and~\ref{sec:SMPOU}, it is necessary that a patch $\Omega_p^\mathcal{R}$, together with the image sets $V_{p'}^{\mathcal{R}'}$ and $V_{p''}^{\mathcal{R}''}$ that contain the respective two nearest-neighboring patches $\Omega_{p'}^{\mathcal{R}'}$ and $\Omega_{p''}^{\mathcal{R}''}$  satisfy the following two conditions: 
\begin{equation}\label{quote2}
\begin{minipage}{0.8\textwidth}
    (a)~The closure $\Clrl{\overline{\Omega}}{\Omega_p^\mathcal{R}}$ is contained entirely in the relative open set $\Intrl{\overline{\Omega}}{V_p^\mathcal{R}\cap\overline{\Omega}} \cup \Intrl{\overline{\Omega}}{V_{p'}^{\mathcal{R}'}\cap\overline{\Omega}} \cup \Intrl{\overline{\Omega}}{{V}_{p''}^{\mathcal{R}''}\cap\overline{\Omega}}$. Further, (b) the intersection of the closure $\Clrl{\overline{\Omega}}{\Omega_p^\mathcal{R}}$ with the tangential interior edges of $V_{p}^{\mathcal{R}}$, $V_{p'}^{\mathcal{R}'}$, and $V_{p''}^{\mathcal{R}''}$ (see eq.~\eqref{eq:EptauR}), as well as the intersection of $\Clrl{\overline{\Omega}}{\Omega_p^\mathcal{R}}$ with the discontinuity curve of the associated boundary projections $\mathcal{P}_{p'}^{\mathcal{R}'}$ or $\mathcal{P}_{p''}^{\mathcal{R}''}$ (which arises only in the case $\mathcal{R}' = \mathcal{C}_2$ or $\mathcal{R}'' = \mathcal{C}_2$; see \eqref{quote3}), must be empty.
\end{minipage}
\end{equation}
Given that $\Clrl{\overline{\Omega}}{\Omega_p^\mathcal{R}}$, the tangential interior edges~\eqref{eq:EptauR}, and the discontinuity curves~\eqref{quote3} are all compact sets, the validity of~\eqref{quote2} implies that an alternate version of~\eqref{quote2}, wherein $\Clrl{\overline{\Omega}}{\Omega_p^\mathcal{R}}$ is replaced by a suitably small neighborhood of $\Omega_p^\mathcal{R}$ (in the topology of~$\overline{\Omega}$), is also valid.  
In detail, defining, for any set $A \subset \overline{\Omega}$, the open $\varepsilon$-neighborhood of $A$  in the topology of~$\overline{\Omega}$ given by
\begin{equation}\label{eq:epsn}
    N_\varepsilon(A;\overline{\Omega})
    =\{\mathbf{x}\in\overline{\Omega}:\mathrm{dist}(\mathbf{x},A)<\varepsilon\},
\end{equation}
the validity of~\eqref{quote2} for each one of the (finitely many) patches $\Omega_p^\mathcal{R}$  implies that, once again letting $\Omega_{p'}^{\mathcal{R}'}$ and $\Omega_{p''}^{\mathcal{R}''}$ denote the respective two nearest-neighboring patches of  $\Omega_p^\mathcal{R}$, the following statement is valid as well:
\begin{equation}\label{eq:epsnquote2}
\begin{minipage}{0.8\textwidth}
There exists a $\varepsilon_0>0$ such that, for all $\varepsilon<\varepsilon_0$ and for all $(p,\mathcal{R})$,
(a)~$N_\varepsilon\left({\Omega_p^\mathcal{R}};\overline{\Omega}\right) \subset \Intrl{\overline{\Omega}}{V_p^\mathcal{R}} \cup \Intrl{\overline{\Omega}}{V_{p'}^{\mathcal{R}'}} \cup \Intrl{\overline{\Omega}}{{V}_{p''}^{\mathcal{R}''}}$; and,
(b)~$N_\varepsilon\left({\Omega_p^\mathcal{R}};\overline{\Omega}\right)$ does not intersect either the tangential interior edges of $V_{p}^{\mathcal{R}}$, $V_{p'}^{\mathcal{R}'}$, and $V_{p''}^{\mathcal{R}''}$, or the discontinuity curve of the associated boundary projections $\mathcal{P}_{p'}^{\mathcal{R}'}$ or $\mathcal{P}_{p''}^{\mathcal{R}''}$.
\end{minipage}
\end{equation}
As shown in Sec.~\ref{subsec:smoothPOU}, eq.~\eqref{eq:epsnquote2} guarantees that the functions $w_p^\mathcal{R}$ in the partition-of-unity set~\eqref{eq:POU}, constructed in Sec.~\ref{sec:SMPOU}, are smooth. 

The definitions of the patches $\Omega_p^\mathcal{R}$ are given in eqs.~\eqref{eq:omegapR}--\eqref{OmegaC2},
where the patches $\widetilde{\Omega}_p^\mathcal{R}$ in parameter space are defined separately for each one of the patch types $\mathcal{R} = \mathcal{S}$, $\mathcal{C}_1$ and $\mathcal{C}_2$ in terms of the parameters $\eta_p^\mathcal{S}$, $\xi_p^{\mathcal{C}_1}$, $\eta_p^{\mathcal{C}_1}$,  $\xi_p^{\mathcal{C}_2}$, and $\eta_p^{\mathcal{C}_2}$. Clearly,~\eqref{eq:gammapReq} is satisfied by definition while~\eqref{quote2} is satisfied for sufficiently small $\eta_p^\mathcal{S}$, $\xi_p^{\mathcal{C}_1}$, $\eta_p^{\mathcal{C}_1}$,  $\xi_p^{\mathcal{C}_2}$, and $\eta_p^{\mathcal{C}_2}$ .

\begin{remark}\label{rmk:intedges}
    In view of the requirement~\eqref{quote2}, for each image ${V_p^\mathcal{R}}$, we define two types of ``interior edges'': the ``transverse'' interior edges $E_{p, \nu}^\mathcal{R} \subset V_p^\mathcal{R} \cap \overline{\Omega}$ and the ``tangential'' interior edges $E_{p, \tau}^\mathcal{R} \subset V_p^\mathcal{R} \cap \overline{\Omega}$. To do this we note that the parameter space domain $\widetilde{V}_p^\mathcal{R}$ is a polygon whose boundary consists of a finite set of edges. The transverse interior edges are those whose image under the mapping $\mathcal{M}_p^\mathcal{R}$ lie inside $\overline{\Omega}$ and intersect $\Gamma$. Specifically, the transverse interior edges are given by
    \begin{equation}\label{eq:EpnuR}
         E_{p, \nu}^\mathcal{R} = \mathcal{M}_p^\mathcal{R}(\widetilde{E}_{p, \nu}^\mathcal{R}) \subset V_p^\mathcal{R} \cap \overline{\Omega}
     \end{equation}
     where $\widetilde{E}_{p, \nu}^\mathcal{R} \subset \widetilde{V}_p^\mathcal{R}$ is given by 
     \begin{equation}\label{eq:EtpnuS}
         \widetilde{E}_{p, \nu}^\mathcal{S} = \{0, 1\} \times [0, 1]\quad \textrm{for} \quad \mathcal{R} = \mathcal{S};
     \end{equation}
      \begin{equation}\label{eq:EtpnuC1}
         \widetilde{E}_{p, \nu}^{\mathcal{C}_1} = \left(\left[\frac{1}{2}, 1\right]\times \{0\}\right) \cup \left(\{0\}\times\left[\frac{1}{2}, 1\right]\right) \quad \textrm{for} \quad \mathcal{R} = \mathcal{C}_1;
     \end{equation}
     and
      \begin{equation}\label{eq:EtpnuC2}
         \widetilde{E}_{p, \nu}^{\mathcal{C}_2} = \left(\left[0, 1\right] \times \{1\}\right) \cup  \left(\{1\} \times \left[0, 1\right] \right) \quad \textrm{for} \quad \mathcal{R} = \mathcal{C}_2.
     \end{equation}
     Tangential interior edges, by contrast, are those whose images under $\mathcal{M}_p^\mathcal{R}$ lie inside $\overline{\Omega}$ but do not intersect~$\Gamma$:   
    \begin{equation}\label{eq:EptauR}
         E_{p, \tau}^\mathcal{R} = \mathcal{M}_p^\mathcal{R}(\widetilde{E}_{p, \tau}^\mathcal{R}) \subset V_p^\mathcal{R} \cap \Omega
     \end{equation}
     where $\widetilde{E}_{p, \tau}^\mathcal{R} \subset \widetilde{V}_p^\mathcal{R}$ is given by 
     \begin{equation}
         \widetilde{E}_{p, \tau}^\mathcal{S} = [0, 1] \times \{1\}\quad \textrm{for} \quad \mathcal{R} = \mathcal{S};
     \end{equation}
      \begin{equation}
         \widetilde{E}_{p, \tau}^{\mathcal{C}_1} = \left([0, 1] \times \{1\}\right) \cup \left(\{1\} \times [0, 1]\right) \quad \textrm{for} \quad \mathcal{R} = \mathcal{C}_1;
     \end{equation}
     and
      \begin{equation}
         \widetilde{E}_{p, \tau}^{\mathcal{C}_2} = \emptyset \quad \textrm{for} \quad \mathcal{R} = \mathcal{C}_2;
     \end{equation}
     as illustrated in Fig.~\ref{fig:VpS}. Clearly, with these definitions, we have
     \begin{equation}\label{eq:VpRbound}
         \Bndrl{\overline{\Omega}}{V_p^\mathcal{R}\cap \overline{\Omega}}= E_{p, \tau}^\mathcal{R}\cup E_{p, \nu}^\mathcal{R}
     \end{equation}
     and 
     \begin{equation}\label{eq:VpRsm}
        \Intrl{\overline{\Omega}}{V_p^\mathcal{R} \cap \overline{\Omega}} = \left(V_p^{\mathcal{R}} \cap \overline{\Omega}\right) \setminus (E_{p, \tau}^\mathcal{R}\cup E_{p, \nu}^\mathcal{R}).
     \end{equation}
\end{remark}
\begin{remark}\label{rmk:ptoGamma}
    The boundary projection $\mathcal{P}_p^\mathcal{R}\left(\mathbf{x}\right):V_p^\mathcal{R}\cap \Omega \rightarrow \Cl{\Gamma_p^\mathcal{R}}$ of a point $\mathbf{x}\in V_p^\mathcal{R}\cap\Omega$ is defined by
    \begin{equation}
        \mathcal{P}_p^{\mathcal{R}} = \left(\mathcal{M}_p^{\mathcal{R}} \circ \widetilde{\mathcal{P}}_p^\mathcal{R} \circ \left(\mathcal{M}_p^\mathcal{R}\right)^{-1}\right)(\mathbf{x})
    \end{equation}
    where the parameter space boundary projections $\widetilde{\mathcal{P}}_p^\mathcal{R}: \widetilde{V}_p^\mathcal{R} \rightarrow \Cl{\widetilde{\Gamma}^\mathcal{R}}$ are defined by the following selections of points:
    \begin{equation}\label{pS}
         \widetilde{\mathcal{P}}_p^\mathcal{S}\left(\xi, \eta\right) = (\xi, 0)\quad \mbox{for} \quad\mathcal{R} = \mathcal{S};
    \end{equation}
    \begin{equation}\label{pC1}
        \widetilde{\mathcal{P}}_p^{\mathcal{C}_1}\left(\xi, \eta\right) = 
        \begin{cases}
            (\frac{1}{2}, \eta)\quad&\textrm{for}\quad \xi<\frac{1}{2}\\
            (\xi, \frac{1}{2}) \quad &\textrm{for}\quad \eta <\frac{1}{2}\\
            \left(\frac{1}{2}, \frac{1}{2}\right)\quad&\textrm{otherwise}
        \end{cases}
        \quad \mbox{for} \quad\mathcal{R} = \mathcal{C}_1;
    \end{equation}
    and
    \begin{equation}\label{pC2}
        \widetilde{\mathcal{P}}_p^{\mathcal{C}_2}\left(\xi, \eta\right) = 
        \begin{cases}
            (\xi, 0)\quad&\textrm{for}\quad \xi \leq \eta\\
            (0, \eta)\quad&\textrm{for}\quad \xi > \eta\\
        \end{cases}
        \quad \mbox{for} \quad\mathcal{R} = \mathcal{C}_2.
    \end{equation}
    
    In each case, the boundary projection $\mathcal{M}_p^\mathcal{R}(\xi^*, \eta^*)$ of a point $\mathbf{x} = \mathcal{M}_p^\mathcal{R}(\xi, \eta) \in V_p^\mathcal{R}\cap \Omega$ is either a point on $\Gamma_p^\mathcal{R}$ that lies along a coordinate curve in $V_p^\mathcal{R}$ (that is, that lies on the image in real space of a curve where $\xi$ or $\eta$ are constant) that also contains the point $\mathbf{x}$ or, if no coordinate curves that contain the point $\mathbf{x}$ intersect the boundary, a situation that may only occur in the case $\mathcal{R} = \mathcal{C}_1$ (and, then, for $\xi, \eta > \frac{1}{2}$), $\mathcal{M}_p^{\mathcal{C}_1}(\xi^*, \eta^*)$ equals the corner point $\mathcal{M}_p^{\mathcal{C}_1}(\tfrac{1}{2}, \tfrac{1}{2})$. In regard to $\mathcal{C}_2$ patches on the other hand, it is important to note that
    \begin{equation}\label{quote3}
        \begin{minipage}{0.8\textwidth}
            the $\mathcal{R} = \mathcal{C}_2$ map $\widetilde{\mathcal{P}}_p^{\mathcal{C}_2}$ exhibits a discontinuity along the curve $\xi = \eta$ for $0\leq \xi, \eta \leq 1$.
        \end{minipage}
    \end{equation} 
    This line, which is explicitly incorporated into the construction of patches adjacent to each $\mathcal{C}_2$ patch (see condition~\eqref{quote2}), will be referred to as the ``discontinuity curve'' of the corresponding $\mathcal{C}_2$-type patch.
\end{remark}

\subsection{Patchwise matching domains and matching meshes}\label{subsec:patchdisc}
The near-boundary, patchwise-matching domains $\Omega_p^{\mathcal{R},\m}$ and meshes $H_p^{\mathcal{R},\m}$ are constructed in what follows on the basis of the parametrized patches $\Omega_p^\mathcal{R}$ introduced in Sec.~\ref{subsec:2d-fc-domain-decomp}, and using the discretization parameters specified in Rem.~\ref{rmk:omga_cnst} below.
\begin{remark}\label{rmk:omga_cnst}
    For future reference, for each patch $\Omega_p^\mathcal{R}$, the patchwise matching domain $\Omega_p^{\mathcal{R}, \m}$ and its corresponding patchwise matching mesh $H_p^{\mathcal{R}, \m}$ are defined by taking into account its intended uniform-grid BTZ-related discretization within the associated parameter-space domain. In detail, the parameter-space domain for the patch $\Omega_p^{\mathcal{R}, \m}$ is defined by means of appropriately chosen number $n_{\xi, p}^\mathcal{R}\geq 2$ and $n_{\eta, p}^\mathcal{R}\geq 2$ of discretization points in the $\xi$ and $\eta$ parameter-space directions (where $n_{\xi, p}^\mathcal{R}$ and $n_{\eta, p}^\mathcal{R}$ are assumed odd for $\mathcal{R} = \mathcal{C}_1$) together with the associated mesh-sizes $h_{\xi, p}^\mathcal{R} = 1/(n_{\xi, p}^\mathcal{R}-1)$ and $h_{\eta, p}^\mathcal{R} = 1/(n_{\eta, p}^\mathcal{R}-1)$. The mesh-sizes $h_{\xi, p}^\mathcal{R}$ and $h_{\eta, p}^\mathcal{R}$ must be selected to yield a uniform discretization of the interval $[0,1]$ that includes both endpoints. In particular, the use of such mesh-sizes, along with the requirement that $n_{\xi, p}^\mathcal{R}$ and $n_{\eta, p}^\mathcal{R}$ be odd when $\mathcal{R} = \mathcal{C}_1$, guarantees that both endpoints of $\Gamma_p^\mathcal{R}$ (eq.~\eqref{eq:gammaRp}) and the corner point of $\Gamma_p^\mathcal{R}$ if $\mathcal{R} = \mathcal{C}_1$ are contained in $H_p^{\mathcal{R}, \m}$. 
\end{remark}

With reference to Rem.~\ref{rmk:omga_cnst}, for each patch $\Omega_p^\mathcal{R}$ (cf.~\eqref{omegaS}, \eqref{OmegaC1}, \eqref{OmegaC2}) the mesh sizes $h_{\xi,p}^\mathcal{R}$ and $h_{\eta,p}^\mathcal{R}$ are chosen sufficiently small so that, for every discrete value of the boundary parameter along $\widetilde{\Gamma}^\mathcal{R}$ (either $\xi$ or $\eta$, depending on the patch type and, in some cases, on the portion of the patch under consideration), the image under $\mathcal{M}_p^\mathcal{R}$ of the corresponding 1D-BTZ matching interval~\eqref{eq:Im} in the direction transverse to the boundary is entirely contained in $\Omega_p^\mathcal{R}$. In detail we assume that
\begin{equation}\label{eq:Sdh}
    (d-1)h_{\eta, p}^\mathcal{S} < \eta_p^\mathcal{S}, \quad \textrm{for}\quad \mathcal{R} = \mathcal{S};
\end{equation}
\begin{equation}\label{eq:C1dh}
    \frac{1}{2}+(d-1)h_{\xi, p}^{\mathcal{C}_1} < \xi_p^{\mathcal{C}_1}\quad \textrm{and}\quad \frac{1}{2}+(d-1)h_{\eta, p}^{\mathcal{C}_1} < \eta_p^{\mathcal{C}_1}, \quad \textrm{for}\quad \mathcal{R} = \mathcal{C}_1;
\end{equation}
and
\begin{equation}\label{eq:C2dh}
    (d-1)h_{\xi, p}^{\mathcal{C}_2}<\xi_p^{\mathcal{C}_2}\quad \textrm{and}\quad (d-1)h_{\eta, p}^{\mathcal{C}_2} < \eta_p^{\mathcal{C}_2}, \quad \textrm{for}\quad \mathcal{R} = \mathcal{C}_2;
\end{equation}
clearly the selection of such mesh-sizes can be accomplished while satisfying the requirements imposed per Rem.~\ref{rmk:omga_cnst}. 
Thus, the patchwise matching domains $\Omega_p^{\mathcal{R}, \m}$ are defined as
\begin{equation}\label{Omegam}
    \Omega_p^{\mathcal{R}, \m} = \mathcal{M}_p^\mathcal{R}(\widetilde{\Omega}_p^{\mathcal{R}, \m}) \subset \Omega_p^\mathcal{R}
\end{equation}
where the corresponding parameter-space domains are given by
 \begin{equation}\label{omegaSm}
    {\widetilde{\Omega}}^{\mathcal{S}, \m}_p = (0, 1) \times \left[0, (d-1)h_{\eta, p}^\mathcal{S}\right),
\end{equation}
\begin{equation}\label{OmegaC1m}
    \widetilde{\Omega}^{\mathcal{C}_1, \m}_p = \left(\left(0, 
    \frac{1}{2}\right] \times \left[\frac{1}{2},  \frac{1}{2}+(d-1)h_{\eta, p}^{\mathcal{C}_1}\right)\right) \cup \left(\left[\frac{1}{2}, \frac{1}{2}+(d-1)h_{\xi, p}^{\mathcal{C}_1}\right) \times \left(0, \frac{1}{2}\right]\right),
\end{equation}
and
\begin{equation}\label{OmegaC2m}
    \widetilde{\Omega}^{\mathcal{C}_2, \m}_p = \left(\left[0, 
    1\right) \times \left[0, (d-1)h_{\eta, p}^{\mathcal{C}_2}\right)\right) \cup \left(\left[0, (d-1)h_{\xi, p}^{\mathcal{C}_2}\right) \times \left[0, 1\right)\right).
\end{equation}
The associated patchwise matching meshes $H_p^{\mathcal{R}, \m}$ are then constructed by uniformly discretizing the closure $\Clrl{\overline{\Omega}}{\Omega_p^{\mathcal{R}, \m}}$ in parameter space:
\begin{equation}\label{eq:HpRm}
    {H}_p^{\mathcal{R}, \m} = \mathcal{M}_p^\mathcal{R}(\widetilde{H}_p^{\mathcal{R}, \m}) \subset \Omega_p^{\mathcal{R}, \m}
\end{equation}
with
\begin{equation}\label{eq:Hmtilde}
    \widetilde{H}_p^{\mathcal{R}, \m} = \{(\xi_i,\eta_j)\ :\ (\xi_i, \eta_j) = (ih_{\xi, p}^\mathcal{R}, jh_{\eta, p}^\mathcal{R})\ ,\ (i,j)\in \mathcal{I}_p^\mathcal{R}\}.
\end{equation}
Here, letting $\llbracket a, b\rrbracket = [a, b] \cap \mathbb{Z}$, for $\mathcal{R} = \mathcal{S}$ we have set
\begin{equation}
\mathcal{I}_p^\mathcal{S} = \llbracket 0, n_{\xi, p}^\mathcal{S}-1\rrbracket \times \llbracket 0, d-1\rrbracket,
\end{equation}
 and, similarly, calling $m_{\xi, p}^{\mathcal{C}_1}=\frac{n_{\xi, p}^{\mathcal{C}_1}-1}{2}\in \mathbb{N}$ and $ m_{\eta, p}^{\mathcal{C}_1} = \frac{n_{\eta, p}^{\mathcal{C}_1}-1}{2}\in \mathbb{N}$,
\begin{equation}\label{eq:IpC1}
\begin{aligned}
\mathcal{I}_p^{\mathcal{C}_1} = &\left(\left\llbracket m_{\xi, p}^{\mathcal{C}_1}, m_{\xi, p}^{\mathcal{C}_1} + d-1\right\rrbracket \times \left\llbracket 0, m_{\eta, p}^{\mathcal{C}_1} + d-1\right\rrbracket\right)\\&\qquad\qquad\qquad\qquad\qquad \cup \left(\left\llbracket 0, m_{\xi, p}^{\mathcal{C}_1}\right\rrbracket \times \left\llbracket m_{\eta, p}^{\mathcal{C}_1}, m_{\eta, p}^{\mathcal{C}_1}+d-1\right\rrbracket\right) 
\end{aligned}
\end{equation}
and
\begin{equation}
    \mathcal{I}_p^{\mathcal{C}_2} = \left(\llbracket 0, n_{\xi, p}^{\mathcal{C}_2}-1\rrbracket \times \llbracket 0, d-1\rrbracket\right) \cup \left(\llbracket0, d-1\rrbracket \times \llbracket 0, n_{\eta, p}^{\mathcal{C}_2}-1\rrbracket\right).
\end{equation}

\begin{figure}[tbh]
\centering
\includegraphics[width=0.875\linewidth,]{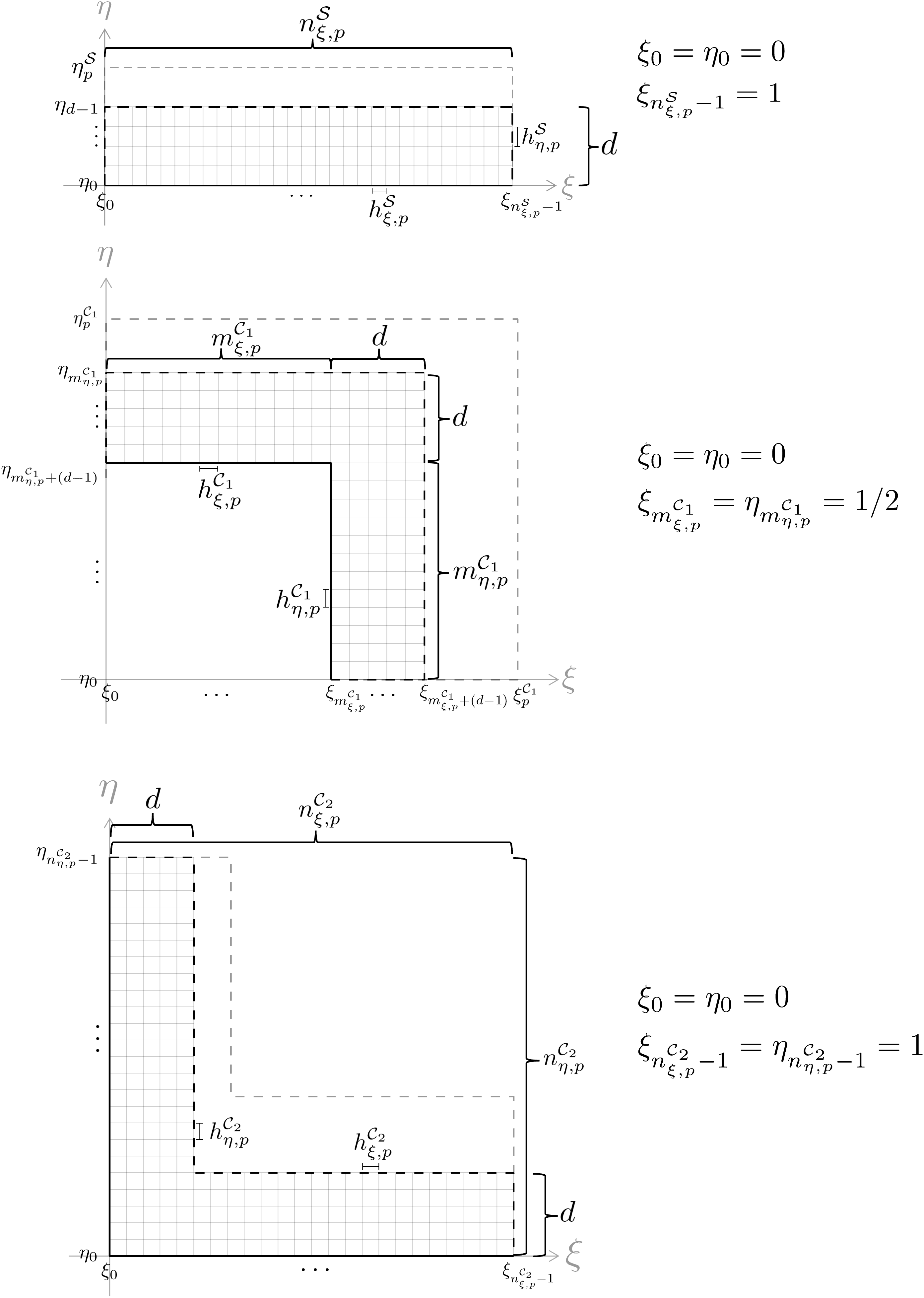}
    \caption{Parameter-space patchwise matching meshes $\widetilde{H}_p^{\mathcal{R}, \m}$
for $\mathcal{R}=\mathcal{S}$ (top panel), $\mathcal{R}=\mathcal{C}_1$ (middle panel),
and $\mathcal{R}=\mathcal{C}_2$ (bottom panel). In each case, mesh lines (gray) are
superimposed on the parameter-space patchwise matching domains
$\widetilde{\Omega}_p^{\mathcal{R},m}$ (not labeled to avoid visual clutter). The parameter-space mesh points are given by
the intersections of the gray grid lines.
}
    \label{fig:Omegae}
\end{figure}

\section{Partition of unity (POU)}\label{sec:SMPOU}
\subsection{Overall POU construction strategy}\label{subsec:POUoverall}
Using the notation introduced in~\eqref{eq:epsn}, let $\varepsilon_0>0$ be as in~\eqref{eq:epsnquote2}, and fix any $\varepsilon\in(0,\varepsilon_0)$. The neighborhood~$W^\m \subset \overline{\Omega}$ in~\eqref{Umhat} is subsequently defined as  
\begin{equation}\label{eq:Wmepsn}
    W^\m
    = \left(\bigcup_{p=1}^{P_\mathcal{S}} N_{\varepsilon}\!\left(\Omega_p^{\mathcal{S}, \m}; \overline{\Omega}\right)\right)
    \cup\left(\bigcup_{p=1}^{P_{\mathcal{C}_1}} N_{\varepsilon}\!\left(\Omega_p^{\mathcal{C}_1, \m}; \overline{\Omega}\right)\right)
    \cup\left(\bigcup_{p=1}^{P_{\mathcal{C}_2}} N_{\varepsilon}\!\left(\Omega_p^{\mathcal{C}_2, \m}; \overline{\Omega}\right)\right).
\end{equation}
The extended matching domains~$\widehat{\Omega}_p^{\mathcal{R}, \m}$ forming the relative open covering in~\eqref{Umhat} are then defined by  
\begin{equation}\label{eq:omegapRhat}
    \widehat{\Omega}_{p}^{\mathcal{R}, \m}
    = \Intrl{\overline{\Omega}}{V_p^\mathcal{R} \cap \overline{\Omega}} \cap W^\m,
\end{equation}
so that each~$\widehat{\Omega}_p^{\mathcal{R}, \m}$ contains all points of neighboring sets $N_{\varepsilon}\left(\Omega_{p'}^{\mathcal{R}', \m}; \overline{\Omega}\right)$ (that are also contained in the parametrization image~$V_p^{\mathcal{R}}$) as well as all points in the intersection of the neighborhood~$N_{\varepsilon}\left(\Omega_p^{\mathcal{R}, \m}; \overline{\Omega}\right)$ with $\Intrl{\overline{\Omega}}{V_p^\mathcal{R} \cap \overline{\Omega}}$. In view of the definition~\eqref{eq:Wmepsn} and condition~(a) of~\eqref{eq:epsnquote2} the union of the sets~$\Intrl{\overline{\Omega}}{V_p^\mathcal{R} \cap \overline{\Omega}}$ taken over all~$(p,\mathcal{R})$ form a covering of $W^\m$; hence, the equality in~\eqref{Umhat} follows directly.

Taking into account the preceding definitions and requirements, the window functions $w_p^\mathcal{R}(\mathbf{x})$ are then constructed through the following two-step procedure (the first step is then detailed in the subsequent section):
\begin{itemize}
    \item[(i)] Selection of the parameters in eqs.~\eqref{eq:wpgs}--\eqref{eq:wpgc2} below that determine a certain pre-POU set $\mathcal{U}^\Gamma$---which consist of 1D, unnormalized, non-negative smooth window functions $w_p^{\mathcal{R},\Gamma}$ ($\mathcal{R} = \mathcal{S},\mathcal{C}_1\mbox{ or } \mathcal{C}_2 \mbox{ and } 1\leq p\leq P_\mathcal{R}$), where $w_p^{\mathcal{R},\Gamma}$ is defined over the boundary $\Gamma$; it is supported on the set $\Clrl{\Gamma}{\Gamma_p^\mathcal{R}}$ (see eq.~\eqref{eq:gammaRp}); in the case $\mathcal{R} = \mathcal{C}_2$ it satisfies $w_p^{\mathcal{R},\Gamma} = 1$ in a certain neighborhood of the corresponding corner point; it satisfies $w_p^{\mathcal{R},\Gamma} = 1$ over the set 
    \begin{equation}
        \{\mathbf{x}:\mathbf{x} \in \Gamma_p^\mathcal{R}, \mathbf{x}\notin \Gamma_{p'}^{\mathcal{R}'} \textrm{ for all }p' \neq p\};
    \end{equation} and where, collectively, the functions $w_p^{\mathcal{R},\Gamma}$ satisfy the relation
    \begin{equation}  
    \sum_{\stackrel{\mathcal{R} = \mathcal{S}, \mathcal{C}_1, \mathcal{C}_2}{p = 1, \dots, P_\mathcal{R}}} w_p^{\mathcal{R},\Gamma}(\mathbf{x})\ne 0 \mbox{ for all } \mathbf{x}\in\Gamma.
    \end{equation}
 
    \item[(ii)] Construction of the POU set $\mathcal{U}$ in~\eqref{eq:POU}---by first obtaining the unnormalized 2D window functions $\widehat{w}_p^{\mathcal{R}}(\mathbf{x})$:
    \begin{equation}\label{wpRhat}
        \widehat{w}_p^{\mathcal{R}}(\mathbf{x}) =
        \begin{cases}
            w_p^{\mathcal{R}, \Gamma}\left(\mathcal{P}_p^\mathcal{R}(\mathbf{x})\right), \quad &\mathbf{x} \in \widehat{\Omega}_p^\mathcal{R},\\
            0, \quad & \mathbf{x} \notin \widehat{\Omega}_p^\mathcal{R}
        \end{cases}
    \end{equation}
    (which incorporate the window functions $w_p^{\mathcal{R}, \Gamma}$ introduced in point~(i) above together with the boundary projection $\mathcal{P}_p^\mathcal{R}$ described in Rem.~\ref{rmk:ptoGamma}), and then defining the normalized window functions $w_p^\mathcal{R} \in\mathcal{U}$ given by
    \begin{equation}\label{wpR}
        w_p^\mathcal{R}(\mathbf{x}) = \widehat{w}_p^\mathcal{R}(\mathbf{x})\left(\sum_{\stackrel{\mathcal{R}' = \mathcal{S}, \mathcal{C}_1, \mathcal{C}_2}{p = 1, \dots, P_\mathcal{R}}}\widehat{w}_{p'}^{\mathcal{R}'}\left(\mathbf{x}\right)\right)^{-1}.
    \end{equation}
\end{itemize}
The numerical algorithm only requires the construction of the POU functions $w_p^\mathcal{R}(\mathbf{x})$ in~\eqref{wpR} for $\mathbf{x}$ in the matching mesh $H^\m$ (eq.~\eqref{eq:Hm}). Since equation~\eqref{wpR} involves procedures carried out on pairs of extended matching domains, $\widehat{\Omega}_p^{\mathcal{R}, \m}$ and $\widehat{\Omega}_{p'}^{\mathcal{R}', \m}$, the construction of the POU functions~\eqref{wpRhat} requires the inversion of the corresponding parametrizations, $\mathcal{M}_p^{\mathcal{R}}$ and $\mathcal{M}_{p'}^{\mathcal{R}'}$. These inversions are performed by means of Newton's method.

Clearly, eq.~\eqref{wpR} guarantees the relation~\eqref{w_norm} is satisfied for all $\mathbf{x}\in H^\m$. Since, additionally, each $w_p^\mathcal{R}$ is a smooth function over the set $W^\m \supset \overline{U^\m}\supset H^\m$ (as shown in Sec.~\ref{subsec:smoothPOU}), the construction of the POU $\mathcal{U}$ is complete.

\subsection{Selection of parameters that determine pre-POU set}\label{subsec:prePOU}

The construction of the smooth one-dimensional pre-POU functions $w_p^{\mathcal{R},\Gamma}$, as required in point~(i) above, utilizes the relative open covering~\eqref{eq:gammaRpcov} of $\Gamma$, and proceeds by employing a one-dimensional ``transition function'' $\widehat{w}^{\Gamma}(\tau)$ that smoothly transitions from $0$ to $1$ as $\tau$ increases over the interval $[0,1]$, and which equals 0 and 1 at $\tau = 0$ and $\tau = 1$, respectively.
Throughout this paper, we employ a transition function based on the complementary error function
\begin{equation}
    \erfc(x) = \frac{2}{\sqrt{\pi}}\int_x^\infty e^{-t^2} dt,
\end{equation}
namely, the transition function
\begin{equation}\label{what}
    \widehat{w}^{\Gamma}(\tau) = \frac{1}{2}\erfc(6(-2\tau+1))
\end{equation}
which, as desired, yields a smooth transition from 0 to 1 in the interval $\tau \in [0, 1]$, and which satisfies $|\widehat{w}^{\Gamma}(0)|< 2\times 10^{-16}$ and $|\widehat{w}^{\Gamma}(1) - 1|< 2\times 10^{-16}$---so that the endpoint values $\widehat{w}^{\Gamma}(0)$ and $\widehat{w}^{\Gamma}(1)$ are within machine precision of 0 and 1, respectively, in double-precision arithmetic.

To  construct the window function $w_p^{\mathcal{S}, \Gamma}$, we first note that the intersections $\Gamma_p^\mathcal{S} \cap \Gamma_{q'}^{\mathcal{R}'}$ and $\Gamma_p^\mathcal{S} \cap \Gamma_{q''}^{\mathcal{R}''}$  of a given boundary portion $\Gamma_p^\mathcal{S}$ and its two neighboring boundary portions $\Gamma_{q'}^{\mathcal{R}'}$ and $\Gamma_{q''}^{\mathcal{R}''}$ are given by 
\begin{equation}\label{eq:gammaSoverlap}
    \Gamma_p^\mathcal{S} \cap \Gamma_{q'}^{\mathcal{R}'} = \mathcal{M}_p^\mathcal{S}\left(\left(0, \xi_{\ell, p}^\mathcal{S}\right)\times \{0\}\right)\quad\textrm{and}\quad \Gamma_p^\mathcal{S} \cap \Gamma_{q''}^{\mathcal{R}''} = \mathcal{M}_p^\mathcal{S}\left(\left(\xi_{r, p}^\mathcal{S}, 1\right)\times \{0\}\right)
\end{equation}
for certain real values $\xi_{\ell, p}^\mathcal{S}$ and $\xi_{r, p}^\mathcal{S}$ satisfying $0 < \xi_{\ell, p}^\mathcal{S} < \xi_{r, p}^\mathcal{S} < 1$. Analogously, to construct the window functions $w_p^{\mathcal{C}_1, \Gamma}$ and $w_p^{\mathcal{C}_2, \Gamma}$ we employ the relations
\begin{equation}
    \Gamma_p^{\mathcal{C}_1} \cap \Gamma_{q'}^{\mathcal{R}'} = \mathcal{M}_p^{\mathcal{C}_1}\left(\left(0, \xi_{\ell, p}^{\mathcal{C}_1}\right)\times\left\{\frac{1}{2}\right\}\right),\quad \Gamma_p^{\mathcal{C}_1} \cap \Gamma_{q''}^{\mathcal{R}''} = \mathcal{M}_p^{\mathcal{C}_1}\left(\left\{\frac{1}{2}\right\}\times\left(0, \eta_{d, p}^{\mathcal{C}_1}\right)\right)
\end{equation}
and
\begin{equation}
    \Gamma_p^{\mathcal{C}_2} \cap \Gamma_{q'}^{\mathcal{R}'} = \mathcal{M}_p^{\mathcal{C}_2}\left(\left(\xi_{r, p}^{\mathcal{C}_2}, 1\right)\times\left\{0\right\}\right),\quad \Gamma_p^{\mathcal{C}_2} \cap \Gamma_{q''}^{\mathcal{R}''} = \mathcal{M}_p^{\mathcal{C}_2}\left(\left\{0\right\}\times\left(\eta_{u, p}^{\mathcal{C}_2}, 1\right)\right)
\end{equation}
for some parameter-space values $ \xi_{\ell, p}^{\mathcal{C}_1}$,  $\eta_{d, p}^{\mathcal{C}_1}$, $ \xi_{r, p}^{\mathcal{C}_2}$, and $\eta_{u, p}^{\mathcal{C}_2}$ satisfying $0 < \xi_{\ell, p}^{\mathcal{C}_1}, \eta_{d, p}^{\mathcal{C}_1} < \frac{1}{2}$ and $0 < \xi_{r, p}^{\mathcal{C}_2}, \eta_{u, p}^{\mathcal{C}_2} < 1$.

\begin{figure}[H]
\centering
\includegraphics[width=1\linewidth,]{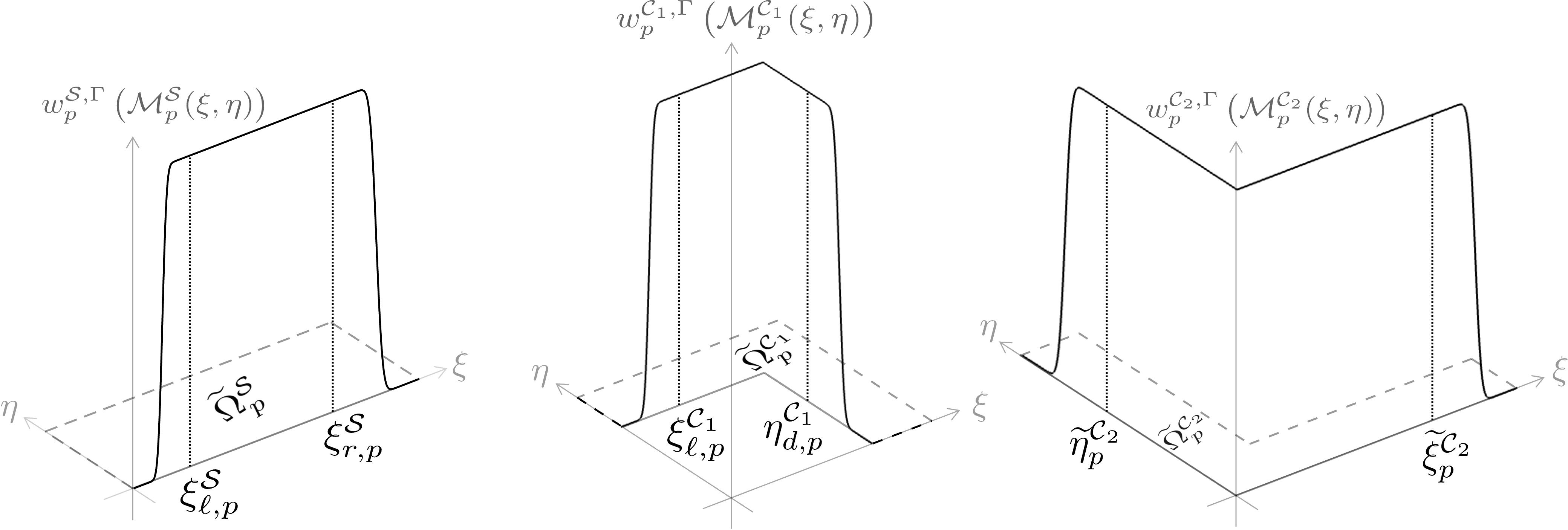}
\caption{Parameter-space depiction of the pre-POU window functions $w_p^{\mathcal{R}, \Gamma}$ for $\mathcal{R} = \mathcal{S}$ (left panel), $\mathcal{R} = \mathcal{C}_1$ (middle panel), and $\mathcal{R} = \mathcal{C}_2$ (right panel). Specifically, in each case the composite function $w_p^{\mathcal{R}, \Gamma} \left(\mathcal{M}_p^\mathcal{R} (\xi,\eta)\right)$ is displayed, with $(\xi,\eta)$ over the parameter space boundary portion $\widetilde{\Gamma}_p^\mathcal{R}$.
}
    \label{fig:wpRgamma}
\end{figure}

Taking into account these relations and using the transition function~\eqref{what}, the value of the pre-POU window function $w_p^{\mathcal{R}, \Gamma}(\mathbf{x})$ for a given pair $(p, \mathcal{R})$ is defined as follows: it is set to zero for all points $\mathbf{x} \notin \Gamma_p^\mathcal{R}$, while for points $\mathbf{x} = \mathcal{M}_p^\mathcal{R}(\xi, \eta) \in \Gamma_p^\mathcal{R}$, the value of $w_p^{\mathcal{R}, \Gamma}(\mathbf{x})$ is prescribed by the relations
\begin{equation}\label{eq:wpgs}
    w_p^{\mathcal{S}, \Gamma}(\mathbf{x}) =
    \begin{cases}
        \widehat{w}^\Gamma \left(\frac{\xi}{\xi_{\ell, p}^\mathcal{S}}\right), \quad &(\xi, \eta) \in \left[0, \xi_{\ell, p}^\mathcal{S}\right] \times \{0\},\\
        1, \quad & (\xi, \eta) \in \left[\xi_{\ell, p}^\mathcal{S}, \xi_{r, p}^\mathcal{S}\right]\times \{0\},\\
        \widehat{w}^\Gamma \left(\frac{\xi-1}{\xi_{r, p}^\mathcal{S}-1}\right), \quad &(\xi, \eta) \in \left[\xi_{r, p}^\mathcal{S}, 1\right]\times \{0\},\\
    \end{cases}
\end{equation}
\begin{equation}\label{eq:wpgc1}
    w_p^{\mathcal{C}_1, \Gamma}(\mathbf{x}) =
    \begin{cases}
        \widehat{w}^\Gamma \left(\frac{\xi}{\xi_{\ell, p}^{\mathcal{C}_1}}\right), \quad &(\xi, \eta) \in \left[0, \xi_{\ell, p}^{\mathcal{C}_1}\right] \times \left\{\frac{1}{2}\right\},\\
        1, \quad & (\xi, \eta) \in \left(\left[\xi_{\ell, p}^{\mathcal{C}_1}, \frac{1}{2}\right]\times \left\{\frac{1}{2}\right\}\right) \cup \left(\left\{\frac{1}{2}\right\}\times \left[\eta_{d, p}^{\mathcal{C}_1}, \frac{1}{2}\right]\right),\\
        \widehat{w}^\Gamma \left(\frac{\eta}{\eta_{d, p}^{\mathcal{C}_1}}\right), \quad &(\xi, \eta) \in\left\{\frac{1}{2}\right\} \times \left[0, \eta_{d, p}^{\mathcal{C}_1}\right], \\
    \end{cases}
\end{equation}
and, letting $\widetilde{\xi}_p^{\mathcal{C}_2} = \max \{\xi_{r, p}^{\mathcal{C}_2}, (d-1)h_{\xi, p}^{\mathcal{C}_2}+\varepsilon, (d-1)h_{\eta, p}^{\mathcal{C}_2}+\varepsilon\}$ and $\widetilde{\eta}_p^{\mathcal{C}_2} = \max \{\eta_{u, p}^{\mathcal{C}_2}, (d-1)h_{\xi, p}^{\mathcal{C}_2} + \varepsilon, (d-1)h_{\eta, p}^{\mathcal{C}_2} + \varepsilon\}$ (where $\varepsilon$ is the same as in~\eqref{eq:Wmepsn}),\\
\begin{equation}\label{eq:wpgc2}
    w_p^{\mathcal{C}_2, \Gamma}(\mathbf{x}) =
    \begin{cases}
        \widehat{w}^\Gamma \left(\frac{\xi-1}{\widetilde{\xi}_p^{\mathcal{C}_2}-1}\right), \quad &(\xi, \eta) \in \left[\widetilde{\xi}_p^{\mathcal{C}_2}, 1\right] \times \left\{0\right\},\\
        1, \quad & (\xi, \eta) \in \left(\left[0, \widetilde{\xi}_p^{\mathcal{C}_2}\right]\times \left\{0\right\}\right) \cup \left(\left\{0\right\}\times \left[0, \widetilde{\eta}_p^{\mathcal{C}_2}\right]\right),\\
        \widehat{w}^\Gamma \left(\frac{\eta-1}{\widetilde{\eta}_p^{\mathcal{C}_2}-1}\right), \quad &(\xi, \eta) \in\left\{0\right\} \times \left[\widetilde{\eta}_p^{\mathcal{C}_2}, 1\right], \\
    \end{cases}
\end{equation}
for $\mathcal{R} = \mathcal{S}$, $\mathcal{C}_1$, and $\mathcal{C}_2$, respectively; see Fig.~\ref{fig:wpRgamma}.
The specific form for the arguments $\widetilde{\xi}_p^{\mathcal{C}_2}$ and $\widetilde{\eta}_p^{\mathcal{C}_2}$ to which $\widehat{w}^\Gamma $ is applied in~\eqref{eq:wpgc2} requires clarification; it ensures the unnormalized 2D window functions~\eqref{wpRhat} in the case $\mathcal{R}=\mathcal{C}_2$ are constant in the square region
\begin{equation}\label{eq:QpR}
    \mathcal{Q}_{p, \varepsilon}^{\mathcal{C}_2} = \mathcal{M}_p^{\mathcal{C}_2}\left(\left[ 0, \max\left\{ (d-1)h_{\xi,p}^{\mathcal{C}_2}+\varepsilon, (d-1)h_{\eta,p}^{\mathcal{C}_2} +\varepsilon\right\} \right]^2\right),
\end{equation}
a property that will be used in the next subsection to guarantee the smoothness of~\eqref{wpR}.
It is easy to check that the pre-POU functions thus defined satisfy the conditions specified by point~(i) in Sec.~\ref{subsec:POUoverall}. 

\subsection{Smoothness of the POU window functions}\label{subsec:smoothPOU}
Once the construction of the one-dimensional pre-POU functions over $\Gamma$ is complete, the full partition of unity~\eqref{eq:POU} is assembled and evaluated at each point $\mathbf{x} \in H^\m$ on the basis of eqs.~\eqref{wpRhat} and~\eqref{wpR}. The present section shows that, in view of the underlying condition~\eqref{eq:epsnquote2}, the POU functions~\eqref{wpR} that result from the constructions in Secs.~\ref{subsec:POUoverall} and~\ref{subsec:prePOU} are indeed smooth in $W^\m$. To do this it suffices to verify that (i)~The un-normalized window function $\widehat{w}_p^\mathcal{R}$ (eq.~\eqref{wpRhat}) is itself smooth for each pair $(p, \mathcal{R})$; and that, (ii)~The sum of all un-normalized window functions satisfies
\begin{equation}\label{whatsum}
    \sum_{\stackrel{\mathcal{R} = \mathcal{S}, \mathcal{C}_1, \mathcal{C}_2}{1 \leq p \leq P_\mathcal{R}}} \widehat{w}_p^{\mathcal{R}}(\mathbf{x}) \neq 0 \quad\textrm{for all} \quad \mathbf{x} \in {W^\m}.
\end{equation}
It is relevant to note that the condition~\eqref{eq:epsnquote2}, which is used below in this section, was included precisely to guarantee the validity of points~(i) and~(ii).

We first check that Point~(ii) is satisfied. Since the functions $\widehat{w}_p^{\mathcal{R}}(\mathbf{x})$ are non-negative, to establish Point~(ii)  it suffices to show that for each $\mathbf{x} \in W^\m$ there exists a pair $(p, \mathcal{R})$ such that $\widehat{w}_p^{\mathcal{R}}(\mathbf{x})\ne 0$. Further, in view of~\eqref{eq:Wmepsn}, it suffices to show that for each $\mathbf{x} \in N_\varepsilon\left(\Omega_p^{\mathcal{R}, \m}; \overline{\Omega}\right)$ we have $\widehat{w}_{p'}^{\mathcal{R}'}(\mathbf{x})\ne 0$ for some pair $(p',\mathcal{R}')$. To this end, we first note that, by construction, the one-dimensional pre-POU function $w_p^{\mathcal{R},\Gamma}$ is positive on the relative open boundary portion~$\Gamma_p^\mathcal{R}$. Moreover, by Rem.~\ref{rmk:ptoGamma},
$\mathcal{P}_p^\mathcal{R}\left(\Intrl{\overline{\Omega}}{V_{p'}^{\mathcal{R}'}}\right)
= \Gamma_{p'}^{\mathcal{R}'},$
and therefore, in view of~\eqref{wpRhat},
\begin{equation}\label{eq:wprne0}
    \textrm{for each $(p', \mathcal{R}')$ we have}\quad\widehat{w}_{p'}^{\mathcal{R}'}(\mathbf{x}) \ne 0 \quad \textrm{for}\quad \mathbf{x} \in \Intrl{\overline{\Omega}}{V_{p'}^{\mathcal{R}'}}.
\end{equation}
Condition~(a) of~\eqref{eq:epsnquote2} guarantees that any point $\mathbf{x}\in N_\varepsilon\left(\Omega_p^\mathcal{R}; \overline{\Omega}\right)$ lies either in the interior $\Intrl{\overline{\Omega}}{V_p^\mathcal{R}}$ of the image of its own parametrization or in the interior $\Intrl{\overline{\Omega}}{V_{p'}^{\mathcal{R}'}}$ of the image associated with a neighboring patch~$\Omega_{p'}^{\mathcal{R}'}$. Consequently, by~\eqref{eq:wprne0}, either $\widehat{w}_p^{\mathcal{R}}(\mathbf{x}) \neq 0$ or $\widehat{w}_{p'}^{\mathcal{R}'}(\mathbf{x}) \neq 0$ necessarily holds, and point~(ii) thus follows, as desired.

To verify Point~(i), on the other hand, we first note that, for each pair $(p,\mathcal{R})$, the function $\widehat{w}_p^\mathcal{R}(\mathbf{x})$ vanishes for $\mathbf{x}\not\in \widehat{\Omega}_p^{\mathcal{R}, \m}$ (eq.~\eqref{wpRhat}), and, thus, it suffices to establish the smoothness of $\widehat{w}_p^\mathcal{R}(\mathbf{x})$ for $\mathbf{x}$ in the set $\Clrl{\overline{\Omega}}{\widehat{\Omega}_p^{\mathcal{R}, \m}} \cap W^\m = \left(\widehat{\Omega}_p^{\mathcal{R}, \m} \cup \Bndrl{\overline{\Omega}}{\widehat{\Omega}_p^{\mathcal{R}, \m}}\right) \cap W^\m$. We first do this for $\mathbf{x}\in \widehat{\Omega}_p^{\mathcal{R}, \m} \subset \Intrl{\overline{\Omega}}{V_p^\mathcal{R}}$. Since, according to~\eqref{wpRhat}, $\widehat{w}_p^\mathcal{R}$ is defined as the composition of the smooth 1D function $w_p^{\mathcal{R},\Gamma}$ with the boundary projection $\mathcal{P}_p^\mathcal{R}$, the smoothness of $\widehat{w}_p^\mathcal{R}$ depends in part on the smoothness/non-smoothness properties of $\mathcal{P}_p^\mathcal{R}$---which we consider in what follows for each one of the three patch types.  

In the case $\mathcal{R} = \mathcal{S}$, the boundary projection $\mathcal{P}_p^\mathcal{S}$ is given by the composition \[ \mathcal{P}_p^\mathcal{S} = \mathcal{M}_p^{\mathcal{S}} \circ \widetilde{\mathcal{P}}_p^\mathcal{S} \circ \left(\mathcal{M}_p^\mathcal{S}\right)^{-1}, \] as stated in Rem.~\ref{rmk:ptoGamma}. But, by Rem.~\ref{rmk:MpRsmooth}, both $\mathcal{M}_p^\mathcal{S}(\mathbf{x})$ and its inverse are smooth functions, and by~\eqref{pS}, so is  $\widetilde{\mathcal{P}}_p^\mathcal{S}$---and, thus, $\widehat{w}_p^\mathcal{S}$ is smooth in $\widehat{\Omega}_p^\mathcal{S} \subset V_p^\mathcal{R}$. In the case $\mathcal{R} = \mathcal{C}_1$, a similar argument involving forward and inverse mappings shows that the projection $\mathcal{P}_p^{\mathcal{C}_1}$ is smooth in $\widehat{\Omega}_p^{\mathcal{C}_1}$ except at the image under $\mathcal{M}_p^{\mathcal{C}_1}$ of the sets
\begin{equation}\label{pC1sm}
    \left\{\left(\frac{1}{2},\eta\right):\frac{1}{2} \leq\eta\leq 1\right\}, \quad \textrm{and}\quad \left\{\left(\xi,\frac{1}{2}\right):\frac{1}{2} \leq\xi\leq 1\right\},
\end{equation}
but in view of eq.~\eqref{eq:wpgc1} and~\eqref{pC1}, $\widehat{w}_p^{\mathcal{C}_1}$ is clearly constant in a neighborhood of these image sets, and, thus, the smoothness of this function throughout $\widehat{\Omega}_p^{\mathcal{C}_1}$ follows. In the case $\mathcal{R} = \mathcal{C}_2$, finally, an argument based on the forward and inverse mappings similar to the ones above shows that  $\widetilde{\mathcal{P}}_p^{\mathcal{C}_2}$ is smooth except along the discontinuity curve~\eqref{quote3}. But condition~(b) in~\eqref{eq:epsnquote2} ensures that the $\varepsilon$-neighborhood $N_\varepsilon\left(\Omega_{p'}^{\mathcal{R}', \m}; \overline{\Omega}\right)$ of any neighboring patchwise matching domain $\Omega_{p'}^{\mathcal{R}', \m}$ does not intersect this curve, and thus, in view of~\eqref{eq:omegapRhat} and the description following that equation we have
\begin{equation}
    \widehat{\Omega}_{p}^{\mathcal{C}_2, \m} \cap \{(\xi, \eta): \xi=\eta\}
    \subseteq 
    N_\varepsilon\left(\Omega_p^{\mathcal{C}_2, \m}; \overline{\Omega}\right) \cap \{(\xi, \eta): \xi=\eta\}.
\end{equation}
In view of definitions~\eqref{OmegaC2m} and~\eqref{eq:QpR}, further,
\begin{equation}
    N_\varepsilon\left(\Omega_p^{\mathcal{C}_2, \m}; \overline{\Omega}\right)  \cap \{(\xi, \eta): \xi=\eta\}\subseteq 
    \mathcal{Q}_{p, \varepsilon}^{\mathcal{C}_2}.
\end{equation}
Since $\widehat{w}_p^{\mathcal{C}_2}$ is constant on $\mathcal{Q}_{p, \varepsilon}^{\mathcal{C}_2}$ (see~\eqref{eq:QpR} and accompanying discussion), it follows that $\widehat{w}_p^{\mathcal{C}_2}$ is smooth on $\widehat{\Omega}_p^{\mathcal{C}_2}$. We have thus shown that, as claimed, $\widehat{w}_p^{\mathcal{R}}$ is smooth in $\widehat{\Omega}_p^{\mathcal{R}}$ for all three patch types, $\mathcal{R} = \mathcal{S}$,  $\mathcal{C}_1$, and $\mathcal{C}_2$.

It remains for us to establish the smoothness of $\widehat{w}_p^{\mathcal{R}}$ on the intersection of $W^\m$ with $\Bndrl{\overline{\Omega}}{\widehat{\Omega}_p^{\mathcal{R}, \m}}$. To do this, taking into account~\eqref{eq:EpnuR}--\eqref{eq:EtpnuC2}, we first show that 
\begin{equation}\label{eq:UmeqEpnu}
    \Bndrl{\overline{\Omega}}{\widehat{\Omega}_p^{\mathcal{R}, \m}} \cap W^\m \subset E_{p,\nu}^\mathcal{R} \cap W^\m
\end{equation}
Indeed, by~\eqref{eq:omegapRhat}, $\Bndrl{\overline{\Omega}}{\widehat{\Omega}_p^{\mathcal{R}, \m}} = \Bndrl{\overline{\Omega}}{V_p^\mathcal{R}\cap W^\m}$. But, for any sets $A$ and $B$ we have $\partial (A \cap B) \subset \partial A \cup \partial B$, and thus, $\Bndrl{\overline{\Omega}}{\widehat{\Omega}_p^{\mathcal{R}, \m}} \subset  \Bndrl{\overline{\Omega}}{V_p^\mathcal{R}}\cup \Bndrl{\overline{\Omega}}{W^\m}$. In view of Rem.~\ref{rmk:intedges}, it follows that 
\begin{equation}
\begin{aligned}
    \Bndrl{\overline{\Omega}}{\widehat{\Omega}_p^{\mathcal{R}, \m}} \cap W^\m\subset  \Bndrl{\overline{\Omega}}{V_p^\mathcal{R}} \cap W^\m &= \left(E_{p,\tau}^\mathcal{R} \cup E_{p,\nu}^\mathcal{R}\right) \cap W^\m \\&= \left(E_{p,\tau}^\mathcal{R} \cap W^\m\right) \cup \left(E_{p,\nu}^\mathcal{R} \cap W^\m\right)
\end{aligned}
\end{equation}
where the second to last equality follows from~\eqref{eq:VpRbound}. By part (b) of condition~\eqref{eq:epsnquote2}, the neighborhood 
$N_{\varepsilon}\left(\Omega_p^{\mathcal{R},\m};\overline{\Omega}\right)$ and each neighborhood
$N_{\varepsilon}\left(\Omega_{p'}^{\mathcal{R}',\m};\overline{\Omega}\right)$ 
associated with patches $\Omega_{p'}^{\mathcal{R}'}$ neighboring $\Omega_{p}^{\mathcal{R}}$ are all
disjoint from the tangential interior edges $E_{p,\tau}^\mathcal{R}$ of $V_p^\mathcal{R}$. 
It thus follows that $E_{p,\tau}^\mathcal{R} \cap W^\m = \emptyset$ and~\eqref{eq:UmeqEpnu} results directly. To complete the proof we establish the smoothness of the function $\widehat{w}_p^\mathcal{R}$ at all points in $E_{p,\nu}^\mathcal{R} \cap W^\m$. Clearly, considering once again~\eqref{wpRhat} and Rems.~\ref{rmk:intedges} and~\ref{rmk:ptoGamma} it follows that the function $\widehat{w}_p^\mathcal{R}(\mathbf{x})$ smoothly approaches zero as $\mathbf{x}$ approaches $E_{p,\nu}^\mathcal{R} \cap W^\m$, since $w_p^{\mathcal{R}, \Gamma}$ smoothly vanishes at the projection points $\mathcal{P}_p^\mathcal{R}(E_{p,\nu}^\mathcal{R})$. Since $\widehat{w}_p^{\mathcal{R}}$ is zero outside $\widehat{\Omega}_p^{\mathcal{R}}$, this completes the proof.

\section{Patchwise extension domains, meshes, and blending-to-zero}\label{SM3}
\subsection{Patchwise extension domains and meshes}\label{subsec:emeshes}
For all pairs $(p, \mathcal{R})$ in~\eqref{exteriordomain}, the extension domains $\Omega_p^{\mathcal{R}, \e}\subset V_p^\mathcal{R}$ are defined as the images under the parametrization $\mathcal{M}_p^\mathcal{R}$ of the corresponding parameter space domains
\begin{equation}\label{eq:omegatildepSe}
    \widetilde{\Omega}^{\mathcal{S}, \e}_p = (0, 1)\times(-Ch_{\eta, p}^\mathcal{S}, 0] \quad \mbox{for}\quad \mathcal{R}=\mathcal{S},
\end{equation}
\begin{equation}\label{eq:omegatildepC1e}
    \widetilde{\Omega}^{\mathcal{C}_1, \e}_p = \left(\left(0, \frac{1}{2}\right]\times\left(\frac{1}{2}-Ch_{\eta, p}^{\mathcal{C}_1}, \frac{1}{2}\right]\right)\cup\left(\left(\frac{1}{2}-Ch_{\xi, p}^{\mathcal{C}_1}, \frac{1}{2}\right]\times\left(0, \frac{1}{2}\right]\right)\quad \mbox{for}\quad\mathcal{R}=\mathcal{C}_1,
\end{equation}
and
\begin{equation}\label{eq:omegatildepC2e}
    \widetilde{\Omega}^{\mathcal{C}_2, \e}_p = \left(\left[0,1\right)\times\left(-Ch_{\eta, p}^{\mathcal{C}_2}, 0\right]\right)\cup\left(\left(-Ch_{\xi, p}^{\mathcal{C}_2}, 0\right]\times\left[0, 1\right)\right)\quad \mbox{for}\quad\mathcal{R}=\mathcal{C}_2.
\end{equation}
which are solely determined by the mesh sizes $h_{\xi, p}^\mathcal{R}$ and $h_{\eta, p}^\mathcal{R}$ of the patchwise matching mesh~\eqref{eq:HpRm} and the BTZ parameter $C$ introduced in Sec.~\ref{sec:fc_1d}. In view of~\eqref{eq:param}, here it is assumed that the previously defined discretization parameters $h_{\xi, p}^\mathcal{R}$ and $h_{\eta, p}^\mathcal{R}$ are sufficiently small, so that 
\begin{equation}
    \widetilde{\Omega}_p^{\mathcal{R}, \e} \subset \widetilde{V}_p^\mathcal{R}.
    \label{eq:omegasubV}
\end{equation}

Having constructed the patchwise extension domains $\Omega_p^{\mathcal{R}, \e}$, the patchwise extension mesh $H^{\mathcal{R},\e}_{p, \rho}$ in~\eqref{exteriormesh} defined in what follows provides a discretization of the closure $\Cl{\Omega^{\mathcal{R}, \e}_p}$. In view of~\eqref{eq:gammaRp_param} and~\eqref{omegaSm}--\eqref{OmegaC2m}, the 2D extension values are obtained by an iterated application of the leftward $\ell$-BTZ procedure (Rem.~\ref{rmk:l-btz}) along coordinate curves in parameter space $\widetilde{V}_p^\mathcal{R}$ on the basis of known function values in the matching set $H_p^{\mathcal{R}, \m}$. 
The corresponding patchwise extension mesh $H_p^{\mathcal{R}, \e}$ is then defined as the image under the parametrization $\mathcal{M}_p^\mathcal{R}$ of equispaced parameter space meshes, each one of which is defined on the basis of the patchwise matching mesh~\eqref{eq:HpRm} as well as the BTZ parameters $C$, $\rho$, and $C_\rho$ ($C_\rho = \rho C$) introduced in Sec.~\ref{sec:fc_1d}. 
To provide the details of these constructions  we let
\begin{equation*}
\begin{aligned}
    &\mathcal{H}(\xi_{\ell}, \xi_{r}, \eta_{d}, \eta_u, n_{\xi}, n_{\eta})\\
    &\qquad=  \left\{\left(i\frac{\xi_r-\xi_\ell}{n_\xi-1}+\xi_\ell, j \frac{\eta_u-\eta_d}{n_\eta-1}+\eta_d\right):\ 0\leq i \leq n_{\xi}-1, 0\leq j \leq n_{\eta}-1\right\}
\end{aligned}
\end{equation*}
denote the uniform grid defined on the rectangle $[\xi_\ell, \xi_r]\times[\eta_d, \eta_u]$ with $n_\xi$ nodes in the $\xi$ direction and $n_\eta$ nodes in the $\eta$ direction, and we define the parameter-space extension meshes $\widetilde{H}_{p, \rho}^{\mathcal{R}, \e}$; once such meshes are introduced, which we do below in what follows, the physical-space extension meshes are given by
\begin{equation}\label{eq:HpRe}
    H_{p, \rho}^{\mathcal{R}, \e} = \mathcal{M}_p^\mathcal{R}(\widetilde{H}_{p, \rho}^{\mathcal{R}, \e}) \subset \Omega_p^{\mathcal{R}, \e}.
\end{equation}

\begin{figure}[phtb]
\centering
\includegraphics[width=0.55\linewidth,]{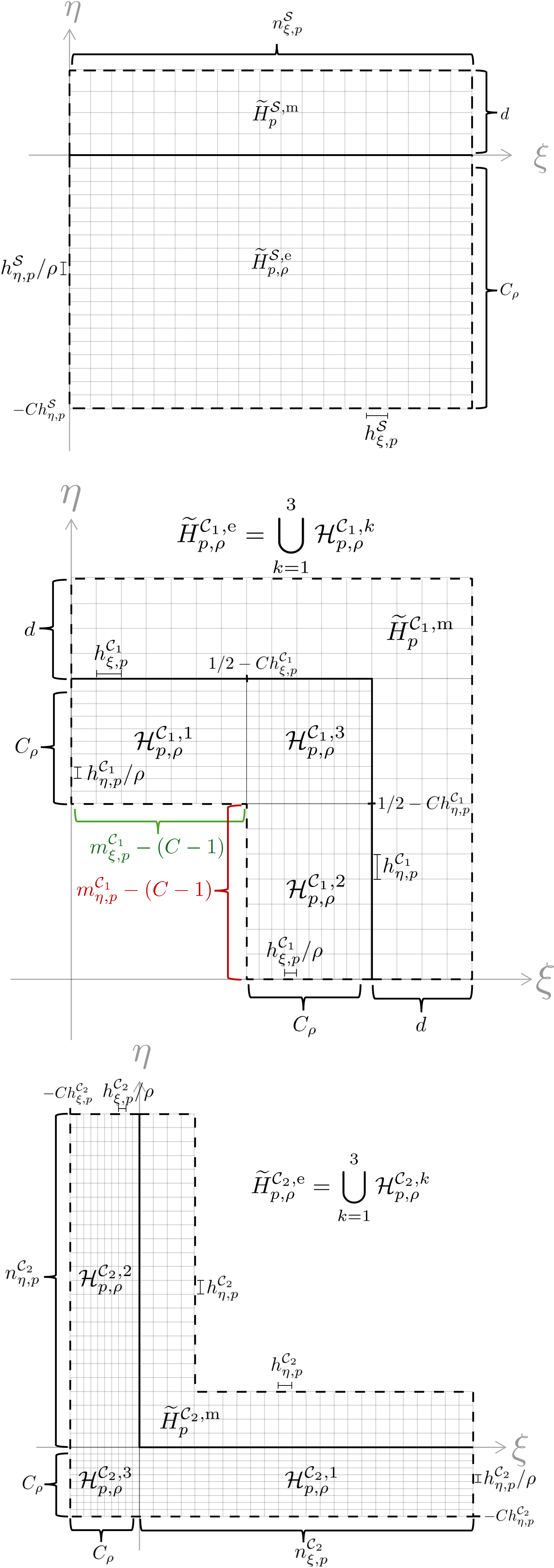}
\caption{
Parameter-space patchwise matching meshes $\widetilde{H}_p^{\mathcal{R}, \m}$ and extension meshes $\widetilde{H}_p^{\mathcal{R}, \e}$ for $\mathcal{R} = \mathcal{S}$ (top panel), $\mathcal{R} = \mathcal{C}_1$ (middle panel), and $\mathcal{R} = \mathcal{C}_2$ (bottom panel); cf. Fig.~\ref{fig:Omegae}. In each case, mesh lines (gray) are superimposed on the corresponding matching domains $\widetilde{\Omega}_p^{\mathcal{R}, \m}$ and extension domains $\widetilde{\Omega}_p^{\mathcal{R}, \e}$ (not labeled to avoid visual clutter). The solid black lines represent the interfaces between $\widetilde{\Omega}_p^{\mathcal{R}, \e}$ and $\widetilde{\Omega}_p^{\mathcal{R}, \m}$.
}
    \label{fig:OmegapRe}
\end{figure}

The parameter-space extension mesh $\widetilde{H}_{p, \rho}^{\mathcal{R}, \e}$ for the patch type $\mathcal{R} = \mathcal{S}$ is given by
\begin{equation}\label{BTZ_patch_mesh_S}
    \widetilde{H}_{p, \rho}^{\mathcal{S}, \e} = \mathcal{H}(0, 1, -Ch_{\eta, p}^\mathcal{S}, 0, n_{\xi, p}^\mathcal{S}, C_\rho+1).
\end{equation}
The corresponding meshes for the patch types $\mathcal{R} = \mathcal{C}_1$ and $\mathcal{C}_2$, in turn, are given by
\begin{equation}\label{BTZ_patch_mesh_C}
    \widetilde{H}_{p, \rho}^{\mathcal{R}, \e} = \bigcup_{k=1}^3 \mathcal{H}_{p, \rho}^{\mathcal{R}, k},
\end{equation}
where, calling $m_{\xi, p}^{\mathcal{C}_1}=\frac{n_{\xi, p}^{\mathcal{C}_1}-1}{2}$ and $ m_{\eta, p}^{\mathcal{C}_1} = \frac{n_{\eta, p}^{\mathcal{C}_1}-1}{2}$, the three ($\rho$-dependent) component grids $\mathcal{H}_{p, \rho}^{\mathcal{R}, k}$ are given by 
\begin{equation}\label{eq:H1pC1}
    \mathcal{H}_{p, \rho}^{\mathcal{C}_1, 1} = \mathcal{H}\left(0, \frac{1}{2}-Ch_{\xi, p}^{\mathcal{C}_1}, -Ch_{\eta, p}^{\mathcal{C}_1}, 0, m_{\xi, p}^{\mathcal{C}_1}-(C-1), C_\rho+1\right),
\end{equation}
\begin{equation}\label{eq:H2pC1}
    \mathcal{H}_{p, \rho}^{\mathcal{C}_1, 2} = \mathcal{H}\left(-Ch_{\xi, p}^{\mathcal{C}_1}, 0, 0, \frac{1}{2}-Ch_{\eta, p}^{\mathcal{C}_1}, C_\rho+1, m_{\eta, p}^{\mathcal{C}_1}-(C-1)\right),\quad\textrm{and}
\end{equation}
\begin{equation}\label{eq:H3pC1}
    \mathcal{H}_{p, \rho}^{\mathcal{C}_1, 3} = \mathcal{H}\left(\frac{1}{2}-Ch_{\xi, p}^{\mathcal{C}_1}, \frac{1}{2},  \frac{1}{2}-Ch_{\eta, p}^{\mathcal{C}_1}, \frac{1}{2},C_\rho+1, C_\rho+1\right)
\end{equation}
 in the case $\mathcal{R} = \mathcal{C}_1$, and they are given  by
\begin{equation}\label{eq:H1pC2}
    \mathcal{H}_{p, \rho}^{\mathcal{C}_2, 1} = \mathcal{H}\left(0, 1, -Ch_{\eta, p}^{\mathcal{C}_2}, 0,n_{\xi, p}^{\mathcal{C}_2}, C_\rho+1\right),
\end{equation}
\begin{equation}\label{eq:H2pC2}
    \mathcal{H}_{p, \rho}^{\mathcal{C}_2, 2} = \mathcal{H}\left(-Ch_{\xi, p}^{\mathcal{C}_2}, 0, 0, 1, C_\rho+1, n_{\eta, p}^{\mathcal{C}_2}\right),\quad\textrm{and}
\end{equation}
\begin{equation}\label{eq:H3pC2}
    \mathcal{H}_{p, \rho}^{\mathcal{C}_2, 3} = \mathcal{H}\left(-Ch_{\xi, p}^{\mathcal{C}_2}, 0, -Ch_{\eta, p}^{\mathcal{C}_2}, 0, C_\rho+1, C_\rho+1\right)
\end{equation}
in the case $\mathcal{R}=\mathcal{C}_2$. Fig.~\ref{fig:OmegapRe} displays the meshes $H_{p, \rho}^{\mathcal{R}, \e}$ overlaid on $\Omega_p^{\mathcal{R}, \e}$ in the cases $\mathcal{R} = \mathcal{S}, \mathcal{C}_1,$ and $\mathcal{C}_2$.

Given the constructions above and the partition-of-unity-weighted functions $f_p^\mathcal{R}: H^{\m}\rightarrow \mathbb{C}$ computed in Sec.~\ref{subsec:POU}, and more specifically, the values of $f_p^\mathcal{R}$ at points in the patchwise matching mesh $H_p^{\mathcal{R}, \m}$, we can obtain the patchwise BTZ grid functions $f_p^{\mathcal{R}, \bt}: H^{\e} \rightarrow \mathbb{C}$. The specific procedures for constructing $f_p^{\mathcal{R}, \bt}$ in the cases $\mathcal{R} = \mathcal{S}, \mathcal{C}_1,$ and $\mathcal{C}_2$ are presented in detail in Secs.~\ref{subsec:SpatchBTZ}--\ref{subsec:C2patchBTZ} below. To facilitate the discussion, we introduce interval notation for subsets of the matching grids $\widetilde{H}_p^{\mathcal{R}, \m}$ in the $\xi$- and $\eta$-directions, following the conventions of eq.~\eqref{eq:Hmtilde} and its accompanying text:  
 \begin{equation}\label{eq:xivec}
     \xivec{p}{\mathcal{R}}(a, b) = \left\{\xi_i = ih_{\xi, p}^{\mathcal{R}}: i\in \llbracket a, b\rrbracket\right\}
 \end{equation}
and
 \begin{equation}\label{eq:etavec}
     \etavec{p}{\mathcal{R}}(a, b)  = \left\{ \eta_i = ih_{\eta, p}^{\mathcal{R}}: i\in \llbracket a, b\rrbracket\right\}.
 \end{equation}

 \begin{remark}\label{rmk:ftildepR}
In view of~\eqref{wpR}, the windowed function $f_p^\mathcal{R}$ vanishes outside $V_p^\mathcal{R} \cap H^\m$, and consequently the 2D-BTZ functions also vanish in $H^\e \setminus H_p^{\mathcal{R}, \e}$. It therefore suffices to construct $f_p^{\mathcal{R}, \bt}$ on $H_p^{\mathcal{R}, \e} \subset V_p^\mathcal{R}$, which is achieved on the basis of the restriction of the grid function $f_p^\mathcal{R}$ to $H_p^{\mathcal{R}, \m}$. For clarity, it is also convenient to formulate the 2D-BTZ construction in parameter space. To this end, we introduce the parameter-space functions $\widetilde{f}_p^{\mathcal{R}} : \widetilde{H}_p^{\mathcal{R}, \m} \to \mathbb{C}$ and $\widetilde{f}_p^{\mathcal{R}, \bt} : \widetilde{H}_p^{\mathcal{R}, \e} \to \mathbb{C}$ defined by
\begin{equation}
    \widetilde{f}_p^{\mathcal{R}}(\xi, \eta) = f_p^{\mathcal{R}}(\mathcal{M}_p^\mathcal{R}(\xi, \eta)), 
    \qquad
    \widetilde{f}_p^{\mathcal{R}, \bt}(\xi, \eta) = f_p^{\mathcal{R}, \bt}(\mathcal{M}_p^\mathcal{R}(\xi, \eta)).
\end{equation}
Clearly, the BTZ extensions  $f_p^{\mathcal{R}, \bt}$ are obtained in terms of  $\widetilde{f}_p^{\mathcal{R}, \bt}$ by employing the inverse parametrization $\left(\mathcal{M}_p^\mathcal{R}\right)^{-1}$.
\end{remark}

\begin{figure}[tb]
    \centering
    \includegraphics[width=0.4\linewidth]{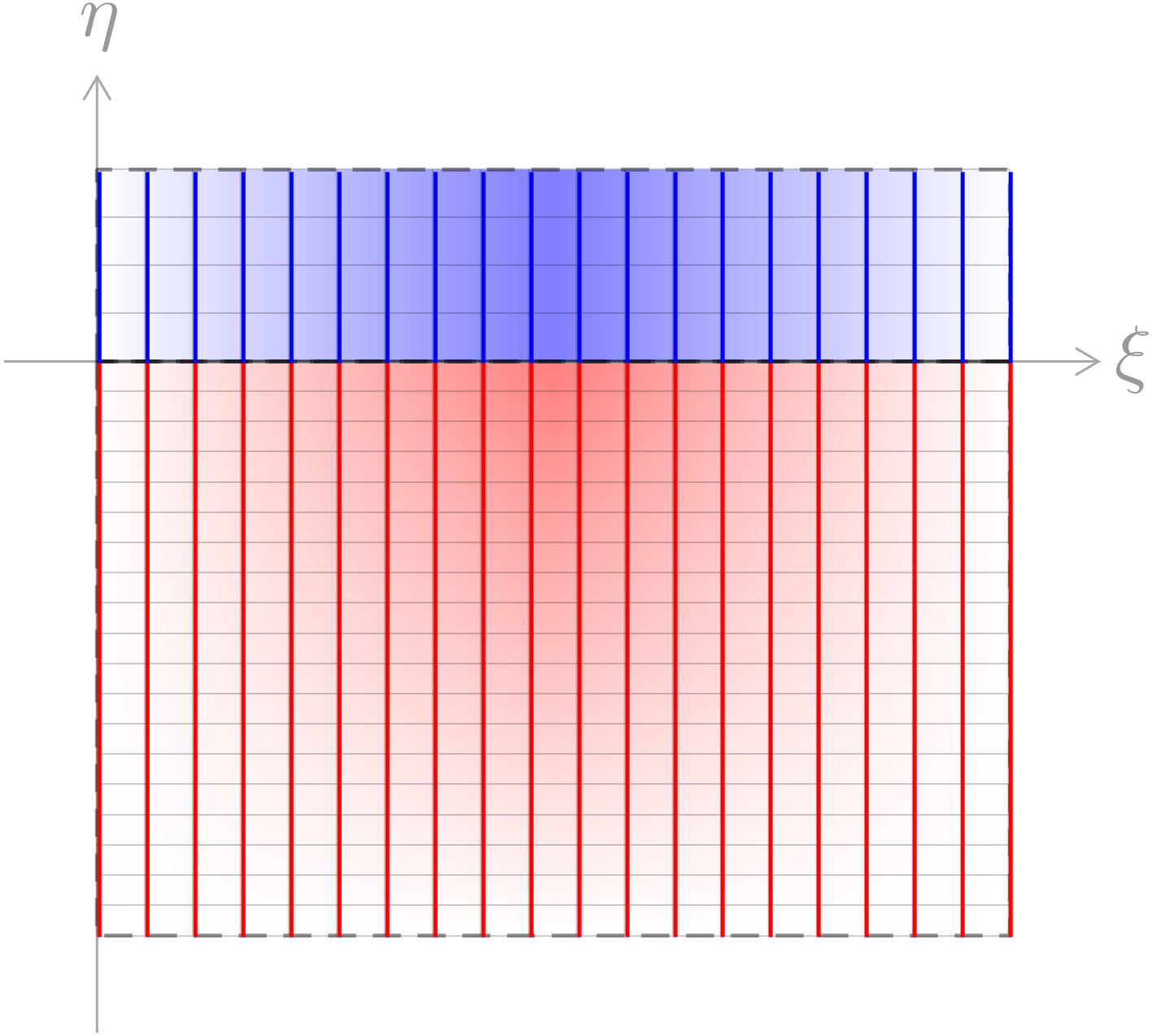}
    \caption{Parameter-space BTZ construction for $\mathcal{R}=\mathcal{S}$; cf. Fig.~\ref{fig:SBTZprelim}. The blue and red lines represent the parameter-space matching and extension meshes $\widetilde{H}_p^{\mathcal{S}, \m}$ and $\widetilde{H}_{p, \rho}^{\mathcal{S}, \e}$, respectively. Color shading is used to indicate the smooth blending-to-zero induced by the windowing functions $w_p^\mathcal{S}$ along the three relevant boundaries of the extension domain $\Omega_p^{\mathcal{S},\e}$ and the two relevant boundaries of the matching domain $\Omega_p^{\mathcal{S},\m}$.}
    \label{fig:SBTZdetail}
\end{figure}

\subsection{$\mathcal{S}$-type patch blending to zero in parameter space}\label{subsec:SpatchBTZ}

Given the discrete windowed grid function $\widetilde{f}_p^\mathcal{S}$ (Rem.~\ref{rmk:ftildepR} with $\mathcal{R} = \mathcal{S}$) defined on the patchwise matching mesh $\widetilde{H}_p^{\mathcal{S}, \m} \subset \widetilde{\Omega}_p^{\mathcal{S}}$ of a given $\mathcal{S}$-type patch $\Omega_p^{\mathcal{S}}$, the desired blended-to-zero grid function $\widetilde{f}_p^{\mathcal{S}, \textrm{b}}$ is defined on the set $\widetilde{H}_{p, \rho}^{\mathcal{S}, \e}$ (eq.~\eqref{BTZ_patch_mesh_S}) by employing 1D $\ell$-BTZ grid function values along the $\eta$ direction of the mesh $\widetilde{H}_p^{\mathcal{S}, \m}$. In detail,  with reference to Rem.~\ref{rmk:l-btz} and eq.~\eqref{eq:xivec}, and setting $D^{\m} = \{x_0, x_1, \dots, x_{d-1}\} = \etavec{p}{\mathcal{S}}(0, d-1)$, we consider the vector
\begin{equation}
    \Phi^{\ell} = \Phi_i^{\ell} = (\widetilde{f}_p^\mathcal{S}(\xi_i, x_0), \dots, \widetilde{f}_p^\mathcal{S}(\xi_i, x_{d-1}))^t, \quad 0 \leq i \leq n_{\xi, p}^\mathcal{S} - 1
\end{equation}
together with the corresponding vector $\Phi_{C_\rho,i}^{\ell b}$ of blending-to-zero grid function values and points $D^{\e}_{C_\rho} = \{\widetilde{x}_{-\rho C}, \widetilde{x}_{-\rho C + 1}, \dots, \widetilde{x}_0\}$ obtained per Rem.~\ref{rmk:refinebtz}, and then, noting that 
\begin{equation}
    \widetilde{H}_p^{\mathcal{S}, \e} = \xivec{p}{\mathcal{S}}(0, n_{\xi, p}^\mathcal{S}-1) \times D^{\e}_{C_\rho}
\end{equation}
we define
\begin{equation}
    \widetilde{f}_p^{\mathcal{S}, \bt}(\xi_i,\widetilde{x}_j)  = \Phi_{C_
    \rho,i}^{\ell b} (\widetilde{x}_j)\quad\mbox{for}\quad (\xi_i,
    \widetilde{x}_j)\in \widetilde{H}_{p, \rho}^{\mathcal{S}, \e}.
\end{equation}

\begin{figure}[htb]
    \centering
    \includegraphics[width=\linewidth]{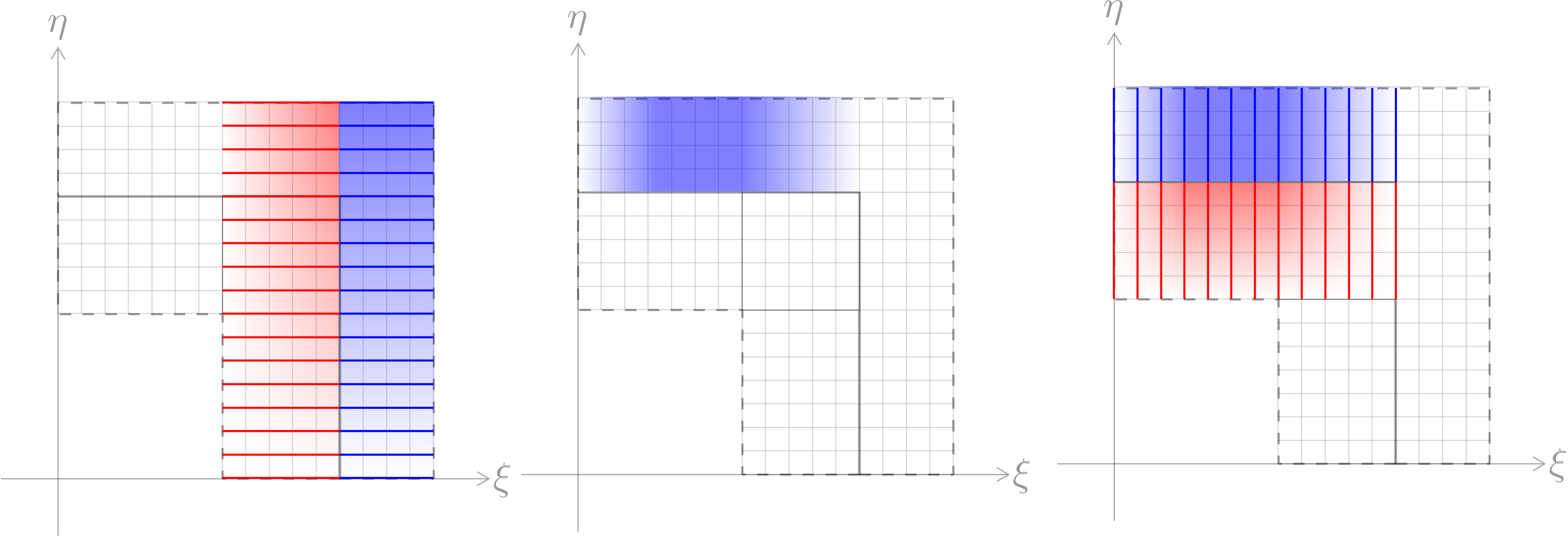}
    \caption{Parameter-space BTZ construction for $\mathcal{R}=\mathcal{C}_1$ (cf. Fig.~\ref{fig:C1BTZprelim}) with refinement factor $\rho=1$ (cf. Fig.~\ref{fig:C1BTZrho}). The blue and red lines represent the matching and extension meshes employed in steps 1, 2, and 3 of Sec.~\ref{subsec:C1patchBTZ}. The combination of the left- and right-panel contributions, which reproduces the original windowed function $\widetilde{f}_p^{\mathcal{C}_1}$ on the matching domain $\Omega_p^{\mathcal{C}_1,\m}$, smoothly blends to zero along the relevant boundaries of the matching and extension domains $\Omega_p^{\mathcal{C}_1,\m}$ and $\Omega_p^{\mathcal{C}_1,\e}$.}
    \label{fig:C1BTZdetail}
\end{figure}
\subsection{$\mathcal{C}_1$-type patch blending to zero in parameter space}\label{subsec:C1patchBTZ}
The blended-to-zero functions for $\mathcal{C}_1$- and $\mathcal{C}_2$-type patches are obtained through applications of the 1D  $\ell$-BTZ procedure along both the $\xi$ and $\eta$ directions. To distinguish between these applications, (and with reference to Rem.~\ref{rmk:l-btz} and~\ref{rmk:refinebtz}), we adopt direction-specific notation:
\begin{itemize}
  \item[--] For blending in the $\xi$ direction, we denote: the matching grid points by $D^{\m, \xi} = \{x_0^\xi, x_1^\xi, \dots, x_{d-1}^\xi\}$; the matching values by $\Phi^{\ell, \xi}$; the refined extension grid points by $D^{\e, \xi}_{C_\rho} = \left\{\widetilde{x}_{-\rho C}^\xi, \widetilde{x}_{-\rho C+1}^\xi, \dots, \widetilde{x}_0^\xi\right\}$; and the blended-to-zero vector by $\Phi_{C_\rho}^{\ell b, \xi}$; and, similarly,

  \item[--] For blending in the $\eta$ direction, we denote: the matching grid points by $D^{\m, \eta} = \{x_0^\eta, x_1^\eta, \dots, x_{d-1}^\eta\}$; the matching values by $\Phi^{\ell, \eta}$; the refined extension grid points by $D^{\e, \eta}_{C_\rho} = \left\{\widetilde{x}_{-\rho C}^\eta, \widetilde{x}_{-\rho C+1}^\eta, \dots, \widetilde{x}_0^\eta\right\}$; and the blended-to-zero vector by $\Phi_{C_\rho}^{\ell b, \eta}$.
\end{itemize}
This notation allows us to easily indicate, within a unified construction setup, the iterated application of the $\ell$-BTZ procedure, in either of the two coordinate directions, as needed for the evaluation of $\mathcal{C}_1$- and $\mathcal{C}_2$-type blending to zero.

The $\mathcal{C}_1$-type blending to zero procedure considered in this section, in particular, proceeds as follows; for clarity the case $\rho = 1$, for which the mesh sizes of the 1D-BTZ matching and extension meshes coincide, is considered first (see Rem.~\ref{rmk:refinebtz}); the general case $\rho > 1$ then proceeds as outlined in Rem.~\ref{rmk:C1nrho}. Given the windowed grid function $\widetilde{f}_p^{\mathcal{C}_1}$ (eq.~\eqref{eq:fpR} with $\mathcal{R} = \mathcal{C}_1$) defined on the parameter-space patchwise-matching mesh $\widetilde{H}_p^{\mathcal{C}_1, \m} \subset \widetilde{\Omega}_p^{\mathcal{C}_1}$ of a $\mathcal{C}_1$-type patch $\Omega_p^{\mathcal{C}_1}$, in accordance with the notation of eq.~\eqref{eq:IpC1} and associated discussion, we introduce the grid point sets
\begin{equation}\label{eq:Atilde}
    \widetilde{\mathcal{A}}_{p}^{\mathcal{C}_1, \m} = 
    \big\{(\xi_i, \eta_j): (i, j) \in 
    \llbracket m_{\xi, p}^{\mathcal{C}_1}, m_{\xi, p}^{\mathcal{C}_1} + d-1\rrbracket\times 
    \llbracket0, m_{\eta, p}^{\mathcal{C}_1} + d-1\rrbracket
    \big\} \subset \widetilde{H}_p^{\mathcal{C}_1, \m},
\end{equation}
and 
\begin{equation}
    \widetilde{\mathcal{B}}_{p}^{\mathcal{C}_1, \m} = 
    \big\{(\xi_i, \eta_j): (i, j) \in 
    \llbracket 0, m_{\xi, p}^{\mathcal{C}_1} \rrbracket \times 
    \llbracket m_{\eta, p}^{\mathcal{C}_1}, m_{\eta, p}^{\mathcal{C}_1} + d-1\rrbracket
    \big\} \subset \widetilde{H}_p^{\mathcal{C}_1, \m},
\end{equation}
that provide the needed discretizations of the first and second sets on the right-hand side of eq.~\eqref{OmegaC1m}.
On the basis of these discretizations the overall BTZ operation is produced by means of the following four step procedure.
\paragraph{Four-step $\mathcal{C}_1$-type BTZ procedure:}
\begin{enumerate}
\item Blending $\widetilde{f}_p^{\mathcal{C}_1}$ along the $\xi$-direction of $\widetilde{\mathcal{A}}_{p}^{\mathcal{C}_1, \m}$ to produce a smooth BTZ function $\widehat{g}_p^{\mathcal{C}_1, \bt}: \widetilde{H}_p^{\mathcal{C}_1, \m} \cup \widetilde{H}_{p, 1}^{\mathcal{C}_1, \e} \to \mathbb{C}$ satisfying $\widehat{g}_p^{\mathcal{C}_1}(\xi, \eta) = \widetilde{f}_p^{\mathcal{C}_1}(\xi, \eta)$ for $(\xi, \eta) \in \widetilde{\mathcal{A}}_{p}^{\mathcal{C}_1, \m}$ and $\widehat{g}_p^{\mathcal{C}_1}(\xi, \eta) = 0$ for $(\xi, \eta) \notin \widetilde{\mathcal{A}}_{p}^{\mathcal{C}_1, \m} \cup \left(D^{\e, \xi}_{C} \times \etavec{p}{\mathcal{C}_1}(0, m_{\eta, p}^{\mathcal{C}_1}+d-1)\right)$.

\item Subtracting the blended-to-zero values obtained per step 1 from the values of $\widetilde{f}_p^\mathcal{R}$ on $\widetilde{H}_p^{\mathcal{C}_1, \m}$ to obtain the discrete function $g_p^{\mathcal{C}_1}: \widetilde{H}_p^{\mathcal{C}_1, \m}\rightarrow \mathbb{C}$.

\item Blending $g_p^{\mathcal{C}_1}$ along the $\eta$-direction of $\widetilde{\mathcal{B}}_{p}^{\mathcal{C}_1, \m}$ and denoting the resulting BTZ function by $g_p^{\mathcal{C}_1, \bt}$.

\item Defining the final $\mathcal{C}_1$-type BTZ function by $\widetilde{f}_p^{\mathcal{C}_1, \bt}$ as the sum of $\widehat{g}_p^{\mathcal{C}_1, \bt}$ and ${g}_p^{\mathcal{C}_1}$ on $\widetilde{H}_{p, 1}^{\mathcal{C}_1, \e}$.

\end{enumerate}
Additional details concerning each one of these algorithmic steps are provided in what follows.

Step~1 utilizes the function values
\begin{equation}
    \Phi^{\ell, \xi}_i = \big(\widetilde{f}_p^{\mathcal{C}_1}(x_0^\xi, \eta_i), \dots, \widetilde{f}_p^{\mathcal{C}_1}(x_{d-1}^\xi, \eta_i)\big)^t,\quad 0 \leq i \leq m_{\eta, p}^{\mathcal{C}_1}+(d-1), 
\end{equation}
of the grid function $\widetilde{f}_p^{\mathcal{C}_1}$ at the 1D matching mesh
$D^{\m, \xi} = \xivec{p}{\mathcal{C}_1}(m_{\xi, p}^{\mathcal{C}_1}, m_{\xi, p}^{\mathcal{C}_1}+d-1 )$ (cf. \eqref{eq:Hmtilde} and~\eqref{eq:IpC1}) to  produce, for each fixed $\eta_i$, the $\ell$-BTZ extension, per Rem.~\ref{rmk:l-btz}, for $0 \leq i \leq m_{\eta, p}^{\mathcal{C}_1}+(d-1)$, obtaining the BTZ extension vector $\Phi_{C, i}^{\ell b, \xi}$ and extension points $D^{\e, \xi}_{C}$.  Noting that
\begin{equation}\label{eq:C1rel1}
    D^{\e, \xi}_{C} \times \etavec{p}{\mathcal{C}_1}(0, m_{\eta, p}^{\mathcal{C}_1}) = \mathcal{H}_{p, 1}^{\mathcal{C}_1, 2} \cup \mathcal{H}_{p, 1}^{\mathcal{C}_1, 3}
\end{equation} 
where $\mathcal{H}_{p, 1}^{\mathcal{C}_1, 2}$ and $\mathcal{H}_{p, 1}^{\mathcal{C}_1, 3}$ (eqs.~\eqref{eq:H2pC1} and~\eqref{eq:H3pC1}) are two of the three component grids that define $\widetilde{H}_{p, 1}^{\mathcal{C}_1, \e}$ and, further, 
\begin{equation}\label{eq:C1rel2}
    D^{\e, \xi}_{C} \times \etavec{p}{\mathcal{C}_1}(m_{\eta, p}^{\mathcal{C}_1}, m_{\eta, p}^{\mathcal{C}_1}+d-1) \subset \widetilde{H}_p^{\mathcal{C}_1, \m},
\end{equation}
Step 1 now constructs the function $\widehat{g}_p^{\mathcal{C}_1, \bt}: \widetilde{H}_p^{\mathcal{C}_1, \m} \cup \widetilde{H}_p^{\mathcal{C}_1, \e} \to \mathbb{C}$ according to
\begin{equation}\label{eq:ghatC1}
    \widehat{g}_p^{\mathcal{C}_1, \bt}(\xi, \eta) =
    \begin{cases}
        \widetilde{f}_p^{\mathcal{C}_1}(\xi, \eta), & (\xi, \eta) \in \widetilde{\mathcal{A}}_{p}^{\mathcal{C}_1, \m}, \\
        \Phi_{C, i}^{\ell b, \xi}(\widetilde{x}_j^\xi), & (\xi, \eta) = (\widetilde{x}_j^\xi, \eta_i)\in D^{\e, \xi}_{C} \times \etavec{p}{\mathcal{C}_1}(0, m_{\eta, p}^{\mathcal{C}_1}+d-1), \\
        0, & \text{otherwise}.
    \end{cases}
\end{equation}
Clearly, the discrete function $\widehat{g}_p^{\mathcal{C}_1, \bt}$ coincides with the restriction to $\widetilde{H}_p^{\mathcal{C}_1, \m} \cup \widetilde{H}_{p, 1}^{\mathcal{C}_1, \e}$ of a smooth function defined in $\widetilde{\Omega}_p^{\mathcal{C}_1, \m} \cup \widetilde{\Omega}_p^{\mathcal{C}_1, \e}$, and the Step 1 construction is thus complete. 

Step 2 then simply defines the discrete function $g_p^{\mathcal{C}_1} :\widetilde{H}_p^{\mathcal{C}_1, \m}\rightarrow \mathbb{C}$ given by
\begin{equation}
    g_p^{\mathcal{C}_1}(\xi, \eta) = \widetilde{f}_p^{\mathcal{C}_1}(\xi, \eta) - \widehat{g}_p^{\mathcal{C}_1, \bt}(\xi, \eta), 
    \quad (\xi, \eta) \in \widetilde{H}_p^{\mathcal{C}_1, \m}.
\end{equation}
By definition, the grid function $g_p^{\mathcal{C}_1}$ vanishes on the set $\widetilde{\mathcal{A}}_{p}^{\mathcal{C}_1, \m}$. Moreover, since it is given by the difference of restrictions of smooth functions to the discrete set $\widetilde{H}_p^{\mathcal{C}_1, \m}$, it is itself the restriction of a smooth function to that set.

Step 3, in turn, blends to zero the discrete values of $g_p^{\mathcal{C}_1}$ on the discrete set $\widetilde{\mathcal{B}}_{p}^{\mathcal{C}_1,\m}$ along the $\eta$-direction. In detail, let
\[
D^{\m,\eta}=\etavec{p}{\mathcal{C}_1}(m_{\eta,p}^{\mathcal{C}_1},m_{\eta,p}^{\mathcal{C}_1}+d-1),
\]
and define the vector
\begin{equation}
    \Phi_i^{\ell,\eta}
    =
    \big(g_p^{\mathcal{C}_1}(\xi_i,x_0^\eta),\dots,g_p^{\mathcal{C}_1}(\xi_i,x_{d-1}^\eta)\big)^t,
    \qquad 0\le i\le m_{\xi,p}^{\mathcal{C}_1},
\end{equation}
whose entries are precisely the values of $g_p^{\mathcal{C}_1}$ at the points of $\widetilde{\mathcal{B}}_{p}^{\mathcal{C}_1,\m}$ with $\xi$-coordinate $\xi_i$. The extension vectors $\Phi_{C,i}^{\ell b,\eta}$ are then obtained at the points $D_C^{\e,\eta}$ as in Rem.~\ref{rmk:refinebtz}. Noting that
\begin{equation}\label{eq:C1rel3}
    \xivec{p}{\mathcal{C}_1}(0,m_{\xi,p}^{\mathcal{C}_1})\times D_C^{\e,\eta}
    =
    \mathcal{H}_{p,1}^{\mathcal{C}_1,1}\cup \mathcal{H}_{p,1}^{\mathcal{C}_1,3},
\end{equation}
Step~3 then produces the blended-to-zero discrete function
\begin{equation}
    g_p^{\mathcal{C}_1,\bt}(\xi,\eta)=
    \begin{cases}
        g_p^{\mathcal{C}_1}(\xi,\eta), & (\xi,\eta)\in \widetilde{\mathcal{B}}_{p}^{\mathcal{C}_1,\m},\\
        \Phi_{C_\rho,i}^{\ell b,\eta}(\widetilde{x}_j^\eta), & (\xi,\eta)=(\xi_i,\widetilde{x}_j^\eta)\in \xivec{p}{\mathcal{C}_1}(0,m_{\xi,p}^{\mathcal{C}_1})\times D_C^{\e,\eta},\\
        0, & \text{otherwise},
    \end{cases}
\end{equation}
which, in analogy with the preceding steps, equals the restriction of a smooth function to the discrete set $\widetilde{H}_p^{\mathcal{C}_1, \m} \cup \widetilde{H}_{p, 1}^{\mathcal{C}_1, \e}$---since for each $\eta \in D^{\m ,\eta}$ the function $g_p^{\mathcal{C}_1}(\xi_i, \eta)$, and, thus, the vectors $\Phi_{C_\rho, i}^{\ell b, \eta}$, smoothly approach the zero vector as the indices $i$ approach the index $m_{\xi, p}^{\mathcal{C}_1}$.

Step~4, finally, defines the final $\mathcal{C}_1$-type BTZ function $\widetilde{f}_p^{\mathcal{C}_1, \bt}: \widetilde{H}_{p, 1}^{\mathcal{C}_1, \e} \rightarrow \mathbb{C}$ by
\begin{equation}
    \widetilde{f}_p^{\mathcal{C}_1, \bt}(\xi, \eta) 
    = g_p^{\mathcal{C}_1, \bt}(\xi, \eta) + \widehat{g}_p^{\mathcal{C}_1, \bt}(\xi, \eta), 
    \quad (\xi, \eta) \in \widetilde{H}_{p, 1}^{\mathcal{C}_1, \e}.
\end{equation}
In analogy with the previous steps, $g_p^{\mathcal{C}_1, \bt} + \widehat{g}_p^{\mathcal{C}_1, \bt}$ is the sum of the restrictions of smooth functions to the discrete set $\widetilde{H}_p^{\mathcal{C}_1, \m} \cup \widetilde{H}_{p, 1}^{\mathcal{C}_1, \e}$, and it is therefore itself the restriction of a smooth function to that set. Since, by definition, the grid function $g_p^{\mathcal{C}_1, \bt} + \widehat{g}_p^{\mathcal{C}_1, \bt}$ coincides with the original values of $\widetilde{f}_p^{\mathcal{C}_1}$ on the matching mesh $\widetilde{H}_p^{\mathcal{C}_1, \m}$, the desired blending-to-zero procedure is complete.
\begin{figure}[H]
\centering
\includegraphics[width=\linewidth,]{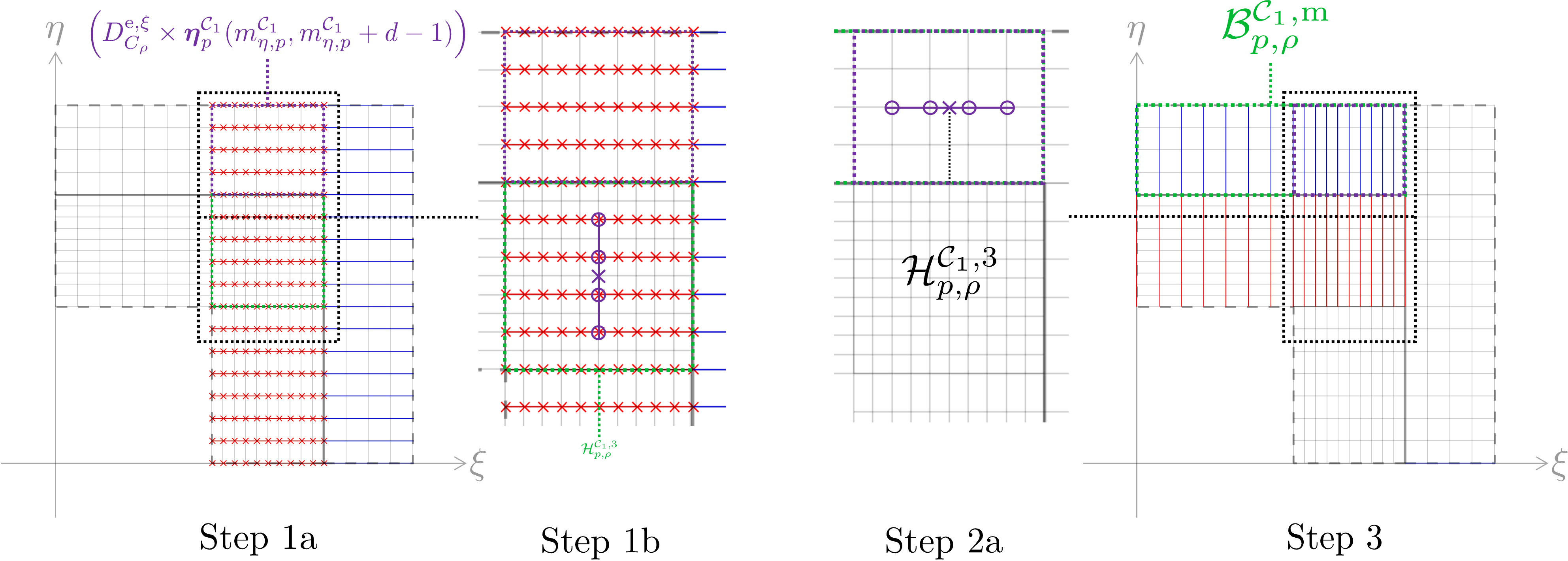}
    \caption{Important elements in the six-step $\rho$-refined $\mathcal{C}_1$-type BTZ procedure described in Rem.~\ref{rmk:C1nrho}.}
    \label{fig:C1BTZrho}
\end{figure}

\begin{remark}\label{rmk:C1nrho}
As noted in Rem.~\ref{rmk:refinebtz},  to avoid loss of accuracy in the interpolation steps, particularly when high precision is sought, for each patch type the BTZ procedure is performed over a $\rho$-refined extension mesh (equal to~\eqref{BTZ_patch_mesh_S} for $\mathcal{S}$-type patches and equal to~\eqref{BTZ_patch_mesh_C} for $\mathcal{C}_1$- and $\mathcal{C}_2$-type patches) within the relevant ($\rho$-independent) extension region~\eqref{eq:omegatildepSe}, \eqref{eq:omegatildepC1e}, or~\eqref{eq:omegatildepC2e}, using an integer refinement factor $\rho > 1$. For clarity, in Steps 1-4 above, the $\mathcal{C}_1$-patch BTZ extension procedure was presented only for the simplest,  unrefined-mesh case, wherein the refinement parameter $\rho$ is set to $\rho = 1$. The corresponding $\mathcal{C}_1$-patch BTZ extension method with $\rho > 1$ proceeds on the basis of a 2D interpolation approach based on the refined 1D-BTZ procedure (Rem.~\ref{rmk:refinebtz}), as described below. This approach substantially mitigates the compounded interpolation error that arises from the successive one-dimensional interpolation steps inherent in the 2D-BTZ algorithm for $\mathcal{C}_1$-type patches. More precisely, the discrete extension $\widetilde{f}_p^{\mathcal{C}_1, \bt}$ on the $\rho$-refined mesh $\widetilde{H}_{p, \rho}^{\mathcal{C}_1, \e}$ ($\rho > 1$, see Fig.~\ref{fig:OmegapRe}), is obtained through a six-step modified version, detailed in what follows, of the four-step procedure introduced earlier in this section.
\paragraph{Six-step $\rho$-refined $\mathcal{C}_1$-type BTZ procedure:}
\begin{enumerate}
    \item [1a.] Blend $\widetilde{f}_p^{\mathcal{C}_1}$ along the $\xi$-direction of $\widetilde{\mathcal{A}}_p^{\mathcal{C}_1, \m}$ using the refined 1D-BTZ procedure (Rem.~\ref{rmk:refinebtz}) to obtain extension values on the $\rho$-refined Step-1a extension mesh $D^{\e, \xi}_{C_\rho} \times \etavec{p}{\mathcal{C}_1}(0, m_{\eta, p}^{\mathcal{C}_1}+d-1)$ and let 
    \begin{equation}
        \widetilde{H}_{p, \rho}^{\mathcal{C}_1, \m} = \widetilde{H}_{p}^{\mathcal{C}_1, \m} \cup \left( D^{\e, \xi}_{C_\rho} \times \etavec{p}{\mathcal{C}_1}(m_{\eta, p}^{\mathcal{C}_1}, m_{\eta, p}^{\mathcal{C}_1}+d-1)\right)
    \end{equation}
    denote the union of the matching mesh $\widetilde{H}_p^{\mathcal{C}_1, \m}$ with the portion of the Step-1a extension mesh that intersects with the matching domain $\widetilde{\Omega}_p^{\mathcal{C}_1, \m}$ (Step 1a panel of Fig.~\ref{fig:C1BTZrho}).

    \item[1b.] Define the blended-to-zero function $\widehat{g}_p^{\mathcal{C}_1, \bt}: \widetilde{H}_{p, \rho}^{\mathcal{C}_1, \m} \cup \widetilde{H}_{p}^{\mathcal{C}_1, \e} \rightarrow \mathbb{C}$ so that it coincides with the original function values on $\widetilde{\mathcal{A}}_p^{\mathcal{C}_1, \m}$; it coincides with the BTZ extension values on the $\rho$-refined Step-1a extension mesh; it interpolates the BTZ extension values from the $\rho$-refined Step-1a extension mesh onto $\mathcal{H}_{p, \rho}^{\mathcal{C}_1, 3} \subset \widetilde{H}_{p, \rho}^{\mathcal{C}_1, \e}$ (where, as shown in the Step 1b panel of Fig.~\ref{fig:C1BTZrho}, the grid $\mathcal{H}_{p, \rho}^{\mathcal{C}_1, 3}$ is more finely refined in the $\eta$-direction, thereby requiring interpolation); and vanishes everywhere else.

    \item[2a.] Interpolate the original function values $\widetilde{f}_p^\mathcal{R}$ from the coarse matching mesh $\widetilde{H}_{p}^{\mathcal{C}_1, \m}$ onto the target mesh $\widetilde{H}_{p, \rho}^{\mathcal{C}_1, \m}$ as shown in the Step 2a panel of Fig.~\ref{fig:C1BTZrho}, which is refined in the $\xi$-direction in a subset of the matching domain (outlined in purple in the Step 2a panel of Fig.~\ref{fig:C1BTZrho}).

    \item[2b.] Subtract $\widehat{g}_p^{\mathcal{C}_1, \bt}$, obtained in Step 1, from the interpolated function values from Step 2a on $\widetilde{H}_{p, \rho}^{\mathcal{C}_1, \m}$ to obtain the discrete function $g_p^{\mathcal{C}_1}: \widetilde{H}_{p, \rho}^{\mathcal{C}_1, \m} \rightarrow \mathbb{C}$.

    \item[3.] Blend to zero the grid function $g_p^{\mathcal{C}_1}$ along the $\eta$-direction of $\widetilde{\mathcal{B}}_{p, \rho}^{\mathcal{C}_1, \m}$, where
    \begin{equation}
        \widetilde{\mathcal{B}}_{p, \rho}^{\mathcal{C}_1, \m} = \widetilde{\mathcal{B}}_{p}^{\mathcal{C}_1, \m} \cup \left( D^{\e, \xi}_{C_\rho} \times \etavec{p}{\mathcal{C}_1}(m_{\eta, p}^{\mathcal{C}_1}, m_{\eta, p}^{\mathcal{C}_1}+d-1)\right),
    \end{equation}
    to obtain blending-to-zero values on the $\rho$-refined extension mesh 
    \begin{equation}\label{eq:C1BTZrhos3e}
        \left(\xivec{p}{\mathcal{C}_1}(0, m_{\xi, p}^{\mathcal{C}_1} -C) \cup D_{C_\rho}^{\e, \xi}\right) \times D_{C_\rho}^{\e, \eta},
    \end{equation}
    as illustrated in the Step 3 panel of Fig.~\ref{fig:C1BTZrho}, and define the blended-to-zero function $g_p^{\mathcal{C}_1, \bt}: \widetilde{H}_{p, \rho}^{\mathcal{C}_1, \m} \cup \widetilde{H}_{p, \rho}^{\mathcal{C}_1, \e} \rightarrow \mathbb{C}$ to coincide with the values of $g_p^{\mathcal{C}_1}$ on $\widetilde{H}_{p, \rho}^{\mathcal{C}_1, \m}$; to coincide with the BTZ extension values, just obtained, on the mesh~\eqref{eq:C1BTZrhos3e}; and to equal zero everywhere else in $\widetilde{H}_{p, \rho}^{\mathcal{C}_1, \m} \cup \widetilde{H}_{p, \rho}^{\mathcal{C}_1, \e}$.

    \item[4] Define the final $\mathcal{C}_1$-type BTZ function $\widetilde{f}_p^{\mathcal{C}_1, \bt}$ as the sum of $\widehat{g}_p^{\mathcal{C}_1, \bt}$ and ${g}_p^{\mathcal{C}_1, \bt}$ on $\widetilde{H}_{p, \rho}^{\mathcal{C}_1, \e}$.
\end{enumerate}
\end{remark}

\begin{figure}
    \centering
    \includegraphics[width=0.9\linewidth]{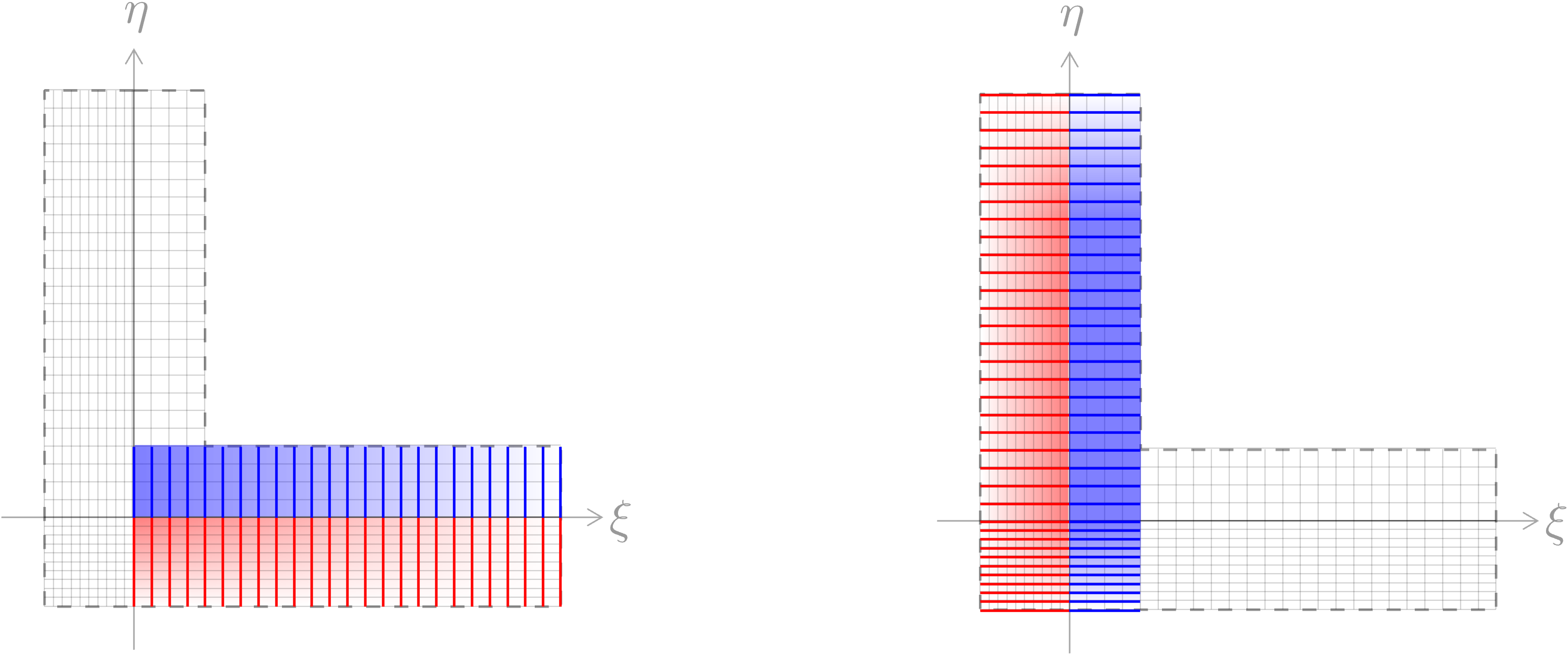}
    \caption{
    Parameter-space BTZ construction for $\mathcal{R}=\mathcal{C}_2$; cf. Fig.~\ref{fig:C2BTZprelim}. The blue and red lines represent the parameter-space matching and extension meshes $\widetilde{H}_p^{\mathcal{C}_2, \m}$ and $\widetilde{H}_{p, \rho}^{\mathcal{C}_2, \e}$, respectively. Color shading is used to indicate the smooth blending-to-zero induced by the windowing functions $w_p^{\mathcal{C}_2}$ along the relevant boundaries of the extension and matching domains $\Omega_p^{\mathcal{C}_2,\e}$ and $\Omega_p^{\mathcal{C}_2,\m}$.}
    \label{fig:C2BTZdetail}
\end{figure}
\subsection{$\mathcal{C}_2$-type patch blending to zero in parameter space}\label{subsec:C2patchBTZ}
Finally, in the $\mathcal{C}_2$-type patch case the BTZ algorithm proceeds as follows. Given the windowed grid function $\widetilde{f}_p^{\mathcal{C}_2}$ defined on the patchwise matching mesh $\widetilde{H}_p^{\mathcal{C}_2, \m} \subset \widetilde{\Omega}_p^{\mathcal{C}_2}$ of a given $\mathcal{C}_2$-type patch $\Omega_p^{\mathcal{C}_2}$, the corresponding blended-to-zero grid function $\widetilde{f}_p^{\mathcal{C}_2, \textrm{b}}$ is first obtained on the component grid $\mathcal{H}_{p, \rho}^{\mathcal{C}_2, 1}$ (eq.~\eqref{eq:H1pC2}) in the same manner as for the $\mathcal{S}$-type patch. Specifically, with reference to Rem.~\ref{rmk:l-btz}, using the vector 
\begin{equation}
    \Phi^{\ell, \eta} = \Phi^{\ell, \eta}_i = (f^{\mathcal{C}_2}_p(\xi_i, x_0^\eta), \dots, f^{\mathcal{C}_2}_p(\xi_i, x_{d-1}^\eta) )^t
\end{equation}
of windowed grid function values over the 1D matching set $D^{\m, \eta} = \etavec{p}{\mathcal{C}_2}(0, d-1)$ together with the corresponding blended-to-zero grid function $\Phi_{C_{\rho}, i}^{\ell b, \eta}$ (obtained per Rem.~\ref{rmk:refinebtz}) over  the extension  grid $D^{\e, \eta}_{C_\rho}$, and noting that
\begin{equation}
    \mathcal{H}_{p, \rho}^{\mathcal{C}_2, 1} = \{\xi_0, \xi_1, \dots, \xi_{n_{\xi, p}^{\mathcal{C}_2}-1}\} \times \{\widetilde{x}_{-\rho C}^\eta, \widetilde{x}_{-\rho C + 1}^\eta, \dots, \widetilde{x}_0^\eta\},
\end{equation}
we define the $\mathcal{C}_2$ blending to zero function
\begin{equation}
    \widetilde{f}_p^{\mathcal{C}_2, \bt}(\xi_i, \widetilde{x}_j^\eta) = \Phi^{\ell b, \eta}_{C_\rho, i}(\widetilde{x}_j^\eta) \quad \textrm{for} \quad (\xi_i, \widetilde{x}_j^{\eta})\in \mathcal{H}_{p, \rho}^{\mathcal{C}_2, 1}.
\end{equation}
The blended-to-zero grid function values over the component grids $\mathcal{H}_{p, \rho}^{\mathcal{C}_2, 2}$ and $\mathcal{H}_{p, \rho}^{\mathcal{C}_2, 3}$ (eq.~\eqref{eq:H2pC2} and~\eqref{eq:H3pC2}) are then obtained by applying 1D $\ell$-BTZ blending along the $\xi$ direction of the appropriate mesh: $H_p^{\mathcal{C}_2, \m}$ for $\mathcal{H}_{p, \rho}^{\mathcal{C}_2, 2}$ and $\mathcal{H}_{p, \rho}^{\mathcal{C}_2, 1}$ for $\mathcal{H}_{p, \rho}^{\mathcal{C}_2, 3}$. In detail, letting $D^{\m, \xi} = \xivec{p}{\mathcal{C}_2}(0, d-1)$,  for each $\widetilde{\eta} \in \etavec{p}{\mathcal{C}_2}(0, n_{\xi, p}^{\mathcal{C}_2}-1) \cup D^{\e, \eta}_{C_\rho}$ we consider the vector
\begin{equation}
    \Phi^{\ell, \xi} = \Phi^{\ell, \xi}_{\widetilde{\eta}} = (\widetilde{f}_p^{\mathcal{C}_2}(x_0^\xi, \widetilde{\eta}), \dots, \widetilde{f}_p^{\mathcal{C}_2}(x^\xi_{d-1}, \widetilde{\eta}))^t
\end{equation}
together with the corresponding blending-to-zero grid function $\Phi_{C_\rho, \eta}^{\ell \bt}$ over the grid $D_{C_\rho}^{\e, \xi}$. Then, noting that
\begin{equation}
    \mathcal{H}_{p, \rho}^{\mathcal{C}_2, 2}\cup \mathcal{H}_{p, \rho}^{\mathcal{C}_2, 3} = D_{C_\rho}^{\e, \xi} \times \left(\etavec{p}{\mathcal{C}_2}(0, n_{\xi, p}^{\mathcal{C}_2}-1) \cup D^{\e, \eta}_{C_\rho}\right),
\end{equation}
we define 
\begin{equation}
    \widetilde{f}_p^{\mathcal{C}_2, \bt}(\widetilde{x}_j^\xi, \eta) = \Phi_{C_\rho, \widetilde{\eta}}^{\ell \bt}(\widetilde{x}_j^\xi) \quad \textrm{for} \quad (\widetilde{x}_j^\xi, \widetilde{\eta}) \in \mathcal{H}_{p, \rho}^{\mathcal{C}_2, 2}\cup \mathcal{H}_{p, \rho}^{\mathcal{C}_2, 3}
\end{equation}
---which provides the desired $\mathcal{C}_2$-type extension of the given grid function $\widetilde{f}_p^{\mathcal{C}_2}$.

\section{Interpolation onto the Cartesian mesh}\label{SM4}

\subsection{Step 1: Identifying $ H \cap \Omega_p^{\mathcal{R}, e} $}\label{subsec:id-int}
This step in the interpolation algorithm consists of constructing, for each $(p, \mathcal{R})$, the set of patchwise interpolation target points $H \cap \Omega_p^{\mathcal{R}, \e}$---that is, the set of global Cartesian grid points that require interpolation from the extension patch mesh $H_p^{\mathcal{R}, \e}$ defined in eq.~\eqref{eq:HpRe}. 

To carry out this step, the boundary of $\Omega_p^{\mathcal{R}, \e}$ is first approximated by a sufficiently fine polygonal curve. An efficient two-stage algorithm---described below---is then used to generate a Boolean mask over $H$, which is set to \texttt{true} precisely at those points in $H$ that lie within the resulting polygonal domain. The initial mask is set to false at all points in $H$, and the algorithm begins by obtaining all the intersections between the Cartesian discretization lines and the sides of the polygonal approximation--- and then, flagging grid points that lie immediately to the left of the intersection points. The algorithm then performs a row-wise scan across the Cartesian mesh and applies the even–odd rule to mark as \texttt{true} all the points in the interior. (The even-odd rule~\cite{galetzka2017simple} determines if a point lies inside or outside a polygon by casting a ray from the point in any direction and counting edge intersections.)

To describe the algorithm in detail, let
\begin{align*}
    E &= \{e_1, e_2, \dots, e_m\}\\
    &= \left\{\{(x_1, y_1), (x_2, y_2)\}, \{(x_2, y_2), (x_3, y_3)\}, \dots, \{(x_m, y_m), (x_{m+1}, y_{m+1})\}\right\}
\end{align*}
denote the set of $m$ polygon edges that approximate the boundary of the relevant domain (with $(x_{m+1}, y_{m+1}) = (x_1, y_1)$). For each $k=1, \dots, m$, the algorithm computes (using notation from~\eqref{eq:mesh})
\begin{equation}
    \widetilde{y}_k = \frac{y_k-b_1}{h}.
\end{equation}
Treating $H$ as a collection of horizontal lines
\begin{equation}
    \{y = jh+b_1:0\leq j\leq n_y-1\},
\end{equation}
and, assuming, without loss of generality that, for a given $k$ value we have $y_{k+1} \geq y_{k}$, each $e_k$ intersects the lines
$y = jh+b_1$ for $\lfloor \widetilde{y}_k\rfloor < j\leq \lfloor\widetilde{y}_{k+1}\rfloor$ at points of the form $(\widehat{x}_{j, k}, jh+b_1)$ where
\begin{equation}
    \widehat{x}_{j, k} = x_k+\frac{(jh+b_1)-y_{k}}{y_{k+1}-y_k}(x_{k+1}-x_{k}).
\end{equation}
The Cartesian grid point directly to the left of $(\widehat{x}_{j, k}, jh+b_1)$ is
\begin{equation}\label{eq:bound_flag}
    \left(\left\lfloor \frac{\widehat{x}_{j, k}-a_1}{h} \right\rfloor h+a_1, jh+b_1\right),
\end{equation}
and the corresponding boolean mask value is marked as desired: if the value is originally \texttt{false}, it is set to \texttt{true}; if it is originally \texttt{true}, it is set to \texttt{false}. Note that~\eqref{eq:Rsup} ensures that~\eqref{eq:bound_flag} is indeed contained in $H$.

With the grid points lying immediately to the left of boundary segments now flagged, the remainder of the boolean mask can be completed using the aforementioned even--odd rule. Specifically, for each horizontal row of the Cartesian mesh, we initialize a logical variable \texttt{interior} as \texttt{false}. We then iterate from left to right across the row and proceed according to the following four cases:
\begin{itemize}
    \item[(i)] If \texttt{interior} is \texttt{false} and the current mask value is \texttt{false}, no action is taken.
    \item[(ii)] If \texttt{interior} is \texttt{false} and the current mask value is \texttt{true}, the mask value is reset to \texttt{false} (as it lies to the left of the boundary, which must be the exterior), and \texttt{interior} is set to \texttt{true}.
    \item[(iii)] If \texttt{interior} is \texttt{true} and the current mask value is \texttt{false}, the mask value is set to \texttt{true} (as it lies in the interior).
    \item[(iv)] If \texttt{interior} is \texttt{true} and the current mask value is \texttt{true}, the mask value is unchanged, and \texttt{interior} is reset to \texttt{false} (indicating the crossing of another boundary).
\end{itemize}
The completed mask then directly identifies the set of interpolation target points $H \cap \Omega_p^{\mathcal{R}, \e}$. This method was found to be orders of magnitude faster than MATLAB's \texttt{inpolygon} function, a performance gain that proved essential for the overall efficiency of the 2D-FC algorithm.

Since this algorithm relies on a polygonal approximation of the boundary $\Gamma$, it may misidentify points located very close to the boundary of $\Omega_p^{\mathcal{R}, \e}$. In practice, this is only relevant for points near $\Omega_p^{\mathcal{R}, \e} \cap \Gamma$, as $f_p^{\mathcal{R}, \bt}$ vanishes at the other patch boundary segments. To address this, the algorithm can first slightly extrude $\Omega_p^{\mathcal{R}, \e}$ toward the interior of $\Omega$ and construct a polygonal approximation of the enlarged domain. Then, in the step described in Sec.~\ref{subsec:computing-inv-im}, where the corresponding parameter-space coordinates are computed, the algorithm precisely re-checks whether the points lie within the parameter-space patchwise extension domain $\widetilde{\Omega}_p^{\mathcal{R}, \e}$, discarding any points as appropriate.

\subsection{Step 2: Computing $ \left(\mathcal{M}_p^{\mathcal{R}}\right)^{-1}({H}_p^{\mathcal{R}}) $}\label{subsec:computing-inv-im}

To interpolate at each point $ \mathbf{x} \in H \cap \Omega_p^{\mathcal{R}, e}$ in parameter space, it is necessary to compute its corresponding coordinate $ (\xi, \eta) \in \widetilde{\Omega}_p^{\mathcal{R}, \e} $ via inversion of the parameterization map: $\mathcal{M}_p^{\mathcal{R}}(\xi, \eta) = \mathbf{x}$.
This inverse mapping $ \left(\mathcal{M}_p^{\mathcal{R}}\right)^{-1}({H}_p^{\mathcal{R}}) $ is computed numerically using Newton’s method. Accurate and efficient inversion requires careful selection of initial guesses to ensure rapid convergence. To do this we adopt the following strategy:
\begin{itemize}
    \item[(i)] First, construct $\widetilde{H}_p^{\mathcal{R}, \e}$ and map these points forward via $ \mathcal{M}_p^{\mathcal{R}} $ to obtain $ H_p^{\mathcal{R}, e} $.
    \item[(ii)] For each $ \mathbf{y} \in {H}_p^{\mathcal{R}, \e} $, identify nearby Cartesian grid points in $H \cap \Omega_p^{\mathcal{R}, \e}$. Specifically, determine the four nearest Cartesian neighbors (using floor and ceiling operators on $ \mathbf{y} $) and check whether any of these points have known inverse mappings.
    \item[(iii)] If one or more nearby points have unknown inverse coordinates $ (\xi, \eta) $, use $\mathbf{y}$ (which has a known inverse mapping) as the initial guess for Newton’s method.
    \item[(iv)] Apply Newton’s method to refine the inverse to a prescribed tolerance.
    \item[(v)] Iterate through all points in $H \cap \Omega_p^{\mathcal{R}, \e}$ using some (arbitrary) ordering of this set. For each point lacking a valid inverse mapping, examine its eight surrounding Cartesian neighbors. If at least one neighbor possesses a valid computed $(\xi, \eta)$ value, use that value as the initial guess for the inverse map and apply Newton's method to compute inverse; otherwise skip the current point and continue to the next point in the order used.
    \item[(vi)] Repeat the above procedure iteratively starting from~(v) until a valid inverse mapping has been assigned to all points.

\end{itemize}
This localized reuse of known inverses as initial guesses significantly reduces the number of Newton iterations required and improves the efficiency of the mapping inversion.

\section{Evaluation of the homogeneous solution for the Poisson solver}\label{sec:poisson-eval}
With reference to eq.~\eqref{eq:rescaled-ie}, the solution of the problem~\eqref{eq:def_laplace} is given by
\begin{equation}\label{eq:poisson_sol}
    v(\mathbf{x}) = \sum_{j}\int_0^1 
    K\bigl(\mathbf{x}, \mathbf{r}_j(w(s'))\bigr)
    \widehat{\varphi}\bigl(\mathbf{r}_j(w(s'))\bigr)\,
    J_j(w(s'))\,ds', \quad \mathbf{x} \in \Omega.
\end{equation}
The function $v(\mathbf{x})$ may, in principle, be evaluated at arbitrary target points \(\mathbf{x}\in\Omega\) by direct application of the trapezoidal rule once discrete values of the boundary density are available. As the target point approaches the boundary \(\Gamma\), however, the integrand in~\eqref{eq:poisson_sol} develops increasingly sharp near-singular behavior, which leads to a pronounced loss of accuracy when a fixed boundary quadrature mesh is employed. To address this difficulty, we adapt the refined-window integration (RWI) and iterated refined-window integration (IRWI) algorithms introduced in~\cite{bruno_two-dimensional_2022}, with modifications appropriate to the present setting involving geometries containing corner points. In contrast to~\cite{bruno_two-dimensional_2022}, where the RWI and IRWI methods are applied directly to the boundary density (which is not smooth at corner points), here we instead apply these methods to the rescaled density $\widehat{\varphi}$, whose discrete values are obtained by discretizing the integral eq.~\eqref{eq:rescaled-ie} by employing the Nystr\"om method described in~\cite[Sec.~12.2]{book-kress-lie-2014}. Once computed on the boundary discretization, $\widehat{\varphi}$ is represented globally along $\Gamma$ by means of a truncated Fourier expansion, which can be refined efficiently and with high accuracy via the FFT. This representation enables trapezoidal-rule integration on arbitrarily refined boundary meshes, thus controlling the integration error. Specifically, given a prescribed integration tolerance $\varepsilon_{\mathrm{int}}$, the IRWI algorithm determines, for each target point $\mathbf{x}$, the minimal refinement factor $R$ needed to resolve the near-singular behavior of the kernel in~\eqref{eq:poisson_sol}. As in~\cite{bruno_two-dimensional_2022}, the required refinement is estimated by computing the difference between successive approximations obtained with refinement factors $R$ and $R+1$; refinement is increased until this difference falls below $\varepsilon_{\mathrm{int}}$. This adaptive procedure ensures that the integration accuracy is controlled uniformly across the computational domain while avoiding unnecessary refinement away from the boundary. 

Utilizing the near-boundary operator evaluation method just described, we now outline the method used for evaluation of the solution on the Cartesian grid $H$ (see~\eqref{eq:mesh}) (note, however, that the procedure described below can in principle be applied to arbitrary sets of target points within $\Omega$). Noting that, since many points in $H$ may lie arbitrarily close to the boundary, direct application of the IRWI algorithm can be computationally expensive. To mitigate this cost, we follow the approach proposed in the reference~\cite{bruno_two-dimensional_2022} and employ near-boundary interpolation of degree $M-1$ (see Sec.~\ref{subsec:interpolationontomesh}) whenever possible.
In detail, for target points lying within a distance $(M-1)h$ of a smooth portion of the boundary and at least $(M-1)h$ away from any corner point, IRWI is not applied directly. Instead, the solution is first evaluated---using IRWI with tolerance $\varepsilon_{\mathrm{int}}$---at a set of auxiliary interpolation nodes located on a tangential--normal grid adjacent to the boundary. The value of $v$ at the Cartesian grid point is then obtained by iterated 1D polynomial interpolation of degree $M-1$ following the procedure outlined in Sec.~\ref{subsec:interpolationontomesh}. This near-boundary interpolation strategy substantially reduces the number of expensive IRWI evaluations while preserving high-order accuracy in regions where the boundary geometry remains smooth.

For points lying within $(M-1)h$ of a corner point, the reduced regularity induced by the corner geometry precludes the use of near-boundary interpolation. In these regions, the solution is therefore evaluated directly by the IRWI algorithm, using the same tolerance-driven refinement strategy described above. Since the number of such corner-adjacent points is small in practice, the additional computational cost incurred by direct IRWI evaluation remains limited. Finally, for points that are sufficiently far from both the boundary and the corners, the solution is evaluated by direct application of the IRWI algorithm, typically requiring little or no refinement.

In summary, the proposed evaluation strategy combines adaptive IRWI integration with near-boundary interpolation to efficiently and accurately compute the homogeneous solution $v(\mathbf{x})$ throughout the domain. Near smooth portions of the boundary, interpolation-based evaluation significantly reduces the cost associated with repeated near-singular integrations, while direct IRWI evaluation is reserved for the limited number of points near corners, where reduced regularity precludes interpolation. Away from the boundary, the algorithm typically requires little or no refinement. Together, these components yield a robust and accurate scheme whose effectiveness is demonstrated in Exs.~\ref{ex:guitarbase-poison} and~\ref{ex:osci-poisson}.

\bibliographystyle{siamplain}
\bibliography{references}